\documentclass[a4paper,12pt]{article}

\usepackage[margin=0.5in]{geometry}
\usepackage{amsmath}
\usepackage{amssymb}
\usepackage{multicol}
\usepackage{bigints}
\usepackage{tikz-network}
\usepackage{tcolorbox}
\usepackage{setspace}
\usepackage{enumitem}
\usepackage{graphicx}
\usepackage{mathtools}
\usepackage{physics}
\usepackage{amsthm}
\usepackage{amsmath,amssymb,xcolor}
\usepackage{graphics}
\usepackage{adjustbox}
\usepackage{txfonts}
\usepackage{multirow}
\usepackage[square,numbers]{natbib}
\title{Rehan-Lanel Indices of Graphs}
\usepackage{authblk}
\author[1,2]{D.C. Gunawardhana}
\author[2]{ G.H.J. Lanel}
\affil[1]{Department of Mathematical Sciences, Faculty of Applied Sciences, South Eastern University of Sri Lanka}
\affil[2]{Department of Mathematics, Faculty of Applied Sciences, University of Sri Jayewardenapura}

\usepackage{longtable}
\usepackage[english]{babel}
\newtheorem{theorem}{Theorem}[]
\newtheorem{proposition}{Proposition}
\newtheorem{corollary}{Corollary}[proposition]

\newtheorem{definition}{Definition}

\providecommand{\keywords}[1]
{
	\small	
	\textbf{\textit{Keywords---}} #1
}

\begin{document}
	
	\maketitle
	\begin{abstract}
		A graph \(G\) consists of vertices \(V(G)\) and edges \(E(G)\). 
		In this paper, we propose four new indices defined and named  as first Rehan-Lanel index of \(G\) \((RL_1)\), second Rehan-Lanel index of \(G\) \((RL_2)\), second Rehan-Lanel index of \(G\) \((RL_3)\) and fourth Rehan-Lanel index of \(G\) \((RL_4)\) . The degrees of the vertices \(u, v \in V(G)\) are denoted by  \(d_G(u)\) and \(d_G(v)\).  Based on these new indices and the definitions of Revan degree, Domination degree, Banhatti degree, Temperature of a vertex, KV indices, we subsequently  introduced an additional 448 indices/exponentials and computed results for first four new indices of each subsequent definitions, for the standard graphs such as \(r-\) regular graph, complete graph, cycle, path and compete bipartite graph.In addition, we performed calculations for the Wheel graph, Sunflower graph, and French Windmill graph. Furthermore, using exponential of a degree of a vertex, the centrality concept, we introduced another 8 indices. Furthermore, we defined new degree called Chandana-Lanel degree of a vertex of a graph(CL degree). Using this degree, new 6 indices were defined. Also, we defined the index called the Heronian Rehan-Lanel index using the heronian mean of two number. These novel 462 indices would be advantageous in QSPR/QSAR studies.
	\end{abstract}
	
	\footnote{* Corresponding author.\\
		Email addresses: rehan@seu.ac.lk(D.C. Gunawardhana), ghjlanel@sjp.ac.lk(Dr. G.H.J. Lanel) }
	
	\section{Introduction}
	A topological index is a numerical quantity that describes the topology of a graph.
	Topological indices of molecular graphs play a well-known vital role in Chemical Graph Theory, and such indices can be used to study the chemical characteristics of pharmaceuticals. A topological index numerically describes the molecule form and is used in the advancement of qualitative structure-activity relationships (QSARs). Several degree-based graph indices have been published in the literature \cite{article71} and have found some applications in \cite{article69, article70} .\\
	
	There are 3 types of topological indices as follows \cite{article72}, namely degree based, distance based and spectral based indices.

	Degree-based indices are notably popular among the various types of indices. Topological indices are very useful in chemical validation, QSAR/QSPR studies, and pharmaceutical drug formulation. Degree-based indices have been thoroughly examined and are thought to be associated with different features of the studied chemical compounds. This study will further explore degree-based indicators.
	
	While researching paraffins, an alkene with a high boiling point, Harold Wiener established the concept of topological indices \cite{article73}. The Wiener index, the first topological index, was introduced by him in 1947. He designated it as path number. As research in chemical graph theory progressed, the path number was given the name Wiener index. Because of its unique theoretical properties and broad variety of applications, the Wiener index is the most explored molecular topological index in chemical graph theory \cite{article72, article74}.\\
	
	\textbf{Milan Randic} \cite{article75} proposed the first realistic degree-based topological index in his influential study published in 1975. Randic termed it \textit{branching index} at beginning. Later it was renamed connectivity index. Nowadays it is popular as the Randic index.
	In 1972, the total \(\pi\)-electron energy of a molecular graph was first demonstrated to depend on its structure \cite{a2} . It can be calculated by adding the squares of the degrees at each vertex. Later, this sum was dubbed the first Zagreb index, and substantial research was conducted on it. \textbf{I. Gutman and N. Trinanjtic} introduced the first and second Zagreb indices \cite{article76}. Another topological index, defined as the total of the cubed degrees of the vertices of a network, was also demonstrated to relate \(\pi\)- electron energy in the same publication. This index, however, was never thoroughly researched at the time.
	
	Furtula and Gutman \cite{article77} restudied this index in 2015 to establish some fundamental features and named it the ''Forgotten topological index" or "F-index.
	
	In \cite{article79},	Vukicevic \& Furtula proposed geometrical-arithmetic index in 2009. Similarly more than 3000 topological indices are available in the literature.  In this study, we hope to introduce few more indices to chemical graph theory.
	
	
	\section{Methodology}
	The first and second Zagreb indices of a molecular graph G are defined as, 
	
	\[M_1(G)=\sum_{uv\in E(G)}d_G(u)+d_G(v)\]
	
	and 
	\[M_2(G)=\sum_{uv\in E(G)}d_G(u)d_G(v)\]
	
	Based on the definitions of the Zagreb indices and their wide applications,
	Kulli \cite{article81} introduced the first Gourava index of a molecular graph as follows:
	
	The first Gourava index of a graph G is defined as,
	\[GO_1(G)=\sum_{uv\in E(G)}d_G(u)+d_G(v)+d_G(u)d_G(v)\]
	
	Using the first and second Zagreb indices, Furtula
	and Gutman \cite{article82} introduced forgotten topological index
	(also called F-index) which was defined as
	\[F(G)=\sum_{u\in V(G)}d^3_G(u)=\sum_{uv\in E(G)}d^2_G(u)+d^2_G(v).\]

	Inspired by these definitions, here we hope to introduce two new indices by combining the remaining indices.
	
	In addition, It is expected to define  the new topological indices based on the definitions of Revan degree, Banhatti degree, temperature, domination degree, neighborhood of a vertex of a graph,  KV indices and Multiplicative indices. , 

	\section{Results}
	In this paper, we hope to introduce  two fresh indices and based on those two indices and some other definitions of a graph such as E-Banhatti degree, Revan degree, Domination degree, Temperature, KV indices and Neighborhood degree sum, we present another new indices/exponentials.\\

	\subsection{Rehan-Lanel indices}

		We can write \(d_G(u)-d_G(v)+d_G(u)d_G(v)=d_G(u)+d_G(v)(d_G(u)-1)>0\)	for all \(d_G(u) \ge 1.\)
	\begin{definition}
		Let	a graph \(G\) consists of a set of vertices $V(G)$ and a set of edges $E(G)$. We propose four new indices called 'the first Rehan-Lanel index' and 'the second Rehan-Lanel index', 'the third Rehan-Lanel index' and 'the fourth Rehan-Lanel index' and they are defined as,
		\begin{center}
			\begin{enumerate}[label=\roman*.]
				\item \(RL_1(G)=\sum\limits_{uv\in E(G)}\left(  d^2_G(u)+d^2_G(v)+d_G(u)d_G(v)\right)\)
				\item \(RL_2(G)=\sum\limits_{uv\in E(G)}\left(  d^2_G(u)+d^2_G(v)-d_G(u)d_G(v)\right).\)
				 \item \(RL_3(G)=\sum\limits_{uv\in E(G)}\left(  d_G(u)-d_G(v)+d_G(u)d_G(v)\right)\)
				\item \(RL_4(G)=\sum\limits_{uv\in E(G)}\left(  |d_G(u)-d_G(v)|d_G(u)d_G(v)\right).\)
			\end{enumerate}
		\end{center}

	\end{definition}

	Using these four  indices and the definitions of multiplicative indices, Revan degree, Banhatti degree, , Temperature and neighborhoods of vertices, we propose more new indices and relevant exponential in the Table \ref{table1a},  \ref{table1b}, \ref{table1c} and  \ref{table1d}.
	
	\begin{longtable}[c]{| c | l | l |}
		\caption{Rehan- Lanel indices/exponentials\label{table1a}}\\
		
		\hline
		\multicolumn{1}{|c}{\textbf{S.N.}} & 
		\multicolumn{1}{|c}{\textbf{Name of the index}} & \multicolumn{1}{|c|}{\textbf{Formula}} \\\hline
		
		\endfirsthead
		
		\hline
		\endlastfoot

		1&	\multirow{2}{5cm}{The first  hyper RL index} & \(HRL_1(G)=\sum\limits_{uv\in E(G)}\left(  d^2_G(u)+d^2_G(v)+d_G(u)d_G(v)\right)^2 \)\\[3ex]
		
		2&	\multirow{2}{5cm}{The second  hyper RL index} & \(HRL_2(G)=\sum\limits_{uv\in E(G)}\left(  d^2_G(u)+d^2_G(v)-d_G(u)d_G(v)\right)^2 \)\\[3ex]
		
		3&	\multirow{2}{5cm}{The first  Inverse RL index} & \(IRL_1(G)=\sum\limits_{uv\in E(G)}\dfrac{1}{\left(  d^2_G(u)+d^2_G(v)+d_G(u)d_G(v)\right)}  \)\\[3ex]
		
		4&	\multirow{2}{5cm}{The second  Inverse RL index} & \(IRL_2(G)=\sum\limits_{uv\in E(G)}\dfrac{1}{\left(  d^2_G(u)+d^2_G(v)-d_G(u)d_G(v)\right) } \)\\[3ex]
		
		5&	\multirow{2}{5cm}{The first general RL index} & \(GRL_1(G)=\sum\limits_{uv\in E(G)}\left(  d^2_G(u)+d^2_G(v)+d_G(u)d_G(v)\right)^a\)\\[3ex]
		
		6&	\multirow{2}{5cm}{The second general RL index} & \(GRL_2(G)=\sum\limits_{uv\in E(G)}\left(  d^2_G(u)+d^2_G(v)-d_G(u)d_G(v)\right)^a \)\\[3ex]
		7&	\multirow{2}{5cm}{The first  RL exponential} & \(RL_1(G,x)=\sum\limits_{uv\in E(G)}x^{\left(  d^2_G(u)+d^2_G(v)+d_G(u)d_G(v)\right)} \)\\[3ex]
		8&	\multirow{2}{6cm}{The second  RL exponential} & \(RL_2(G,x)=\sum\limits_{uv\in E(G)}x^{\left(  d^2_G(u)+d^2_G(v)-d_G(u)d_G(v)\right) }\)\\[3ex]
		9&	\multirow{2}{5cm}{The first  hyper RL exponential} & \(HRL_1(G,x)=\sum\limits_{uv\in E(G)}x^{\left(  d^2_G(u)+d^2_G(v)+d_G(u)d_G(v)\right)^2} \)\\[3ex]
		
		10&	\multirow{2}{6cm}{The second  hyper RL exponential} & \(HRL_2(G,x)=\sum\limits_{uv\in E(G)}x^{\left(  d^2_G(u)+d^2_G(v)-d_G(u)d_G(v)\right)^2 }\)\\[3ex]
		\hline
		11&	\multirow{2}{5cm}{The first  Inverse RL exponential} & \(IRL_1(G,x)=\sum\limits_{uv\in E(G)}x^{\dfrac{1}{\left(  d^2_G(u)+d^2_G(v)+d_G(u)d_G(v)\right)} } \)\\[3ex]
		12&	\multirow{2}{6cm}{The second  Inverse RL exponential} & \(IRL_2(G,x)=\sum\limits_{uv\in E(G)}x^{\dfrac{1}{\left(  d^2_G(u)+d^2_G(v)-d_G(u)d_G(v)\right) }} \)\\[3ex]
		13&	\multirow{2}{5cm}{The first general RL exponential} & \(GRL_1(G,x)=\sum\limits_{uv\in E(G)}x^{\left(  d^2_G(u)+d^2_G(v)+d_G(u)d_G(v)\right)^a}\)\\[3ex]
		
		14&	\multirow{2}{6cm}{The second general RL exponential} & \(GRL_2(G,x)=\sum\limits_{uv\in E(G)}x^{\left(  d^2_G(u)+d^2_G(v)-d_G(u)d_G(v)\right)^a} \)\\[3ex]

		
	\end{longtable}
	
	\begin{longtable}[c]{| c | l | l |}
		\caption{Multiplicative Rehan- Lanel indices/exponentials\label{table1b}}\\

		\hline
		\multicolumn{1}{|c}{\textbf{S.N.}} & 
		\multicolumn{1}{|c}{\textbf{Name of the index}} & \multicolumn{1}{|c|}{\textbf{Formula}} \\\hline
		\endfirsthead

		\hline
		\endlastfoot

		1&	\multirow{2}{6cm}{The first multiplicative  RL index} & \(MRL_1(G)=\prod\limits_{uv\in E(G)}\left(  d^2_G(u)+d^2_G(v)+d_G(u)d_G(v)\right)\)\\[3ex]
		
		2&	\multirow{2}{6cm}{The second multiplicative  RL index} & \(MRL_2(G)=\prod\limits_{uv\in E(G)}\left(  d^2_G(u)+d^2_G(v)-d_G(u)d_G(v)\right)\)\\[3ex]
		
		3&	\multirow{2}{6cm}{The first multiplicative hyper RL index} & \(MHRL_1(G)=\prod\limits_{uv\in E(G)}\left(  d^2_G(u)+d^2_G(v)+d_G(u)d_G(v)\right)^2 \)\\[3ex]
		
		4&	\multirow{2}{5cm}{The second multiplicative hyper RL index} & \(MHRL_2(G)=\prod\limits_{uv\in E(G)}\left(  d^2_G(u)+d^2_G(v)-d_G(u)d_G(v)\right)^2 \)\\[3ex]
		5&	\multirow{2}{5cm}{The first  multiplicative inverse RL index} & \(MIRL_1(G)=\prod\limits_{uv\in E(G)}\dfrac{1}{\left(  d^2_G(u)+d^2_G(v)+d_G(u)d_G(v)\right)}  \)\\[3ex]
		6&	\multirow{1}{7cm}{The second multiplicative inverse RL index} & \(MIRL_2(G)=\prod\limits_{uv\in E(G)}\dfrac{1}{\left(  d^2_G(u)+d^2_G(v)-d_G(u)d_G(v)\right) } \)\\[3ex]

		7&	\multirow{1}{7cm}{The first multiplicative general RL index} & \(MGRL_1(G)=\prod\limits_{uv\in E(G)}\left(  d^2_G(u)+d^2_G(v)+d_G(u)d_G(v)\right)^a\)\\[3ex]
		
		8&	\multirow{1}{7cm}{The second multiplicative general RL index} & \(MGRL_2(G)=\prod\limits_{uv\in E(G)}\left(  d^2_G(u)+d^2_G(v)-d_G(u)d_G(v)\right)^a \)\\[3ex]
		
		9&	\multirow{1}{6cm}{The first multiplicative  RL exponential} & \(MRL_1(G,x)=\prod\limits_{uv\in E(G)}x^{\left(  d^2_G(u)+d^2_G(v)+d_G(u)d_G(v)\right) }\)\\[3ex]
		
		10&	\multirow{1}{6cm}{The second multiplicative RL exponential} & \(MRL_2(G,x)=\prod\limits_{uv\in E(G)}x^{\left(  d^2_G(u)+d^2_G(v)-d_G(u)d_G(v)\right)} \)\\[3ex]
		11&	\multirow{1}{7cm}{The first multiplicative  hyper RL exponential} & \(MHRL_1(G,x)=\prod\limits_{uv\in E(G)}x^{\left(  d^2_G(u)+d^2_G(v)+d_G(u)d_G(v)\right)^2 }\)\\[3ex]
		
		12&	\multirow{1}{7cm}{The second multiplicative hyper RL exponential} & \(MHRL_2(G,x)=\prod\limits_{uv\in E(G)}x^{\left(  d^2_G(u)+d^2_G(v)-d_G(u)d_G(v)\right)^2} \)\\[3ex]
		13&	\multirow{1}{7cm}{The first multiplicative Inverse RL exponential} & \(MIRL_1(G,x)=\prod\limits_{uv\in E(G)}x^{\dfrac{1}{\left(  d^2_G(u)+d^2_G(v)+d_G(u)d_G(v)\right)} } \)\\[3ex]
		
		14&	\multirow{1}{7cm}{The second multiplicative Inverse RL exponential} & \(MIRL_2(G,x)=\prod\limits_{uv\in E(G)}x^{\dfrac{1}{\left(  d^2_G(u)+d^2_G(v)-d_G(u)d_G(v)\right) }} \)\\[4ex]
		\hline 
		15&	\multirow{1}{7cm}{The first multiplicative general RL exponential} & \(MGRL_1(G,x)=\prod\limits_{uv\in E(G)}x^{\left(  d^2_G(u)+d^2_G(v)+d_G(u)d_G(v)\right)^a}\)\\[3ex]
		16&	\multirow{1}{7cm}{The second multiplicative general RL exponential} & \(MGRL_2(G,x)=\prod\limits_{uv\in E(G)}x^{\left(  d^2_G(u)+d^2_G(v)-d_G(u)d_G(v)\right)^a} \)\\[3ex]

	\end{longtable}
		\begin{longtable}[c]{| c | l | l |}
		\caption{Rehan- Lanel indices/exponentials\label{table1c}}\\
		
		\hline
		\multicolumn{1}{|c}{\textbf{S.N.}} & 
		\multicolumn{1}{|c}{\textbf{Name of the index}} & \multicolumn{1}{|c|}{\textbf{Formula}} \\\hline
		
		\endfirsthead
		
		\hline
		\endlastfoot

		1&	\multirow{2}{5cm}{The third  hyper RL index} & \(HRL_1(G)=\sum\limits_{uv\in E(G)}\left(  d_G(u)-d_G(v)+d_G(u)d_G(v)\right)^2 \)\\[3ex]
		
		2&	\multirow{2}{5cm}{The fourth  hyper RL index} & \(HRL_2(G)=\sum\limits_{uv\in E(G)}\left(  |d_G(u)-d_G(v)|d_G(u)d_G(v)\right)^2 \)\\[3ex]
		3&	\multirow{2}{5cm}{The third  Inverse RL index} & \(IRL_1(G)=\sum\limits_{uv\in E(G)}\dfrac{1}{\left(   d_G(u)-d_G(v)+d_G(u)d_G(v)\right)}  \)\\[3ex]
		
		4&	\multirow{2}{5cm}{The fourth  Inverse RL index} & \(IRL_2(G)=\sum\limits_{uv\in E(G)}\dfrac{1}{\left( |d_G(u)-d_G(v)|d_G(u)d_G(v)\right) } \)\\[3ex]
		5&	\multirow{2}{5cm}{The third general RL index} & \(GRL_1(G)=\sum\limits_{uv\in E(G)}\left(   d_G(u)-d_G(v)+d_G(u)d_G(v)\right)^a\)\\[3ex]
		6&	\multirow{2}{5cm}{The fourth general RL index} & \(GRL_2(G)=\sum\limits_{uv\in E(G)}\left( |d_G(u)-d_G(v)|d_G(u)d_G(v)\right)^a \)\\[3ex]
		7&	\multirow{2}{5cm}{The third  RL exponential} & \(RL_1(G,x)=\sum\limits_{uv\in E(G)}x^{\left(   d_G(u)-d_G(v)+d_G(u)d_G(v)\right)} \)\\[3ex]
		
		8&	\multirow{2}{6cm}{The fourth  RL exponential} & \(RL_2(G,x)=\sum\limits_{uv\in E(G)}x^{\left( |d_G(u)-d_G(v)|d_G(u)d_G(v)\right) }\)\\[3ex]
		9&	\multirow{2}{5cm}{The third  hyper RL exponential} & \(HRL_1(G,x)=\sum\limits_{uv\in E(G)}x^{\left(   d_G(u)-d_G(v)+d_G(u)d_G(v)\right)^2} \)\\[3ex]
		
		10&	\multirow{2}{6cm}{The fourth  hyper RL exponential} & \(HRL_2(G,x)=\sum\limits_{uv\in E(G)}x^{\left( |d_G(u)-d_G(v)|d_G(u)d_G(v)\right)^2 }\)\\[3ex]
		
		11&	\multirow{2}{5cm}{The third  Inverse RL exponential} & \(IRL_1(G,x)=\sum\limits_{uv\in E(G)}x^{\dfrac{1}{\left(   d_G(u)-d_G(v)+d_G(u)d_G(v)\right)} } \)\\[3ex]
		12&	\multirow{2}{6cm}{The fourth  Inverse RL exponential} & \(IRL_2(G,x)=\sum\limits_{uv\in E(G)}x^{\dfrac{1}{\left(  |d_G(u)-d_G(v)|d_G(u)d_G(v)\right) }} \)\\[3ex]
		
		13&	\multirow{2}{5cm}{The third general RL exponential} & \(GRL_1(G,x)=\sum\limits_{uv\in E(G)}x^{\left(   d_G(u)-d_G(v)+d_G(u)d_G(v)\right)^a}\)\\[3ex]
		
		14&	\multirow{2}{6cm}{The fourth general RL exponential} & \(GRL_2(G,x)=\sum\limits_{uv\in E(G)}x^{\left( |d_G(u)-d_G(v)|d_G(u)d_G(v)\right)^a} \)\\[3ex]

		
	\end{longtable}
	
	\begin{longtable}[c]{| c | l | l |}
		\caption{Multiplicative Rehan- Lanel indices/exponentials\label{table1d}}\\

		\hline
		\multicolumn{1}{|c}{\textbf{S.N.}} & 
		\multicolumn{1}{|c}{\textbf{Name of the index}} & \multicolumn{1}{|c|}{\textbf{Formula}} \\\hline
		\endfirsthead

		\hline
		\endlastfoot

		1&	\multirow{2}{6cm}{The third multiplicative  RL index} & \(MRL_1(G)=\prod\limits_{uv\in E(G)}\left(  d_G(u)-d_G(v)+d_G(u)d_G(v)\right)\)\\[3ex]
		
		2&	\multirow{2}{6cm}{The fourth multiplicative  RL index} & \(MRL_2(G)=\prod\limits_{uv\in E(G)}\left( |d_G(u)-d_G(v)|d_G(u)d_G(v)\right)\)\\[3ex]
		\hline 
		3&	\multirow{2}{6cm}{The third multiplicative hyper RL index} & \(MHRL_1(G)=\prod\limits_{uv\in E(G)}\left(  d_G(u)-d_G(v)+d_G(u)d_G(v)\right)^2 \)\\[3ex]
		
		4&	\multirow{2}{5cm}{The fourth multiplicative hyper RL index} & \(MHRL_2(G)=\prod\limits_{uv\in E(G)}\left( |d_G(u)-d_G(v)|d_G(u)d_G(v)\right)^2 \)\\[3ex]
		5&	\multirow{2}{5cm}{The third  multiplicative inverse RL index} & \(MIRL_1(G)=\prod\limits_{uv\in E(G)}\dfrac{1}{\left(  d_G(u)-d_G(v)+d_G(u)d_G(v)\right)}  \)\\[3ex]
		6&	\multirow{1}{7cm}{The fourth multiplicative inverse RL index} & \(MIRL_2(G)=\prod\limits_{uv\in E(G)}\dfrac{1}{\left(  |d_G(u)-d_G(v)|d_G(u)d_G(v)\right) } \)\\[3ex]

		7&	\multirow{1}{7cm}{The third multiplicative general RL index} & \(MGRL_1(G)=\prod\limits_{uv\in E(G)}\left(  d_G(u)-d_G(v)+d_G(u)d_G(v)\right)^a\)\\[3ex]
		
		8&	\multirow{1}{7cm}{The fourth multiplicative general RL index} & \(MGRL_2(G)=\prod\limits_{uv\in E(G)}\left( |d_G(u)-d_G(v)|d_G(u)d_G(v)\right)^a \)\\[3ex]
		
		9&	\multirow{1}{6cm}{The third multiplicative  RL exponential} & \(MRL_1(G,x)=\prod\limits_{uv\in E(G)}x^{\left(  d_G(u)-d_G(v)+d_G(u)d_G(v)\right) }\)\\[3ex]
		
		10&	\multirow{1}{6cm}{The fourth multiplicative RL exponential} & \(MRL_2(G,x)=\prod\limits_{uv\in E(G)}x^{\left( |d_G(u)-d_G(v)|d_G(u)d_G(v)\right)} \)\\[3ex]
		11&	\multirow{1}{7cm}{The third multiplicative  hyper RL exponential} & \(MHRL_1(G,x)=\prod\limits_{uv\in E(G)}x^{\left(  d_G(u)-d_G(v)+d_G(u)d_G(v)\right)^2 }\)\\[3ex]
		
		12&	\multirow{1}{7cm}{The fourth multiplicative hyper RL exponential} & \(MHRL_2(G,x)=\prod\limits_{uv\in E(G)}x^{\left(  |d_G(u)-d_G(v)|d_G(u)d_G(v)\right)^2} \)\\[3ex]
		13&	\multirow{1}{7cm}{The third multiplicative Inverse RL exponential} & \(MIRL_1(G,x)=\prod\limits_{uv\in E(G)}x^{\dfrac{1}{\left( d_G(u)-d_G(v)+d_G(u)d_G(v)\right)} } \)\\[3ex]
		
		14&	\multirow{1}{7cm}{The fourth multiplicative Inverse RL exponential} & \(MIRL_2(G,x)=\prod\limits_{uv\in E(G)}x^{\dfrac{1}{\left(  |d_G(u)-d_G(v)|d_G(u)d_G(v)\right) }} \)\\[3ex]
		
		15&	\multirow{1}{7cm}{The third multiplicative general RL exponential} & \(MGRL_1(G,x)=\prod\limits_{uv\in E(G)}x^{\left(  d_G(u)-d_G(v)+d_G(u)d_G(v)\right)^a}\)\\[3ex]
		
		16&	\multirow{1}{7cm}{The fourth multiplicative general RL exponential} & \(MGRL_2(G,x)=\prod\limits_{uv\in E(G)}x^{\left(  |d_G(u)-d_G(v)|d_G(u)d_G(v)\right)^a} \)\\[3ex]

	\end{longtable}

	We determine the first, second, third and fourth Rehan-Lanel indices \((RL_1 \ \text{and} \ RL_2)\) for several standard graphs, as well as the Wheel graph and Sunflower graph.   In addition, we present the first, second, third and fourth Rehan-Lanel exponentials for the Wheel graph and Sunflower graph. 
	
	\renewcommand\qedsymbol{$\blacksquare$}
	\begin{proposition}
		Let \(G\) is an $r$-regular graph with \(n\) vertices and
		\(r \geq 2\). Then, 
		
		\begin{enumerate}[label=\roman*.]
			
			\item \(RL_1(G)=\dfrac{3nr^3}{2}\)
			\item \(RL_2(G)=\dfrac{nr^3}{2}\)
			 \item \(RL_3(G)=\dfrac{nr^3}{2}\)
			\item \(RL_4(G)=0\).
		\end{enumerate}
		
	\end{proposition}
	
	\begin{proof}
		Let \(G\) is an $r-$regular graph with \(n\) vertices and
		\(r \geq 2\).\\
		Then, \(G\) contains
		\(\dfrac{nr}{2}\)	edges. 
		From definition, we have
		$$ \begin{array}{ll}
			RL_1(G)&=\sum\limits_{uv\in E(G)} (r^2+r^2+r\cdot r)\\
			&=\dfrac{nr}{2} (3r^2)\\
			&=\dfrac{3nr^3}{2}.
		\end{array}$$
		
		$$ \begin{array}{ll}
			RL_2(G)&=\sum\limits_{uv\in E(G)} (r^2+r^2-r\cdot r)\\
			&=\dfrac{nr}{2} (r^2)\\
			&=\dfrac{nr^3}{2}.
		\end{array}$$
		
			$$ \begin{array}{ll}
			RL_3(G)&=\sum\limits_{uv\in E(G)} (r-r+r\cdot r)\\
			&=\dfrac{nr}{2} (r^2)\\
			&=\dfrac{nr^3}{2}.
		\end{array}$$
		
		$$ \begin{array}{ll}
			RL_4(G)&=\sum\limits_{uv\in E(G)} (|r-r|r\cdot r)\\
			&=\dfrac{nr}{2} (0)\\
			&=0.
		\end{array}$$
	\end{proof}
	
	\begin{corollary}
		Let	\(C_n\) be a cycle with \(n \geq 3\) vertices. Then, 
		\begin{enumerate}[label=\roman*.]
			\item \(RL_1(C_n)=12n.\)
			\item \(RL_2(C_n)=4n.\)
			\item \(RL_3(C_n)=4n.\)
			\item \(RL_4(C_n)=0.\)
		\end{enumerate} 
	\end{corollary}
	
	\begin{corollary}
		Let \(K_n\) be a complete graph with \(n\geq 3\)
		vertices. Then,
		\begin{enumerate}[label=\roman*.]
			\item \(RL_1(K_n)=\dfrac{3n(n-1)^3}{2}.\)
			\item \(RL_2(K_n)=\dfrac{n(n-1)^3}{2}.\)
			 \item \(RL_3(K_n)=\dfrac{n(n-1)^3}{2}.\)
			\item \(RL_4(K_n)=0.\)
		\end{enumerate}

	\end{corollary}
	
	\begin{proposition}
		Let \(P_n\) be a path with \(n\geq 3\) vertices. Then,
		\begin{enumerate}[label=\roman*.]
			\item \(RL_1(P_n)= 12n-22.\)
			\item \(RL_2(P_n)=4n-6.\)
			 \item \(RL_3(P_n)= 4n-6.\)
			\item \(RL_4(P_n)=4.\)
		\end{enumerate}
	\end{proposition}
	\begin{proof}
		The proof can be done similarly as above.
	\end{proof}
	\begin{proposition}
		Let $K_{m,n}$ be a complete bipartite graph with
		\(1\leq m\leq n\) and \(n \geq 2\). Then,	
		\begin{enumerate}[label=\roman*.]
			\item \(RL_1(K_{m,n})= mn(m^2+n^2+mn).\)
			\item \(RL_2(K_{m,n})= mn(m^2+n^2-mn).\)
			\item \(RL_3(K_{m,n})= mn(m-n+mn).\)
			\item \(RL_4(K_{m,n})= m^2n^2|m-n|.\)
		\end{enumerate}
		
	\end{proposition}
	\begin{proof}
		The proof can be done similarly as above.
	\end{proof}
	\subsubsection{Wheel Graph}
\begin{figure}[h]
		
		\centering
		\includegraphics[width=0.5\textwidth]{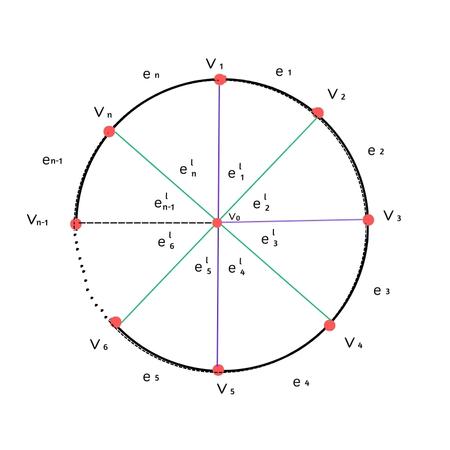}
		\caption{Wheel Graph}
		\label{fig:wh1}
	\end{figure}

	The wheel graph $W_n$ is determined by connecting $K_1$ and $C_n$.
	$K_1$ is the graph of order $1$ and $C_n$ is the cycle graph as shown in the Figure \ref{fig:wh1}.

	In \(W_n\), there are two types of edges as follows \(d(v)=n\) and \(d(v_i)=3\),
	$$\begin{array}{lll}
		E_1    & =\left\lbrace uv \in E(W_n)| d(u)=d(v)=3 \right\rbrace, & |E_1|=n. \\
		E_2&=\left\lbrace uv \in E(W_n)| d(u)=3 ,d(v)=n \right\rbrace , &|E_2|=n. 
	\end{array}$$

	\begin{theorem}
		Let \(W_n\) be a wheel graph. Then,
		
		\begin{enumerate}[label=\roman*.]
			\item \(RL_1(W_n)= n(n^2+3n+36)\)
			\item \(RL_2(W_n)= n(n^2-3n+18)\)
			\item \(RL_3(W_n)= 2n(2n+3)\)
			\item \(RL_4(W_n)= 3n^2|n-3|.\)
		\end{enumerate}
	
	\end{theorem}
	
	\begin{proof}
		$$\begin{array}{ll}
			RL_1(W_n)&= \sum\limits_{uv\in E(G)}\left(  d^2_G(u)+d^2_G(v)+d_G(u)d_G(v)\right)\\
			&= \sum\limits_{uv\in E_1(G)}\left(  d^2_G(u)+d^2_G(v)+d_G(u)d_G(v)\right)+\sum\limits_{uv\in E_2(G)}\left(  d^2_G(u)+d^2_G(v)+d_G(u)d_G(v)\right)\\
			&=n(3^2+3^2+3\cdot 3) + n(n^2+3^2 +n \cdot 3)\\
			&=n(n^2+3n+36).
		\end{array}$$
		
		$$\begin{array}{ll}
			RL_2(W_n)&= \sum\limits_{uv\in E(G)}\left(  d^2_G(u)+d^2_G(v)-d_G(u)d_G(v)\right)\\
			&= \sum\limits_{uv\in E_1(G)}\left(  d^2_G(u)+d^2_G(v)-d_G(u)d_G(v)\right)+\sum\limits_{uv\in E_2(G)}\left(  d^2_G(u)+d^2_G(v)-d_G(u)d_G(v)\right)\\
			&=n(3^2+3^2-3\cdot 3) + n(n^2+3^2 -n \cdot 3)\\
			&=n(n^2-3n+18).
		\end{array}$$
		
		 $$\begin{array}{ll}
			RL_3(W_n)&= \sum\limits_{uv\in E(G)}\left(  d_G(u)-d_G(v)+d_G(u)d_G(v)\right)\\
			&= \sum\limits_{uv\in E_1(G)}\left(  d_G(u)-d_G(v)+d_G(u)d_G(v)\right)+\sum\limits_{uv\in E_2(G)}\left(  d_G(u)-d_G(v)+d_G(u)d_G(v)\right)\\
			&=n(3-3+3\cdot 3) + n(n-3 +n \cdot 3)\\
			&=2n(2n+3).
		\end{array}$$
		
		$$\begin{array}{ll}
			RL_4(W_n)&= \sum\limits_{uv\in E(G)}\left(  |d_G(u)-d_G(v)|d_G(u)d_G(v)\right)\\
			&= \sum\limits_{uv\in E_1(G)}\left(  |d_G(u)-d_G(v)|d_G(u)d_G(v)\right)+\sum\limits_{uv\in E_2(G)}\left(  |d_G(u)-d_G(v)|d_G(u)d_G(v)\right)\\
			&=n(|3-3|3\cdot 3) + n(|n-3|n \cdot 3)\\
			&=3n^2|n-3|.
		\end{array}$$
	\end{proof}
	
	\begin{theorem}
		Let \(W_n\) be a wheel graph. Then,
		
		\begin{itemize}
			\item[i.] \(RL_1(W_n,x)=n(x^{27}+x^{n^2+3n+9})\)
			\item[ii.]\(RL_2(W_n,x)=nx^9(x^{n^2-3n}+1)\)
				\item[iii.] \(RL_3(W_n,x)=n(x^{9}+x^{4n-3})\)
			\item[iv.]\(RL_4(W_n,x)=n(x^{3n|n-3|}+1)\).
		\end{itemize}
	\end{theorem}
	
	\begin{proof}
		$$\begin{array}{ll}
			RL_1(W_n,x)&= \sum\limits_{uv\in E(G)}x^{\left(  d^2_G(u)+d^2_G(v)+d_G(u)d_G(v)\right)}\\
			&= \sum\limits_{uv\in E_1(G)}x^{\left(  d^2_G(u)+d^2_G(v)+d_G(u)d_G(v)\right)}+\sum\limits_{uv\in E_2(G)}x^{\left(  d^2_G(u)+d^2_G(v)+d_G(u)d_G(v)\right)}\\
			&=nx^{(3^2+3^2+3\cdot 3)} + nx^{(n^2+3^2 +n \cdot 3)}\\
			&=n(x^{27}+x^{n^2+3n+9}).
		\end{array}$$
		
		$$\begin{array}{ll}
			RL_2(W_n,x)&= \sum\limits_{uv\in E(G)}x^{\left(  d^2_G(u)+d^2_G(v)-d_G(u)d_G(v)\right)}\\
			&= \sum\limits_{uv\in E_1(G)}x^{\left(  d^2_G(u)+d^2_G(v)-d_G(u)d_G(v)\right)}+\sum\limits_{uv\in E_2(G)}x^{\left(  d^2_G(u)+d^2_G(v)-d_G(u)d_G(v)\right)}\\
			&=nx^{(3^2+3^2-3\cdot 3)} + nx^{(n^2+3^2 -n \cdot 3)}\\
			&=n(x^{9}+x^{n^2-3n+9})\\
			&=nx^9(x^{n^2-3n}+1).
		\end{array}$$
		
		     $$\begin{array}{ll}
			RL_3(W_n,x)&= \sum\limits_{uv\in E(G)}x^{\left(  d_G(u)-d_G(v)+d_G(u)d_G(v)\right)}\\
			&= \sum\limits_{uv\in E_1(G)}x^{\left(  d_G(u)-d_G(v)+d_G(u)d_G(v)\right)}+\sum\limits_{uv\in E_2(G)}x^{\left(  d_G(u)-d_G(v)+d_G(u)d_G(v)\right)}\\
			&=nx^{(3-3+3\cdot 3)} + nx^{(n-3 +n \cdot 3)}\\
			&=n(x^{9}+x^{4n-3}).
		\end{array}$$
		
		$$\begin{array}{ll}
			RL_4(W_n,x)&= \sum\limits_{uv\in E(G)}x^{\left(  |d_G(u)-d_G(v)|d_G(u)d_G(v)\right)}\\
			&= \sum\limits_{uv\in E_1(G)}x^{\left(  |d_G(u)-d_G(v)|d_G(u)d_G(v)\right)}+\sum\limits_{uv\in E_2(G)}x^{\left(  |d_G(u)-d_G(v)|d_G(u)d_G(v)\right)}\\
			&=nx^{(|3-3|\cdot 3\cdot 3)} + nx^{(|n-3|\cdot n \cdot 3)}\\
			&=n(x^{0}+x^{3n|n-3|})\\
			&=n(x^{3n|n-3|}+1).
		\end{array}$$
	\end{proof}
	
	\subsubsection{Sunflower Graph}
\begin{figure}[h]
		
		\centering
		\includegraphics[width=0.5\textwidth]{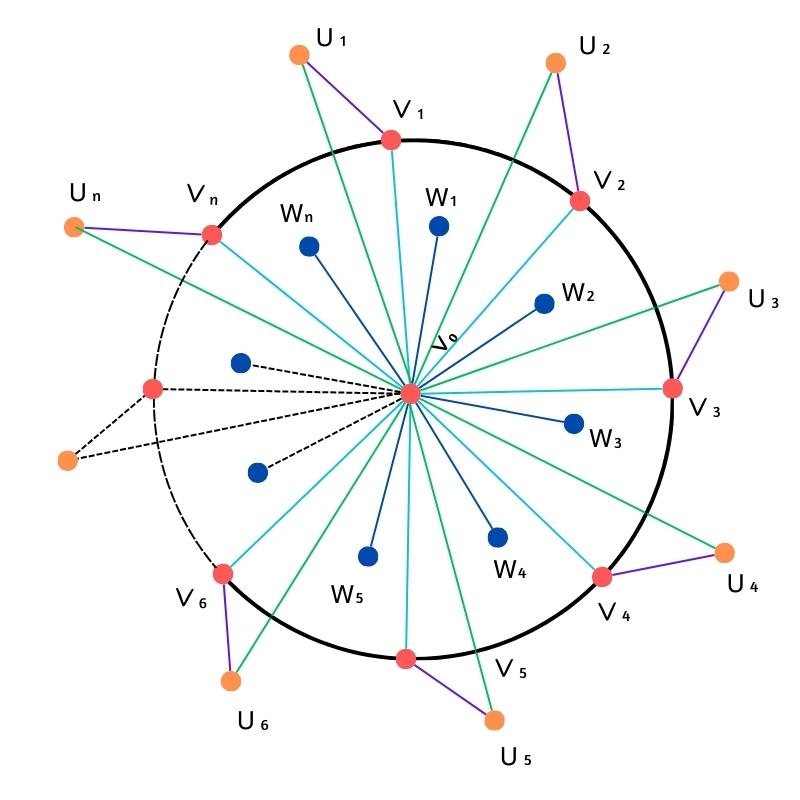}
		\caption{Sunflower Graph}
		\label{fig:sf1}
	\end{figure}

	The sunflower graph \(S\!f_n\) is derived from the flower graph by extending \(n\) distinct edges to the apex of the flower graph, as depicted in Figure \ref{fig:sf1}.

	In \({S\!f}_n\), there are five types of edges as follows,
	$$\begin{array}{lll}
		E_1&=\left\lbrace uv \in E({S\!f}_n)| d(v)=d(v)=4 \right\rbrace,& |E_1|=n.\\
		E_2&=\left\lbrace uv \in E({S\!f}_n)| d(v_0)=3n ,d(v)=4 \right\rbrace , &|E_2|=n.\\
		E_3&=\left\lbrace uv \in E({S\!f}_n)| d(u)=2, d(v)=4 \right\rbrace,& |E_3|=n\\
		E_4&=\left\lbrace uv \in E({S\!f}_n)| d(v_0)=3n ,d(v)=2 \right\rbrace , &|E_4|=n.\\
		E_5&=\left\lbrace uv \in E({S\!f}_n)| d(w)=1 ,d(v_0)=3n \right\rbrace ,& |E_5|=n.
		
	\end{array}$$
	
	\begin{theorem}
		Let \(S\!f_n\) be a sunflower graph. Then,
		
		\begin{enumerate}[label=\roman*.]
			\item \(RL_1(S\!f_n)= n(27n^2+21n+97)\)
			\item \(RL_2(S\!f_n)= n(27n^2-21n+49)\)
			\item \(RL_3(S\!f_n)= n(25n+17)\)
			\item \(RL_4(S\!f_n)= n(12n|3n-4|+6n|3n-2|+3n|3n-1|+16).\)
		\end{enumerate}
	
	\end{theorem}
	
	\begin{proof}
		$$\begin{array}{ll}
			RL_1(S\!f_n)&= \sum\limits_{uv\in E(G)}\left(  d^2_G(u)+d^2_G(v)+d_G(u)d_G(v)\right)\\
			&= \sum\limits_{uv\in E_1(G)}\left(  d^2_G(u)+d^2_G(v)+d_G(u)d_G(v)\right)+\sum\limits_{uv\in E_2(G)}\left(  d^2_G(u)+d^2_G(v)+d_G(u)d_G(v)\right)+ \\ & \sum\limits_{uv\in E_3(G)}\left(  d^2_G(u)+d^2_G(v)+d_G(u)d_G(v)\right)+\sum\limits_{uv\in E_4(G)}\left(  d^2_G(u)+d^2_G(v)+d_G(u)d_G(v)\right)\\ &\sum\limits_{uv\in E_5(G)}\left(  d^2_G(u)+d^2_G(v)+d_G(u)d_G(v)\right)\\
			&=n(4^2+4^2+4\cdot 4) + n((3n)^2+4^2 +3n \cdot 4)+ n(2^2+4^2+2 \cdot 4)+n((3n)^2+2^2+3n \cdot 2)\\ &  n(1^2+(3n)^2+1 \cdot 3n)\\
			&=n(27n^2+21n+97).
		\end{array}$$
		
		$$\begin{array}{ll}
			RL_2(S\!f_n)&= \sum\limits_{uv\in E(G)}\left(  d^2_G(u)+d^2_G(v)-d_G(u)d_G(v)\right)\\
			&= \sum\limits_{uv\in E_1(G)}\left(  d^2_G(u)+d^2_G(v)-d_G(u)d_G(v)\right)+\sum\limits_{uv\in E_2(G)}\left(  d^2_G(u)+d^2_G(v)-d_G(u)d_G(v)\right)+ \\ & \sum\limits_{uv\in E_3(G)}\left(  d^2_G(u)+d^2_G(v)-d_G(u)d_G(v)\right)+\sum\limits_{uv\in E_4(G)}\left(  d^2_G(u)+d^2_G(v)-d_G(u)d_G(v)\right)\\ &+\sum\limits_{uv\in E_5(G)}\left(  d^2_G(u)+d^2_G(v)-d_G(u)d_G(v)\right)\\
			&=n(4^2+4^2-4\cdot 4) + n((3n)^2+4^2 -3n \cdot 4)+ n(2^2+4^2-2 \cdot 4)+n((3n)^2+2^2-3n \cdot 2)\\ & + n(1^2+(3n)^2-1 \cdot 3n)\\
			&=n(27n^2-21n+49).
		\end{array}$$ 
		
		  $$\begin{array}{ll}
			RL_3(S\!f_n)&= \sum\limits_{uv\in E(G)}\left(  d_G(u)+d_G(v)+d_G(u)d_G(v)\right)\\
			&= \sum\limits_{uv\in E_1(G)}\left(   d_G(u)-d_G(v)+d_G(u)d_G(v)\right)+\sum\limits_{uv\in E_2(G)}\left(   d_G(u)-d_G(v)+d_G(u)d_G(v)\right)+ \\ & \sum\limits_{uv\in E_3(G)}\left(   d_G(u)-d_G(v)+d_G(u)d_G(v)\right)+\sum\limits_{uv\in E_4(G)}\left(   d_G(u)-d_G(v)+d_G(u)d_G(v)\right)\\ &\sum\limits_{uv\in E_5(G)}\left(   d_G(u)-d_G(v)+d_G(u)d_G(v)\right)\\
			&=n(4-4+4\cdot 4) + n(3n-4 +3n \cdot 4)+ n(2-4+2 \cdot 4)+n(3n-2+3n \cdot 2)\\ &  n(1-3n+1 \cdot 3n)\\
			&=n(25n+17).
		\end{array}$$
		
		$$\begin{array}{ll}
			RL_4(S\!f_n)&= \sum\limits_{uv\in E(G)}\left(  |d_G(u)-d_G(v)|d_G(u)d_G(v)\right)\\
			&= \sum\limits_{uv\in E_1(G)}\left(  |d_G(u)-d_G(v)|d_G(u)d_G(v)\right)+\sum\limits_{uv\in E_2(G)}\left(  |d_G(u)-d_G(v)|d_G(u)d_G(v)\right)+ \\ & \sum\limits_{uv\in E_3(G)}\left(  |d_G(u)-d_G(v)|d_G(u)d_G(v)\right)+\sum\limits_{uv\in E_4(G)}\left(  |d_G(u)-d_G(v)|d_G(u)d_G(v)\right)\\ &+\sum\limits_{uv\in E_5(G)}\left(  |d_G(u)-d_G(v)|d_G(u)d_G(v)\right)\\
			&=n(|4-4|\cdot 4\cdot 4) + n(|3n-4|\cdot 3n \cdot 4)+ n(|2-4|\cdot 2 \cdot 4)+n(|3n-2|\cdot 3n \cdot 2)\\ & + n(|1-3n|\cdot 1 \cdot 3n)\\
			&=n(12n|3n-4|+6n|3n-2|+3n|3n-1|+16).
		\end{array}$$ 
	\end{proof}

	\begin{proposition}
		Let \(S\!f_n\) be a sunflower graph. Then,
		
		\begin{itemize}
			\item[i.] \(RL_1(S\!f_n,x)= n(x^{48}+x^{9n^2+12n+16}+x^{28}+x^{9n^2+6n+4}+x^{9n^2+3n+1}).\)
			\item[ii.] \( RL_2(S\!f_n,x)=n(x^{16}+x^{9n^2-12n+16}+x^{12}+x^{9n^2-6n+4}+x^{9n^2-3n+1}).\)
				\item[iii.] \(RL_3(S\!f_n,x)= n(x^{16}+x^{15n-4}+x^{6}+x^{9n-2}+x).\)
			\item[iv.] \( RL_4(S\!f_n,x)=n(x^{16}+x^{12n|3n-4|}+x^{6n|3n-2|}+x^{3n|3n-1|}+1).\)
		\end{itemize}

	\end{proposition}
	
	\begin{proof}
		$$\begin{array}{ll}
			RL_1(S\!f_n,x)&= \sum\limits_{uv\in E(G)}x^{\left(  d^2_G(u)+d^2_G(v)+d_G(u)d_G(v)\right)}\\
			&= \sum\limits_{uv\in E_1(G)}x^{\left(  d^2_G(u)+d^2_G(v)+d_G(u)d_G(v)\right)}+\sum\limits_{uv\in E_2(G)}x^{\left(  d^2_G(u)+d^2_G(v)+d_G(u)d_G(v)\right)}+  \sum\limits_{uv\in E_3(G)}x^{\left(  d^2_G(u)+d^2_G(v)+d_G(u)d_G(v)\right)}\\& +\sum\limits_{uv\in E_4(G)}x^{\left(  d^2_G(u)+d^2_G(v)+d_G(u)d_G(v)\right)}+\sum\limits_{uv\in E_5(G)}x^{\left(  d^2_G(u)+d^2_G(v)+d_G(u)d_G(v)\right)}\\
			&=nx^{(4^2+4^2+4\cdot 4)} + nx^{((3n)^2+4^2 +3n \cdot 4)}+ nx^{(2^2+4^2+2 \cdot 4)}+nx^{((3n)^2+2^2+3n \cdot 2)}+ nx^{(1^2+(3n)^2+1 \cdot 3n)}\\
			&=n(x^{48}+x^{9n^2+12n+16}+x^{28}+x^{9n^2+6n+4}+x^{9n^2+3n+1}).
		\end{array}$$
		
		$$\begin{array}{ll}
			RL_2(S\!f_n,x)&= \sum\limits_{uv\in E(G)}x^{\left(  d^2_G(u)+d^2_G(v)-d_G(u)d_G(v)\right)}\\
			&= \sum\limits_{uv\in E_1(G)}x^{\left(  d^2_G(u)+d^2_G(v)-d_G(u)d_G(v)\right)}+\sum\limits_{uv\in E_2(G)}x^{\left(  d^2_G(u)+d^2_G(v)-d_G(u)d_G(v)\right)}+  \sum\limits_{uv\in E_3(G)}x^{\left(  d^2_G(u)+d^2_G(v)-d_G(u)d_G(v)\right)}\\& +\sum\limits_{uv\in E_4(G)}x^{\left(  d^2_G(u)+d^2_G(v)-d_G(u)d_G(v)\right)}+\sum\limits_{uv\in E_5(G)}x^{\left(  d^2_G(u)+d^2_G(v)-d_G(u)d_G(v)\right)}\\
			&=nx^{(4^2+4^2-4\cdot 4)} + nx^{((3n)^2+4^2 -3n \cdot 4)}+ nx^{(2^2+4^2-2 \cdot 4)}+nx^{((3n)^2+2^2-3n \cdot 2)}+ nx^{(1^2+(3n)^2-1 \cdot 3n)}\\
			&=n(x^{16}+x^{9n^2-12n+16}+x^{12}+x^{9n^2-6n+4}+x^{9n^2-3n+1}.)
		\end{array}$$
		
		  $$\begin{array}{ll}
			RL_3(S\!f_n,x)&= \sum\limits_{uv\in E(G)}x^{\left(  d_G(u)-d_G(v)+d_G(u)d_G(v)\right)}\\
			&= \sum\limits_{uv\in E_1(G)}x^{\left(  d_G(u)-d_G(v)+d_G(u)d_G(v)\right)}+\sum\limits_{uv\in E_2(G)}x^{\left(  d_G(u)-d_G(v)+d_G(u)d_G(v)\right)}+  \sum\limits_{uv\in E_3(G)}x^{\left(  d_G(u)-d_G(v)+d_G(u)d_G(v)\right)}\\& +\sum\limits_{uv\in E_4(G)}x^{\left(  d_G(u)-d_G(v)+d_G(u)d_G(v)\right)}+\sum\limits_{uv\in E_5(G)}x^{\left(  d_G(u)-d_G(v)+d_G(u)d_G(v)\right)}\\
			&=nx^{(4^2+4^2+4\cdot 4)} + nx^{((3n)^2+4^2 +3n \cdot 4)}+ nx^{(2^2+4^2+2 \cdot 4)}+nx^{((3n)^2+2^2+3n \cdot 2)}+ nx^{(1^2+(3n)^2+1 \cdot 3n)}\\
			&=n(x^{48}+x^{9n^2+12n+16}+x^{28}+x^{9n^2+6n+4}+x^{9n^2+3n+1}).
		\end{array}$$
		
		$$\begin{array}{ll}
			RL_4(S\!f_n,x)&= \sum\limits_{uv\in E(G)}x^{\left(  |d_G(u)-d_G(v)|d_G(u)d_G(v)\right)}\\
			&= \sum\limits_{uv\in E_1(G)}x^{\left(  |d_G(u)-d_G(v)|d_G(u)d_G(v)\right)}+\sum\limits_{uv\in E_2(G)}x^{\left(  |d_G(u)-d_G(v)|d_G(u)d_G(v)\right)}+  \sum\limits_{uv\in E_3(G)}x^{\left(  |d_G(u)-d_G(v)|d_G(u)d_G(v)\right)}\\& +\sum\limits_{uv\in E_4(G)}x^{\left(  |d_G(u)-d_G(v)|d_G(u)d_G(v)\right)}+\sum\limits_{uv\in E_5(G)}x^{\left(  |d_G(u)-d_G(v)|d_G(u)d_G(v)\right)}\\
			&=nx^{(4^2+4^2-4\cdot 4)} + nx^{((3n)^2+4^2 -3n \cdot 4)}+ nx^{(2^2+4^2-2 \cdot 4)}+nx^{((3n)^2+2^2-3n \cdot 2)}+ nx^{(1^2+(3n)^2-1 \cdot 3n)}\\
			&=n(x^{16}+x^{9n^2-12n+16}+x^{12}+x^{9n^2-6n+4}+x^{9n^2-3n+1}.)
		\end{array}$$
	\end{proof}
	
	By following a similar procedure, we can determine the results for the different indices/exponentials listed in the tables in the section 3.1 . 
	
	\subsection{Banhatti RL indices }
	
	Kulli \cite{article34} defined the Banhatti degree of a vertex \(u\) of a graph \(G\) as, \[ B(u)=\dfrac{d_G(e)}{n-d_G(u)},\]
	where \(\left| V(G)\right|=n \) and the vertex \(u\) and edge \(e\) are incident in \(G\).

	In \cite{article34}	Kulli proposed the first and second E-Banhatti indices of a graph G and they are defined as,
	
	\[EB_1(G)=\sum_{uv\in E(G)}^{}B(u)+B(v) \ \text{and} \ EB_2(G)=\sum_{uv\in E(G)}^{}B(u)B(v). \]

	Inspired by this, we present the first and second E-Banhatti Rehan-Lanel indices.  
	
	\begin{definition}
		The first and second E-Banhatti Rehan-Lanel indices of a graph \(G\) are defined by, 
	\begin{enumerate}[label=\roman*.]
		\item \(BRL_1(G)=\sum\limits_{uv\in E(G)}(B^2(u)+B^2(v)+B(u)B(v))\)
		
		\item \(BRL_2(G)=\sum\limits_{uv\in E(G)}(B^2(u)+B^2(v)-B(u)B(v))\)
		\item \([BRL_3(G)=\sum\limits_{uv\in E(G)}(B(u)-B(v)+B(u)B(v))\)
		\item \(BRL_2(G)=\sum\limits_{uv\in E(G)}(|B(u)-B(v)|B(u)B(v)).\)
	\end{enumerate}	
	\end{definition}

	Expanding on these indices, we moreover propose additional novel indices, as outlined in the Table \ref{tableb1}, \ref{tableb2}, \ref{tableb3} and \ref{tableb4}.\\
	
.
	\begin{longtable}[c]{| c | l | l |}
		\caption{\small E-Banhatti Rehan- Lanel  indices/exponentials \label{tableb1}}\\
		
		\hline
		\multicolumn{1}{|c}{\textbf{S.N.}} & 
		\multicolumn{1}{|c}{\textbf{Name of the index}} & \multicolumn{1}{|c|}{\textbf{Formula}} \\\hline
		\endfirsthead
		
		\hline
		\endlastfoot

		1&	The first hyper E-Banhatti RL index & \(HBRL_1(G)=\sum\limits_{uv\in E(G)}\left(  B^2(u)+B^2(v)+B(u)B(v)\right)^2 \)\\[1ex]
		2&	The second  hyper E- Banhatti RL index & \(HBRL_2(G)=\sum\limits_{uv\in E(G)}\left(  B^2(u)+B^2(v)-B(u)B(v)\right)^2 \)\\[3ex]
		3&	\multirow{1}{6cm}{The first inverse E Banhatti RL index} & \(IBRL_1(G)=\sum\limits_{uv\in E(G)}\dfrac{1}{\left(  B^2(u)+B^2(v)+B(u)B(v)\right)} \)\\[3ex]
		
		4&	\multirow{1}{6cm}{The second inverse E Banhatti RL index} & \(IBRL_2(G)=\sum\limits_{uv\in E(G)}\dfrac{1}{\left(  B^2(u)+B^2(v)-B(u)B(v)\right)} \)\\[3ex]
		
		5&	\multirow{1}{6cm}{The first general E Banhatti RL index} & \(GBRL_1(G)=\sum\limits_{uv\in E(G)}\left(  B^2(u)+B^2(v)+B(u)B(v)\right)^a \)\\[3ex]
		6&	\multirow{1}{6cm}{The second  general E Banhatti RL index} & \(GBRL_2(G)=\sum\limits_{uv\in E(G)}\left(  B^2(u)+B^2(v)-B(u)B(v)\right)^a \)\\[3ex]
		
		7&	\multirow{1}{6cm}{The first E Banhatti RL exponential} & \(BRL_1(G,x)=\sum\limits_{uv\in E(G)}x^{\left(  B^2(u)+B^2(v)+B(u)B(v)\right)} \)\\[3ex]
		
		8&	\multirow{1}{6cm}{The second  E Banhatti RL exponential} & \(BRL_2(G,x)=\sum\limits_{uv\in E(G)}x^{\left(  B^2(u)+B^2(v)-B(u)B(v)\right)} \)\\[3ex]
		
		9&	\multirow{1}{7cm}{The first hyper E Banhatti RL exponential} & \(HBRL_1(G,x)=\sum\limits_{uv\in E(G)}x^{\left(  B^2(u)+B^2(v)+B(u)B(v)\right)^2 }\)\\[3ex]
		
		10&	\multirow{1}{7cm}{The second  hyper E Banhatti RL exponential} & \(HBRL_2(G,x)=\sum\limits_{uv\in E(G)}x^{\left(  B^2(u)+B^2(v)-B(u)B(v)\right)^2} \)\\[3ex]
		11&	\multirow{1}{7cm}{The first inverse E Banhatti RL exponential} & \(IBRL_1(G,x)=\sum\limits_{uv\in E(G)}x^{\frac{1}{\left(  B^2(u)+B^2(v)+B(u)B(v)\right)}} \)\\[3ex]
		
		12&	\multirow{1}{7cm}{The second inverse E Banhatti RL exponential} & \(IBRL_2(G,x)=\sum\limits_{uv\in E(G)}x^{\frac{1}{\left(  B^2(u)+B^2(v)-B(u)B(v)\right)} }\)\\[3ex]
		\hline 
		13&	\multirow{1}{7cm}{The first general E Banhatti RL exponential} & \(GBRL_1(G,x)=\sum\limits_{uv\in E(G)}x^{\left(  B^2(u)+B^2(v)+B(u)B(v)\right)^a} \)\\[3ex]
		
		14&	\multirow{1}{7cm}{The second  general E Banhatti RL exponential} & \(GBRL_2(G,x)=\sum\limits_{uv\in E(G)}x^{\left(  B^2(u)+B^2(v)-B(u)B(v)\right)^a} \)\\[3ex]

	\end{longtable}
	
	\begin{longtable}[c]{| c | l | l |}
		\caption{\small Multiplicative E- Banhatti Rehan- Lanel indices \label{tableb2}}\\
		
		\hline
		\multicolumn{1}{|c}{\textbf{S.N.}} & 
		\multicolumn{1}{|c}{\textbf{Name of the index}} & \multicolumn{1}{|c|}{\textbf{Formula}} \\\hline
		\endfirsthead
		
		\hline
		\endlastfoot

		1&	The first multiplivative E Banhatti RL index & \(BRL_1(G)=\prod\limits_{uv\in E(G)}\left(  B^2(u)+B^2(v)+B(u)B(v)\right) \)\\[3ex]
		
		2&	The second multiplivative E Banhatti RL index & \(BRL_2(G)=\prod\limits_{uv\in E(G)}\left(  B^2(u)+B^2(v)-B(u)B(v)\right) \)\\[3ex]
		
		3&		\multirow{2}{6cm}{The first  multiplicative hyper Banhatti RL index} & \(HBRL_1(G)=\prod\limits_{uv\in E(G)}\left(  B^2(u)+B^2(v)+B(u)B(v)\right)^2 \)\\[5ex]
		
		4&		\multirow{2}{7cm}{The second multiplicative hyper Banhatti RL index} & \(HBRL_2(G)=\prod\limits_{uv\in E(G)}\left(  B^2(u)+B^2(v)-B(u)B(v)\right)^2 \)\\[5ex]
		
		5&	\multirow{1}{7cm}{The first multiplicative inverse Banhatti RL index} & \(IBRL_1(G)=\prod\limits_{uv\in E(G)}\dfrac{1}{\left(  B^2(u)+B^2(v)+B(u)B(v)\right)} \)\\[3ex]
		6&	\multirow{2}{5cm}{The second multiplicative inverse Banhatti RL index} & \(IBRL_2(G)=\prod\limits_{uv\in E(G)}\dfrac{1}{\left(  B^2(u)+B^2(v)-B(u)B(v)\right)} \)\\[6ex]
		
		7&	\multirow{2}{7cm}{The first general multiplicative Banhatti RL index} & \(GBRL_1(G)=\prod\limits_{uv\in E(G)}\left(  B^2(u)+B^2(v)+B(u)B(v)\right)^a \)\\[6ex]
		8&	\multirow{2}{7cm}{The second  general multiplicative Banhatti RL index} & \(GBRL_2(G)=\prod\limits_{uv\in E(G)}\left(  B^2(u)+B^2(v)-B(u)B(v)\right)^a \)\\[6ex]
		
		9&	\multirow{2}{5cm}{The first multiplicative Banhatti RL exponential} & \(BRL_1(G,x)=\prod\limits_{uv\in E(G)}x^{\left(  B^2(u)+B^2(v)+B(u)B(v)\right)} \)\\[3ex]
		
		10&	\multirow{2}{5cm}{The second multiplicative Banhatti RL exponential} & \(BRL_2(G,x)=\prod\limits_{uv\in E(G)}x^{\left(  B^2(u)+B^2(v)-B(u)B(v)\right)} \)\\[3ex]
		
		11&	\multirow{2}{5cm}{The first multiplicative hyper Banhatti RL exponential} & \(HBRL_1(G,x)=\prod\limits_{uv\in E(G)}x^{\left(  B^2(u)+B^2(v)+B(u)B(v)\right)^2 }\)\\[5ex]
		
		12&	\multirow{2}{5cm}{The second multiplicative hyper Banhatti RL exponential} & \(HBRL_2(G,x)=\prod\limits_{uv\in E(G)}x^{\left(  B^2(u)+B^2(v)-B(u)B(v)\right)^2} \)\\[3ex]
		13&	\multirow{2}{5cm}{The first multiplicative inverse Banhatti RL exponential} & \(IBRL_1(G,x )=\prod\limits_{uv\in E(G)}x^{\dfrac{1}{\left(  B^2(u)+B^2(v)+B(u)B(v)\right)}} \)\\[3ex]
		14&	\multirow{2}{5cm}{The second multiplicative inverse Banhatti RL exponential} & \(IBRL_2(G)=\prod\limits_{uv\in E(G)}x^{\dfrac{1}{\left(  B^2(u)+B^2(v)-B(u)B(v)\right)} }\)\\[6ex]
		
		15&	\multirow{2}{5cm}{The first multiplicative general Banhatti RL exponential} & \(GBRL_1(G,x)=\prod\limits_{uv\in E(G)}x^{\left(  B^2(u)+B^2(v)+B(u)B(v)\right)^a} \)\\[6ex]
		\hline
		16&	\multirow{2}{5cm}{The second multiplicative general Banhatti RL exponential} & \(GBRL_2(G,x)= \prod\limits_{uv\in E(G)}\displaystyle{x^{\left(  B^2(u)+B^2(v)-B(u)B(v)\right)^a}} \)\\[6ex]
		
		\hline
		
	\end{longtable}
	\begin{longtable}[c]{| c | l | l |}
		\caption{\small E-Banhatti Rehan- Lanel  indices/exponentials \label{tableb3}}\\
		
		\hline
		\multicolumn{1}{|c}{\textbf{S.N.}} & 
		\multicolumn{1}{|c}{\textbf{Name of the index}} & \multicolumn{1}{|c|}{\textbf{Formula}} \\\hline
		\endfirsthead
		
		\hline
		\endlastfoot

		1&	The third hyper E-Banhatti RL index & \(HBRL_1(G)=\sum\limits_{uv\in E(G)}\left( B(u)-B(v)+B(u)B(v)\right)^2 \)\\[3ex]
		
		2&	The fourth  hyper E- Banhatti RL index & \(HBRL_2(G)=\sum\limits_{uv\in E(G)}\left( |B(u)-B(v)|B(u)B(v)\right)^2 \)\\[3ex]
		
		3&	\multirow{1}{6cm}{The third inverse E Banhatti RL index} & \(IBRL_1(G)=\sum\limits_{uv\in E(G)}\dfrac{1}{\left(  B(u)-B(v)+B(u)B(v)\right)} \)\\[3ex]
		
		4&	\multirow{1}{6cm}{The fourth inverse E Banhatti RL index} & \(IBRL_2(G)=\sum\limits_{uv\in E(G)}\dfrac{1}{\left( |B(u)-B(v)|B(u)B(v)\right)} \)\\[3ex]
		
		5&	\multirow{1}{6cm}{The third general E Banhatti RL index} & \(GBRL_1(G)=\sum\limits_{uv\in E(G)}\left(  B(u)-B(v)+B(u)B(v)\right)^a \)\\[3ex]
		
		6&	\multirow{1}{6cm}{The fourth  general E Banhatti RL index} & \(GBRL_2(G)=\sum\limits_{uv\in E(G)}\left(  |B(u)-B(v)|B(u)B(v)\right)^a \)\\[3ex]
		
		7&	\multirow{1}{6cm}{The third E Banhatti RL exponential} & \(BRL_1(G,x)=\sum\limits_{uv\in E(G)}x^{\left(  B(u)-B(v)+B(u)B(v)\right)} \)\\[3ex]
		
		8&	\multirow{1}{6cm}{The fourth  E Banhatti RL exponential} & \(BRL_2(G,x)=\sum\limits_{uv\in E(G)}x^{\left( |B(u)-B(v)|B(u)B(v)\right)} \)\\[3ex]
		
		9&	\multirow{1}{7cm}{The third hyper E Banhatti RL exponential} & \(HBRL_1(G,x)=\sum\limits_{uv\in E(G)}x^{\left(  B(u)-B(v)+B(u)B(v)\right)^2 }\)\\[3ex]
		
		10&	\multirow{1}{7cm}{The fourth  hyper E Banhatti RL exponential} & \(HBRL_2(G,x)=\sum\limits_{uv\in E(G)}x^{\left( |B(u)-B(v)|B(u)B(v)\right)^2} \)\\[3ex]
		
		11&	\multirow{1}{7cm}{The third inverse E Banhatti RL exponential} & \(IBRL_1(G,x)=\sum\limits_{uv\in E(G)}x^{\frac{1}{\left(  B(u)-B(v)+B(u)B(v)\right)}} \)\\[3ex]
		
		12&	\multirow{1}{7cm}{The fourth inverse E Banhatti RL exponential} & \(IBRL_2(G,x)=\sum\limits_{uv\in E(G)}x^{\frac{1}{\left(  |B(u)-B(v)|B(u)B(v)\right)} }\)\\[3ex]
		
		13&	\multirow{1}{7cm}{The third general E Banhatti RL exponential} & \(GBRL_1(G,x)=\sum\limits_{uv\in E(G)}x^{\left(  B(u)-B(v)+B(u)B(v)\right)^a} \)\\[3ex]
		
		14&	\multirow{1}{7cm}{The fourth  general E Banhatti RL exponential} & \(GBRL_2(G,x)=\sum\limits_{uv\in E(G)}x^{\left(  |B(u)-B(v)|B(u)B(v)\right)^a} \)\\[3ex]

	\end{longtable}
	
	\begin{longtable}[c]{| c | l | l |}
		\caption{\small Multiplicative E- Banhatti Rehan- Lanel indices \label{tableb4}}\\
		
		\hline
		\multicolumn{1}{|c}{\textbf{S.N.}} & 
		\multicolumn{1}{|c}{\textbf{Name of the index}} & \multicolumn{1}{|c|}{\textbf{Formula}} \\\hline
		\endfirsthead
		
		\hline
		\endlastfoot

		1&	The third multiplivative E Banhatti RL index & \(BRL_1(G)=\prod\limits_{uv\in E(G)}\left(   B(u)-B(v)+B(u)B(v)\right) \)\\[3ex]
		
		2&	The fourth multiplivative E Banhatti RL index & \(BRL_2(G)=\prod\limits_{uv\in E(G)}\left(  |B(u)-B(v)|B(u)B(v)\right) \)\\[3ex]
		
		3&		\multirow{2}{6cm}{The third  multiplicative hyper Banhatti RL index} & \(HBRL_1(G)=\prod\limits_{uv\in E(G)}\left(   B(u)-B(v)+B(u)B(v)\right)^2 \)\\[5ex]
		\hline 
		4&		\multirow{2}{7cm}{The fourth multiplicative hyper Banhatti RL index} & \(HBRL_2(G)=\prod\limits_{uv\in E(G)}\left( |B(u)-B(v)|B(u)B(v)\right)^2 \)\\[5ex]
		
		5&	\multirow{1}{7cm}{The third multiplicative inverse Banhatti RL index} & \(IBRL_1(G)=\prod\limits_{uv\in E(G)}\dfrac{1}{\left(   B(u)-B(v)+B(u)B(v)\right)} \)\\[3ex]
		6&	\multirow{2}{5cm}{The fourth multiplicative inverse Banhatti RL index} & \(IBRL_2(G)=\prod\limits_{uv\in E(G)}\dfrac{1}{\left(  |B(u)-B(v)|B(u)B(v)\right)} \)\\[6ex]
		
		7&	\multirow{2}{7cm}{The third general multiplicative Banhatti RL index} & \(GBRL_1(G)=\prod\limits_{uv\in E(G)}\left(   B(u)-B(v)+B(u)B(v)\right)^a \)\\[6ex]
		8&	\multirow{2}{7cm}{The fourth  general multiplicative Banhatti RL index} & \(GBRL_2(G)=\prod\limits_{uv\in E(G)}\left(  |B(u)-B(v)|B(u)B(v)\right)^a \)\\[6ex]
		
		9&	\multirow{2}{5cm}{The third multiplicative Banhatti RL exponential} & \(BRL_1(G,x)=\prod\limits_{uv\in E(G)}x^{\left(   B(u)-B(v)+B(u)B(v)\right)} \)\\[3ex]
		
		10&	\multirow{2}{5cm}{The fourth multiplicative Banhatti RL exponential} & \(BRL_2(G,x)=\prod\limits_{uv\in E(G)}x^{\left(  |B(u)-B(v)|B(u)B(v)\right)} \)\\[3ex]
		
		11&	\multirow{2}{5cm}{The third multiplicative hyper Banhatti RL exponential} & \(HBRL_1(G,x)=\prod\limits_{uv\in E(G)}x^{\left(   B(u)-B(v)+B(u)B(v)\right)^2 }\)\\[5ex]
		
		12&	\multirow{2}{5cm}{The fourth multiplicative hyper Banhatti RL exponential} & \(HBRL_2(G,x)=\prod\limits_{uv\in E(G)}x^{\left(  |B(u)-B(v)|B(u)B(v)\right)^2} \)\\[3ex]
		13&	\multirow{2}{5cm}{The third multiplicative inverse Banhatti RL exponential} & \(IBRL_1(G,x )=\prod\limits_{uv\in E(G)}x^{\dfrac{1}{\left(   B(u)-B(v)+B(u)B(v)\right)}} \)\\[3ex]
		
		14&	\multirow{2}{5cm}{The fourth multiplicative inverse Banhatti RL exponential} & \(IBRL_2(G)=\prod\limits_{uv\in E(G)}x^{\dfrac{1}{\left(  |B(u)-B(v)|B(u)B(v)\right)} }\)\\[6ex]
		
		15&	\multirow{2}{5cm}{The third multiplicative general Banhatti RL exponential} & \(GBRL_1(G,x)=\prod\limits_{uv\in E(G)}x^{\left(   B(u)-B(v)+B(u)B(v)\right)^a} \)\\[6ex]
		16&	\multirow{2}{5cm}{The fourth multiplicative general Banhatti RL exponential} & \(GBRL_2(G,x)= \prod\limits_{uv\in E(G)}\displaystyle{x^{\left(  |B(u)-B(v)|B(u)B(v)\right)^a}} \)\\[6ex]
		
		\hline
		
	\end{longtable}
	
	We compute the first and second Banhatti Rehan-Lanel indices \((BRL_1, BRL_2)\) for few standard graph such as \(r-\) regular graph, cycle, complete graph, path and complete bipartite graph as well as wheel graph and sunflower graph. Additionally, we calculate the first and second Banhatti Rehan-Lanel exponential for the wheel graph.
	\begin{proposition}
		Let \(G\) be r-regular with \(|V(G)|= n\) and \(r\geq 2\). Then
		\begin{enumerate}[label=\roman*.]
			\item \(BRL_1(G)=6nr\left( \dfrac{r-1}{n-r}\right) ^2\)
			\item \(BRL_2(G)=2nr\left( \dfrac{r-1}{n-r}\right) ^2\)
			\item \(BRL_3(G)=\dfrac{2nr(r-1)^2 }{(n-r)^2}\)
			\item \(BRL_4(G)=2nr\left( \dfrac{r-1}{n-r}\right) ^2.\)
		\end{enumerate}
	
	\end{proposition}
	
	\begin{proof}
		Let \(G\) is an $r-$regular graph with \(n\) vertices and
		\(r \geq 2\).

		Then, \(G\) has \(\dfrac{nr}{2}\)	edges. For any edge \(uv=e\) in \(G\), \(d_G(e)=
		d_G(u)+ d_G(u)-2=2r-2.\)\\
		
		From definition we have\\
		
		$$	\begin{array}{rl}
			BRL_1(G)=&\displaystyle \sum\limits_{uv\in E(G)}\left[ B^2(u)+B^2(v) +B(u)B(v)\right] \\
			=&\displaystyle\sum\limits_{uv\in E(G)}\left[ \left( \dfrac{2r-2}{n-r} \right)^2 +\left( \dfrac{2r-2}{n-r} \right)^2  +\left( \dfrac{2r-2}{n-r} \right) \left( \dfrac{2r-2}{n-r} \right)\right]  \\
			=& 12\dfrac{nr}{2}\dfrac{(r-1)^2}{(n-r)^2}\\
			=& \dfrac{6nr(r-1)^2 }{(n-r)^2}.
		\end{array}$$
		
		$$	\begin{array}{rl}
			BRL_2(G)=&\displaystyle\sum\limits_{uv\in E(G)}\left[ B^2(u)+B^2(v) -B(u)B(v)\right] \\
			=&\displaystyle\sum\limits_{uv\in E(G)}\left[ \left( \dfrac{2r-2}{n-r} \right)^2 +\left( \dfrac{2r-2}{n-r} \right)^2  -\left( \dfrac{2r-2}{n-r} \right) \left( \dfrac{2r-2}{n-r} \right)\right]  \\
			=& 12\dfrac{nr}{2}\dfrac{(r-1)^2}{(n-r)^2}\\
			=& \dfrac{2nr(r-1)^2 }{(n-r)^2}.
		\end{array}$$
		
			$$	\begin{array}{rl}
			BRL_3(G)=&\displaystyle \sum\limits_{uv\in E(G)}\left[ B(u)-B(v) +B(u)B(v)\right] \\
			=&\displaystyle\sum\limits_{uv\in E(G)}\left[ \left( \dfrac{2r-2}{n-r} \right) -\left( \dfrac{2r-2}{n-r} \right)  +\left( \dfrac{2r-2}{n-r} \right) \left( \dfrac{2r-2}{n-r} \right)\right]  \\
			=& 4\dfrac{nr}{2}\dfrac{(r-1)^2}{(n-r)^2}\\
			=& \dfrac{2nr(r-1)^2 }{(n-r)^2}.
		\end{array}$$
		
		$$	\begin{array}{rl}
			BRL_4(G)=&\displaystyle\sum\limits_{uv\in E(G)} |B(u)-B(v)|B(u)B(v) \\
			=&\displaystyle\sum\limits_{uv\in E(G)}\left| \left( \dfrac{2r-2}{n-r} \right) -\left( \dfrac{2r-2}{n-r} \right)\right|  \left( \dfrac{2r-2}{n-r} \right) \left( \dfrac{2r-2}{n-r} \right)  \\
			=& 0.
			
		\end{array}$$
	\end{proof}
	
	\begin{corollary}
		Let	\(C_n\)
		be a cycle with \(n \geq 3\) vertices. Then,
		\begin{enumerate}[label=\roman*.]
			\item \(BRL_1(C_n)=12n\left( \dfrac{1}{n-2}\right) ^2\)
			\item \(BRL_2(C_n)=4n\left( \dfrac{1}{n-2}\right) ^2\)
			\item  \(BRL_3(C_n)=4n\left( \dfrac{1}{n-2}\right)\)
			\item\(\ BRL_4(C_n)=0\)
		\end{enumerate}

	\end{corollary}

	\begin{corollary}
		Let \(K_n\) be a complete graph with \(n\geq 3\)
		vertices. Then,
		\begin{enumerate}[label=\roman*.]
			\item \(BRL_1(K_n)=6n(n-1)(n-2)^2\)
			\item \(BRL_2(K_n)=2n(n-1)(n-2)^2\)
			\item \(BRL_3(K_n)=2n(n-1)(n-2)^2 \)
			\item \(BRL_4(K_n)=0.\)
		\end{enumerate}
		
	\end{corollary}

	\begin{proposition}
		Let \(P_n\) be a path with \(n\geq 3\) vertices. Then,
		
		\begin{enumerate}[label=\roman*.]
	\item \(BRL_1(P_n)=\dfrac{2(n-1)^2+2(n-2)^2+2(n-1)(n-2)+12(n-1)^2(n-3)}{(n-1)^2(n-2)^2}\)
	\item \(BRL_2(P_n)= \dfrac{2(n-1)^2+2(n-2)^2-2(n-1)(n-2)+4(n-1)^2(n-3)}{(n-1)^2(n-2)^2}\)
	\item \(BRL_3(P_n)= \dfrac{2\left|n^2-6n+10\right|}{(n-1)(n-2)^2}\)
			\item \(BRL_4(P_n)= \dfrac{4|n|}{(n-1)^2(n-2)^2}\)
		\end{enumerate}

	\end{proposition}

	\begin{proof}
		The proof can be done similarly.\\
	\end{proof}
	
	\begin{proposition}
		For \(K_{m,n}\) with \(1 \leq m \leq n\) and \(n \geq 2\),

		\begin{enumerate}[label=\roman*.]
			\item \(BRL_1(K_{m,n})=\dfrac{(m+n-2)^2\left(m^2+n^2+mn\right) }{mn} \)
			\item \(BRL_2(K_{m,n})=\dfrac{(m+n-2)^2 \left(m^2+n^2-mn\right)}{mn}\)
			\item \(BRL_3(K_{m,n})=2(m+n-2)\left(m-1\right) , m>n\)
			\item \(BRL_4(K_{m,n})=\dfrac{(m+n-2)^4 \left|m-n\right|}{mn}, m>n\).
		\end{enumerate}
		
	\end{proposition}
	\begin{proof}
		$$	\begin{array}{rl}
			BRL_2(K_{m,n})=&\sum\limits_{uv\in E(G)}\left[ B^2(u)+B^2(v) +B(u)B(v)\right] \\
			=&\sum\limits_{uv\in E(G)}\left[ \left( \dfrac{m+n-2}{n} \right)^2 +\left( \dfrac{m+n-2}{m} \right)^2  +\left( \dfrac{m+n-2}{n} \right) \left( \dfrac{m+n-2}{m} \right)\right]  \\
			=& mn\dfrac{(m+n-2)^2}{m^2n^2} \left[ m^2+n^2+mn\right]\\
			=& \dfrac{(m+n-2)^2 \left(m^2+n^2+mn\right) }{mn}. 
		\end{array}$$

		$$	\begin{array}{rl}
			BRL_2(K_{m,n})=&\sum\limits_{uv\in E(G)}\left[ B^2(u)+B^2(v) -B(u)B(v)\right] \\
			=&\sum\limits_{uv\in E(G)}\left[ \left( \dfrac{m+n-2}{n} \right)^2 +\left( \dfrac{m+n-2}{m} \right)^2  -\left( \dfrac{m+n-2}{n} \right) \left( \dfrac{m+n-2}{m} \right)\right]  \\
			=& mn\dfrac{(m+n-2)^2}{m^2n^2} \left[ m^2+n^2-mn\right]\\
			=& \dfrac{(m+n-2)^2 \left(m^2+n^2-mn\right) }{mn}.
		\end{array}$$
		
			$$	\begin{array}{rl}
			BRL_3(K_{m,n})=&\sum\limits_{uv\in E(G)}\left[ B(u)-B(v) +B(u)B(v)\right] \\
			=&\sum\limits_{uv\in E(G)}\left[ \left( \dfrac{m+n-2}{n} \right)-\left( \dfrac{m+n-2}{m} \right)  +\left( \dfrac{m+n-2}{n} \right) \left( \dfrac{m+n-2}{m} \right)\right]  \\
			=& mn\dfrac{(m+n-2)}{mn} \left[ m-n+m+n-2\right]\\
			=& 2(m+n-2) \left(m-1\right) . 
		\end{array}$$

		$$	\begin{array}{rl}
			BRL_4(K_{m,n})=&\sum\limits_{uv\in E(G)}\left| B(u)-B(v)\right|B(u)B(v) \\
			=&\sum\limits_{uv\in E(G)}\left| \left( \dfrac{m+n-2}{n} \right) -\left( \dfrac{m+n-2}{m} \right) \right| \left( \dfrac{m+n-2}{n} \right) \left( \dfrac{m+n-2}{m} \right)  \\
			=& mn\dfrac{(m+n-2)^3\left| m-n\right|}{m^2n^2} \\
			=& \dfrac{(m+n-2)^3 \left|m-n\right| }{mn}.
		\end{array}$$
		
	\end{proof}
	\begin{corollary}
		For \(K_{n,n}\) with \(n \geq 2\),
		
		\begin{enumerate}[label=\roman*.]
			\item \(BRL_1(K_{n,n})=12(n-1)^2 \)
			\item \(BRL_2(K_{n,n})=(n-1)^2\)
			\item \(BRL_3(K_{n,n})=4(n-1)^2 \)
			\item  \(BRL_4(K_{n,n})=0\).
		\end{enumerate}
	\end{corollary}

	\begin{corollary}
		For \(K_{1,n}\) with \(n \geq 2\),
		\begin{enumerate}[label=\roman*.]
			\item \(BRL_1(K_{1,n})=\dfrac{(n-1)^2\left(n^2+n+1\right) }{n} \)
			\item \(BRL_2(K_{1,n})=\dfrac{(n-1)^2 \left(n^2+1-n\right) }{n}\)
			\item \(BRL_3(K_{1,n})=2(n-1)^2 \)
			\item\(BRL_4(K_{1,n})=\dfrac{(n-1)^3\left|1-n\right| }{n}\).
		\end{enumerate}
		
	\end{corollary}

	We compute the first, second, third and fourth Banhatti Rehan-Lanel indices of the
	wheel graph \(W_n\).
	\begin{theorem}
		
		Let \(W_n\) be a wheel graph with \(n+1\) vertices and
		\(2n\) edges. Then, 
		\begin{enumerate}[label=\roman*.]
			\item \(BRL_1(W_n)=\dfrac{n}{(n-2)^2}((n+1)^2(n^2-3n+3)+48)\).	
			\item \(BRL_2(W_n)=\dfrac{n}{(n-2)^2}((n+1)^2(n^2-5n+7)+16)\)
		\item\(BRL_3(W_n)=\dfrac{n}{(n-2)^2}\left[(n+1)(n-2)-(n+1)(n-2)^2+(n+1)^2+16\right]\)
			\item \(BRL_4(W_n)=\dfrac{n|n-1|(n+1)^3}{(n-2)}\)
		\end{enumerate}

	\end{theorem}

	\begin{proof}
		
		There are two distinct types of Banhatti edges that can be defined in the set \(W_n\), and these are determined by the Banhatti degrees of the end vertices of each edge.
		$$\begin{array}{lll}
			BE_1=&\left\lbrace uv \in E(W_n)| \ B(u)=B(v)=\dfrac{4}{n-2} \right\rbrace, &|E_1|=n.   \\
			BE_2=&\left\lbrace uv \in E(W_n)| \ B(u)=\dfrac{n+1}{n-2} ,B(v)=n+1 \right\rbrace ,&|E_2|=n. 
		\end{array}$$

		From definition and by cardinalities of the Banhatti
		edge partition of \(W_n\), we obtain
		$$	\begin{array}{rl}
			BRL_1(W_n)=&\sum\limits_{uv\in E(G)}\left[ B^2(u)+B^2(v) +B(u)B(v)\right] \\
			=&\sum\limits_{uv\in E(G)}\left[ \left(\dfrac{4}{n-2} \right)^2 +\left(\dfrac{4}{n-2} \right)^2  +\left(\dfrac{4}{n-2} \right) \left(\dfrac{4}{n-2} \right)\right] +\sum\limits_{uv\in E(G)}\left[ \left(\dfrac{n+1}{n-2} \right)^2 +\left(n+1 \right)^2  +\left(\dfrac{n+1}{n-2} \right) \left(n+1 \right)\right] \\
			=& n\times 3\left(\dfrac{4}{n-2} \right)^2 +n\left[ \left(\dfrac{n+1}{n-2} \right)^2 +\left(n+1 \right)^2  +\left(\dfrac{n+1}{n-2} \right) \left(n+1 \right)\right]\\
			=& \dfrac{n}{(n-2)^2}\left[ (n+1)^2(n^2-3n+3)+48\right]. 
		\end{array}$$
		
		$$	\begin{array}{rl}
			BRL_2(W_n)=&\sum\limits_{uv\in E(G)}\left[ B^2(u)+B^2(v) -B(u)B(v)\right] \\
			=&\sum\limits_{uv\in E(G)}\left[ \left(\dfrac{4}{n-2} \right)^2 +\left(\dfrac{4}{n-2} \right)^2  -\left(\dfrac{4}{n-2} \right) \left(\dfrac{4}{n-2} \right)\right] +\sum\limits_{uv\in E(G)}\left[ \left(\dfrac{n+1}{n-2} \right)^2 +\left(n+1 \right)^2  -\left(\dfrac{n+1}{n-2} \right) \left(n+1 \right)\right] \\
			=& n\left(\dfrac{4}{n-2} \right)^2 +n\left[ \left(\dfrac{n+1}{n-2} \right)^2 +\left(n+1 \right)^2  -\left(\dfrac{n+1}{n-2} \right) \left(n+1 \right)\right]\\
			=& \dfrac{n}{(n-2)^2}\left[ (n+1)^2(n^2-5n+7)+16\right].
		\end{array}$$
		
		$$	\begin{array}{rl}
			BRL_3(W_n)=&\sum\limits_{uv\in E(G)}\left[ B(u)-B(v) +B(u)B(v)\right] \\
			=&\sum\limits_{uv\in E(G)}\left[ \left(\dfrac{4}{n-2} \right) -\left(\dfrac{4}{n-2} \right) +\left(\dfrac{4}{n-2} \right) \left(\dfrac{4}{n-2} \right)\right] +\sum\limits_{uv\in E(G)}\left[ \left(\dfrac{n+1}{n-2} \right) -\left(n+1 \right)  +\left(\dfrac{n+1}{n-2} \right) \left(n+1 \right)\right] \\
			=& n\left(\dfrac{4}{n-2} \right)^2 +n\left[ \left(\dfrac{n+1}{n-2} \right)(n-2)^2 -\left(n+1 \right)(n-2)^2  +\left(\dfrac{n+1}{n-2} \right) \left(n+1 \right)(n-2)\right]\\
			=& \dfrac{n}{(n-2)^2}\left[ (n+1)(n-2)-(n+1)(n-2)^2+(n+1)^2+16\right]. 
		\end{array}$$
		
		$$	\begin{array}{rl}
			BRL_4(W_n)=&\sum\limits_{uv\in E(G)}\left| B(u)-B(v)\right| B(u)B(v) \\
			=&\sum\limits_{uv\in E(G)}\left| \left(\dfrac{4}{n-2} \right) -\left(\dfrac{4}{n-2} \right) \right|\left(\dfrac{4}{n-2} \right) \left(\dfrac{4}{n-2} \right) +\sum\limits_{uv\in E(G)}\left| \left(\dfrac{n+1}{n-2} \right) -\left(n+1 \right)\right| \left(\dfrac{n+1}{n-2} \right) \left(n+1 \right) \\
			=& n\left(0 \right)+n\left| \left(\dfrac{n+1}{n-2} \right) -\left(n+1 \right)\right|\left(\dfrac{n+1}{n-2} \right) \left(n+1 \right)\\
			=& \dfrac{n|n-1|(n+1)^3}{(n-2)} .
		\end{array}$$
	\end{proof}
	We determine the first, second, third and fourth E-Banhatti RL polynomials of \(W_n\) by applying definitions and cardinalities of the Banhatti edge partition.
	\begin{theorem}
		Let \(W_n\) be a wheel graph. Then,
		
		\begin{enumerate}[label=\roman*.]
			\item \(BRL_1(W_n,x)=nx^{\left[ \dfrac{48}{(n-2)^2}\right]}+nx^{\left[ \dfrac{(n+1)^2(n^2-3n+3)}{(n-2)^2}\right]}\)
			\item \(BRL_2(W_n,x)=nx^{\left[ \dfrac{16}{(n-2)^2}\right]}+nx^{\left[ \dfrac{(n+1)^2(n^2-5n+7)}{(n-2)^2}\right]}\)	
			\item \(BRL_3(W_n,x)=nx^{\left[ \dfrac{16}{(n-2)^2}\right]}+nx^{\left[ \dfrac{4(n+1)}{(n-2)}\right]}\)
			\item \(BRL_4(W_n,x)=n+nx^{\left[ \dfrac{(n+1)^3|n-3|}{(n-2)^2}\right]}\).
		\end{enumerate}
		
	\end{theorem}

	\begin{proof}
		$$\begin{array}{ll}
			BRL(W_n,x)&= nx^{\left[ \left(\dfrac{4}{n-2} \right)^2 +\left(\dfrac{4}{n-2} \right)^2  +\left(\dfrac{4}{n-2} \right) \left(\dfrac{4}{n-2} \right)\right]}+nx^{\left[ \left(\dfrac{n+1}{n-2} \right)^2 +\displaystyle \left(n+1 \right)^2  +\left(\dfrac{n+1}{n-2} \right) \left(n+1 \right)\right]}.\\
			&= nx^{\left[ \dfrac{48}{(n-2)^2}\right]}+nx^{\left[ \dfrac{(n+1)^2(n^2-3n+3)}{(n-2)^2}\right].}
		\end{array}$$
		
		$$\begin{array}{ll}
			BRL(W_n,x)&= nx^{\left[ \left(\dfrac{4}{n-2} \right)^2 +\left(\dfrac{4}{n-2} \right)^2  -\left(\dfrac{4}{n-2} \right) \left(\dfrac{4}{n-2} \right)\right]}+nx^{\left[ \left(\dfrac{n+1}{n-2} \right)^2 +\displaystyle \left(n+1 \right)^2  -\left(\dfrac{n+1}{n-2} \right) \left(n+1 \right)\right]}\\
			&= nx^{\left[ \dfrac{16}{(n-2)^2}\right]}+nx^{\left[ \dfrac{(n+1)^2(n^2-5n+7)}{(n-2)^2}\right].}
		\end{array}$$
		
		 $$\begin{array}{ll}
			BRL_3(W_n,x)&= nx^{\left[ \left(\dfrac{4}{n-2} \right) -\left(\dfrac{4}{n-2} \right)  +\left(\dfrac{4}{n-2} \right) \left(\dfrac{4}{n-2} \right)\right]}+nx^{\left[ \left(\dfrac{n+1}{n-2} \right) -\displaystyle \left(n+1 \right)  +\left(\dfrac{n+1}{n-2} \right) \left(n+1 \right)\right]}.\\
			&= nx^{\left[ \dfrac{16}{(n-2)^2}\right]}+nx^{\left[ \dfrac{4(n+1)}{(n-2)}\right].}
		\end{array}$$
		
		$$\begin{array}{ll}
			BRL_4(W_n,x)&= nx^{\left|\left(\dfrac{4}{n-2} \right) -\left(\dfrac{4}{n-2} \right)\right|  \left(\dfrac{4}{n-2} \right) \left(\dfrac{4}{n-2} \right)}+nx^{\left|\left(\dfrac{n+1}{n-2} \right) -\displaystyle \left(n+1 \right)\right|\left(\dfrac{n+1}{n-2} \right) \left(n+1 \right)}\\
			&= n+nx^{\left[ \dfrac{(n+1)^3|n-3|}{(n-2)^2}\right]}.
		\end{array}$$
	\end{proof}
	
	\begin{theorem}
		Let \(S\!f_n\) be a sunflower graph with \(3n+1\) vertices and \(5n\) edges. Then,
		
		\begin{enumerate}[label=\roman*.]
			\item \(\displaystyle BRL_1(S\!f_n)=\frac{12n}{(n-1)^2}+n(3n+2)^2\left( \frac{(3n-3)^2+(3n+2)(3n-3)+1}{(3n-3)^2}\right) \\ \ +16n\left(\dfrac{1}{(3n-1)^2}+\dfrac{1}{(3n-3)^2}+\dfrac{1}{(3n-1)(3n-3)} \right) + 3n^3\left( 1+\dfrac{1}{(3n-1)^2}+\dfrac{1}{3n-1}\right)+n(3n-1)^2\left(\dfrac{1}{9n^2}+1+\dfrac{1}{3n} \right)\)
			\item \(\displaystyle BRL_1(S\!f_n)=\frac{4n}{(n-1)^2}+n(3n+2)^2\left( \frac{(3n-3)^2-(3n+2)(3n-3)+1}{(3n-3)^2}\right) \\ \ +16n\left(\dfrac{1}{(3n-1)^2}+\dfrac{1}{(3n-3)^2}-\dfrac{1}{(3n-1)(3n-3)} \right) + 3n^3\left( 1+\dfrac{1}{(3n-1)^2}-\dfrac{1}{3n-1}\right)+n(3n-1)^2\left(\dfrac{1}{9n^2}+1-\dfrac{1}{3n} \right)\)
			\item \(\displaystyle BRL_3(S\!f_n)=\dfrac{4n}{(n-1)^2}+\left( \dfrac{3n^2(3n+2)}{(3n-3)}\right) +\left(\dfrac{8n}{(3n-1)(3n-3)} \right) + \left( \dfrac{18n^2}{(3n-1)}\right)\)
			\item \(\displaystyle BRL_4(S\!f_n)=4n(3n+2)^2\left( \dfrac{|3n-4|}{(3n-3)^2}\right) +\left(\dfrac{128n}{(3n-1)^2(3n-3)^2} \right) + 27n^4\left( \dfrac{|3n-2|}{(3n-1)^2}\right)+n(3n-1)^3\left(\dfrac{|1-3n|}{9n^2} \right)\).
		\end{enumerate}

	\end{theorem}
	
	\begin{proof}
		In \(S\!f_n\), five types of Banhatti edges can be defined based on Banhatti degrees of end vertices of each edge follow:
		$$\begin{array}{lll}
			BE_1&=\left\lbrace uv \in E(Sf_n)| B(u)=B(v)=\dfrac{6}{3n-3}=\dfrac{2}{n-1}\right\rbrace, &|E_1|=n.\\
			BE_2&=\left\lbrace uv \in E(Sf_n)| B(v_0)=3n+2 ,B(v)=\dfrac{3n+2}{3n-3} \right\rbrace ,& |E_2|=n.\\
			BE_3&=\left\lbrace uv \in E(Sf_n)| B(u)=\dfrac{4}{3n-1}, B(v)=\dfrac{4}{3n-3} \right\rbrace, &|E_3|=n.\\
			BE_4&=\left\lbrace uv \in E(Sf_n)| B(v_0)=3n ,B(v)=\dfrac{3n}{3n-1} \right\rbrace , &|E_4|=n.\\
			BE_5&=\left\lbrace uv \in E(Sf_n)| B(w)=\dfrac{3n-1}{3n} ,B(v_0)=3n-1 \right\rbrace , &|E_5|=n.
			
		\end{array}$$
		$$\begin{array}{ll}
			BRL_1(Sf_n)&=\sum\limits_{uv\in E(Sf_n)}\left[ B^2(u)+B^2(v) + B(u)B(v)\right]\\
			&=n\left[\left( \dfrac{2}{n-1}\right)^2+\left( \dfrac{2}{n-1}\right)^2+\left( \dfrac{2}{n-1}\right)\left( \dfrac{2}{n-1}\right)   \right]+n\left[\left( 3n+2\right)^2+\left( \dfrac{3n+2}{3n-3}\right)^2+\left( 3n+2\right)\left( \dfrac{4}{3n-3}\right)   \right] \\ &+ n\left[\left( \dfrac{4}{3n-1}\right)^2+\left( \dfrac{4}{3n-3}\right)^2+\left( \dfrac{4}{3n-1}\right)\left( \dfrac{4}{3n-3}\right)   \right]+n\left[\left( 3n\right)^2+\left( \dfrac{3n}{3n-1}\right)^2+\left( 3n\right)\left( \dfrac{3n}{3n-1}\right)   \right]\\ &+  n\left[\left( \dfrac{3n-1}{3n}\right)^2+\left( 3n-1\right)^2+\left( \dfrac{3n-1}{3n}\right)\left( 3n-1\right)   \right]\\
			&=n\dfrac{12}{(n-1)^2}+n(3n+2)^2\left( \dfrac{(3n-3)^2+(3n+2)(3n-3)+1}{(3n-3)^2}\right) \\&+16n\left(\dfrac{1}{(3n-1)^2}+\dfrac{1}{(3n-3)^2}+\dfrac{1}{(3n-1)(3n-3)} \right) + 3n^3\left( 1+\dfrac{1}{(3n-1)^2}+\dfrac{1}{3n-1}\right)\\&+n(3n-1)^2\left(\dfrac{1}{9n^2}+1+\dfrac{1}{3n} \right).  
		\end{array}$$
		
		$$\begin{array}{ll}
			BRL_2(S\!f_n)&=\sum\limits_{uv\in E(Sf_n)}\left[ B^2(u)+B^2(v) - B(u)B(v)\right]\\
			&=n\left[\left( \dfrac{2}{n-1}\right)^2+\left( \dfrac{2}{n-1}\right)^2-\left( \dfrac{2}{n-1}\right)\left( \dfrac{2}{n-1}\right)   \right]+n\left[\left( 3n+2\right)^2+\left( \dfrac{3n+2}{3n-3}\right)^2-\left( 3n+2\right)\left( \dfrac{4}{3n-3}\right)   \right] \\ &+ n\left[\left( \dfrac{4}{3n-1}\right)^2+\left( \dfrac{4}{3n-3}\right)^2-\left( \dfrac{4}{3n-1}\right)\left( \dfrac{4}{3n-3}\right)   \right]+n\left[\left( 3n\right)^2+\left( \dfrac{3n}{3n-1}\right)^2-\left( 3n\right)\left( \dfrac{3n}{3n-1}\right)   \right]\\ &+  n\left[\left( \dfrac{3n-1}{3n}\right)^2+\left( 3n-1\right)^2-\left( \dfrac{3n-1}{3n}\right)\left( 3n-1\right)   \right]\\
			&=\dfrac{4n}{(n-1)^2}+n(3n+2)^2\left( \dfrac{(3n-3)^2-(3n+2)(3n-3)+1}{(3n-3)^2}\right) \\&+16n\left(\dfrac{1}{(3n-1)^2}+\dfrac{1}{(3n-3)^2}-\dfrac{1}{(3n-1)(3n-3)} \right) + 3n^3\left( 1+\dfrac{1}{(3n-1)^2}-\dfrac{1}{3n-1}\right)\\&+n(3n-1)^2\left(\dfrac{1}{9n^2}+1-\dfrac{1}{3n} \right).  
		\end{array}$$
		
		$$\begin{array}{ll}
			BRL_3(Sf_n)&=\sum\limits_{uv\in E(Sf_n)}\left[ B(u)-B(v) + B(u)B(v)\right]\\
			&=n\left[\left( \dfrac{2}{n-1}\right)-\left( \dfrac{2}{n-1}\right)+\left( \dfrac{2}{n-1}\right)\left( \dfrac{2}{n-1}\right)   \right]+n\left[\left( 3n+2\right)-\left( \dfrac{3n+2}{3n-3}\right)+\left( 3n+2\right)\left( \dfrac{4}{3n-3}\right)   \right] \\ &+ n\left[\left( \dfrac{4}{3n-1}\right)-\left( \dfrac{4}{3n-3}\right)+\left( \dfrac{4}{3n-1}\right)\left( \dfrac{4}{3n-3}\right)   \right]+n\left[\left( 3n\right)-\left( \dfrac{3n}{3n-1}\right)+\left( 3n\right)\left( \dfrac{3n}{3n-1}\right)   \right]\\ &+  n\left[\left( \dfrac{3n-1}{3n}\right)-\left( 3n-1\right)+\left( \dfrac{3n-1}{3n}\right)\left( 3n-1\right)   \right]\\
			&=\dfrac{4n}{(n-1)^2}+\left( \dfrac{3n^2(3n+2)}{(3n-3)}\right) +\left(\dfrac{8n}{(3n-1)(3n-3)} \right) + \left( \dfrac{18n^2}{(3n-1)}\right). 
		\end{array}$$
		
		$$\begin{array}{ll}
			BRL_4(S\!f_n)&=\sum\limits_{uv\in E(Sf_n)}\left| B(u)-B(v)\right| B(u)B(v)\\
			&=n\left|\left( \dfrac{2}{n-1}\right)-\left( \dfrac{2}{n-1}\right)\right|\left( \dfrac{2}{n-1}\right)\left( \dfrac{2}{n-1}\right)   +n\left|\left( 3n+2\right)-\left( \dfrac{3n+2}{3n-3}\right)\right|\left( 3n+2\right)\left( \dfrac{4}{3n-3}\right)    \\ &+ n\left|\left( \dfrac{4}{3n-1}\right)-\left( \dfrac{4}{3n-3}\right)\right|\left( \dfrac{4}{3n-1}\right)\left( \dfrac{4}{3n-3}\right)   +n\left|\left( 3n\right)-\left( \dfrac{3n}{3n-1}\right)\right|\left( 3n\right)\left( \dfrac{3n}{3n-1}\right)   \\ &+  n\left|\left( \dfrac{3n-1}{3n}\right)-\left( 3n-1\right)\right|\left( \dfrac{3n-1}{3n}\right)\left( 3n-1\right)   \\
			&=4n(3n+2)^2\left( \dfrac{|3n-4|}{(3n-3)^2}\right) +\left(\dfrac{128n}{(3n-1)^2(3n-3)^2} \right) + 27n^4\left( \dfrac{|3n-2|}{(3n-1)^2}\right)+n(3n-1)^3\left(\dfrac{|1-3n|}{9n^2} \right).  
		\end{array}$$
	\end{proof}
	
	Similarly, the results for other indices/exponentials listed in the above four tables can be obatined.
	\subsection{Revan RL indices}
	Let's say that $G$ is a connected graph that is finite, simple, and has vertex and edge sets denoted by $V(G)$ and \(E(G)\). $d_G(u)$ is the degree of a vertex $u$, and it is determined by the number of vertices that are adjacent to it. In the context of the graph $G$, let \(\Delta (G)\) indicate the greatest (lowest) degree among the vertices of $G$ correspondingly. In the graph $G$, the Revan vertex degree of \(u\) is defined as the expression \(r_G(u)= \Delta(G)+\delta(G)-d_G(u)\). Within this context, the Revan edge that connects the Revan vertices $u$ and $v$ is denoted by the symbol $uv$.
	

	
	Kulli introduced the first and second Revan indices based on the Zagreb indices in \cite{article35} as follows.
	
	\[R_1(G)=\sum_{uv\in E(G)}^{}r_G(u)+r_G(v) \ \text{and} \ R_2(G)=\sum_{uv\in E(G)}^{}r_G(u)r_G(v).  \]
	
	More details and applications of Revan indices can be found in \cite{article64,article62,  article58,article59,article60,article63, article61}.

	Building upon these indices, we present the first, second, third and fourth Rehan-Lanel indices of a graph \(G\) as follows.
	
	\begin{definition}
		The first and second Revan Rehan-Lanel indices of a graph \(G\) are defined by,
		
		\begin{enumerate}[label=\roman*.]
			\item \(RRL_1(G)=\sum_{uv\in E(G)}r_G^2(u)+r_G^2(v)+r_G(u)r_G(v)\)
			\item \(RRL_2(G)=\sum_{uv\in E(G)}r_G^2(u)+r_G^2(v)-r_G(u)r_G(v)\)
			\item \(RRL_3(G)=\sum_{uv\in E(G)}r_G(u)-r_G(v)+r_G(u)r_G(v)\)
			\item \(RRL_4(G)=\sum_{uv\in E(G)}|r_G(u)-r_G(v)|r_G(u)r_G(v).\)
		\end{enumerate}

	\end{definition}

	Furthermore, we introduce an additional set of indices utilising the Revan degree and Rehan-Lanel indices as outlined in the Table \ref{tablec1}, \ref{tablec2}, \ref{tablec3} and \ref{tablec4}.
	
	\begin{longtable}[c]{| c | l | l |}
		\caption{Revan Rehan- Lanel  indices/exponentials \label{tablec1}}\\
		
		\hline
		\multicolumn{1}{|c}{\textbf{S.N.}} & 
		\multicolumn{1}{|c}{\textbf{Name of the index}} & \multicolumn{1}{|c|}{\textbf{Formula}} \\\hline
		\endfirsthead
		\hline
		\endlastfoot

		1&	The first hyper Revan RL index & \(HRRL_1(G)=\sum\limits_{uv\in E(G)}\left(  r_G^2(u)+r_G^2(v)+r_G(u)r_G(v))\right)^2\)\\[3ex]
		
		2&	The second  hyper Revan RL index & \(HRRL_2(G)=\sum\limits_{uv\in E(G)}\left(  r_G^2(u)+r_G^2(v)-r_G(u)r_G(v))\right)^2 \)\\[3ex]
		
		3&	\multirow{1}{5cm}{The first inverse Revan RL index} & \(IRRL_1(G)=\sum\limits_{uv\in E(G)}\dfrac{1}{\left( r_G^2(u)+r_G^2(v)+r_G(u)r_G(v)\right)} \)\\[3ex]
		4&	\multirow{1}{6cm}{The second inverse Revan RL index} & \(IRRL_2(G)=\sum\limits_{uv\in E(G)}\dfrac{1}{\left( r_G^2(u)+r_G^2(v)-r_G(u)r_G(v)\right)} \)\\[3ex]
		
		5&	\multirow{1}{6cm}{The first general Revan RL index} & \(GRRL_1(G)=\sum\limits_{uv\in E(G)}\left(  r_G^2(u)+r_G^2(v)+r_G(u)r_G(v)\right)^a \)\\[3ex]
		
		6&	\multirow{1}{6cm}{The second  general Revan RL index} & \(GRRL_2(G)=\sum\limits_{uv\in E(G)}\left(  r_G^2(u)+r_G^2(v)-r_G(u)r_G(v)\right)^a \)\\[3ex]
		
		7&	\multirow{1}{6cm}{The first Revan RL exponential} & \(RRL_1(G,x)=\sum\limits_{uv\in E(G)}x^{\left(  r_G^2(u)+r_G^2(v)+r_G(u)r_G(v)\right)} \)\\[3ex]
		8&	\multirow{1}{6cm}{The second  Revan RL exponential} & \(RRL_2(G,x)=\sum\limits_{uv\in E(G)}x^{\left( r_G^2(u)+r_G^2(v)-r_G(u)r_G(v)\right)} \)\\[3ex]
		
		9&	\multirow{1}{6cm}{The first general Revan RL exponential} & \(HRRL_1(G,x)=\sum\limits_{uv\in E(G)}x^{\left( r_G^2(u)+r_G^2(v)+r_G(u)r_G(v)\right)^a }\)\\[3ex]
		
		10&	\multirow{1}{6cm}{The second  general Revan RL exponential} & \(\displaystyle HRRL_2(G,x)=\sum\limits_{uv\in E(G)}x^{\left( r_G^2(u)+r_G^2(v)-r_G(u)r_G(v)\right)^a} \)\\[3ex]
		\hline 
		11&	\multirow{1}{6cm}{The first inverse Revan RL exponential} & \(\displaystyle IRRL_1(G,x)=\sum\limits_{uv\in E(G)}x^{\dfrac{1}{\left(  r_G^2(u)+r_G^2(v)+r_G(u)r_G(v)\right)}} \)\\[3ex]
		
		12&	\multirow{1}{6cm}{The second inverse Revan RL exponential} & \(\displaystyle IRRL_2(G,x)=\sum\limits_{uv\in E(G)}x^{\dfrac{1}{\left(  r_G^2(u)+r_G^2(v)-r_G(u)r_G(v)\right)} }\)\\[3ex]
		
		13&	\multirow{1}{6cm}{The first general Revan RL exponential} & \(GRRL_1(G,x)=\sum\limits_{uv\in E(G)}x^{\left( r_G^2(u)+r_G^2(v)+r_G(u)r_G(v)\right)^a} \)\\[3ex]
		
		14&	\multirow{1}{6cm}{The second  general Revan RL exponential} & \(GRRL_2(G,x)=\sum\limits_{uv\in E(G)}x^{\left(  r_G^2(u)+r_G^2(v)-r_G(u)r_G(v)\right)^a} \)\\[3ex]

	\end{longtable}
	
	\begin{longtable}[c]{| c | l | l |}
		\caption{Multiplicative Revan Rehan- Lanel  indices/exponentials \label{tablec2}}\\
		
		\hline
		\multicolumn{1}{|c}{\textbf{S.N.}} & 
		\multicolumn{1}{|c}{\textbf{Name of the index}} & \multicolumn{1}{|c|}{\textbf{Formula}} \\\hline
		\endfirsthead
		
		\hline
		\endlastfoot

		1&	The first multiplicative Revan RL index & \(MRRL_1(G)=\prod\limits_{uv\in E(G)}\left(  r_G^2(u)+r_G^2(v)+r_G(u)r_G(v)\right) \)\\[3ex]
		
		2&	The second multiplicative Revan RL index & \(MRRL_2(G)=\prod\limits_{uv\in E(G)}\left(  r_G^2(u)+r_G^2(v)-r_G(u)r_G(v)\right) \)\\[3ex]
		
		3&		\multirow{2}{4cm}{The first  multiplicative hyper Revan RL index} & \(MHRRL_1(G)=\prod\limits_{uv\in E(G)}\left( r_G^2(u)+r_G^2(v)+r_G(u)r_G(v)\right)^2 \)\\[5ex]
		4&		\multirow{2}{4cm}{The second multiplicative hyper Revan RL index} & \(MHRRL_2(G)=\prod\limits_{uv\in E(G)}\left( r_G^2(u)+r_G^2(v)-r_G(u)r_G(v)\right)^2 \)\\[5ex]
		
		5&	\multirow{2}{4cm}{The first multiplicative inverse Revan RL index} & \(MIRRL_1(G)=\prod\limits_{uv\in E(G)}\dfrac{1}{\left(  r_G^2(u)+r_G^2(v)+r_G(u)r_G(v)\right)} \)\\[6ex]
		
		6&	\multirow{2}{5cm}{The second multiplicative inverse Revan RL index} & \(MIRRL_2(G)=\prod\limits_{uv\in E(G)}\dfrac{1}{\left(  r_G^2(u)+r_G^2(v)-r_G(u)r_G(v)\right)} \)\\[6ex]
		
		7&	\multirow{2}{5cm}{The first general multiplicative Revan RL index} & \(MGRRL_1(G)=\prod\limits_{uv\in E(G)}\left(  r_G^2(u)+r_G^2(v)+r_G(u)r_G(v)\right)^a \)\\[6ex]
		8&	\multirow{2}{5cm}{\small The second  general multiplicative Revan RL index} & \(MGRRL_2(G)=\prod\limits_{uv\in E(G)}\left(  r_G^2(u)+r_G^2(v)-r_G(u)r_G(v)\right)^a \)\\[6ex]
		
		9&	\multirow{2}{5cm}{\small The first multiplicative Revan RL exponential} & \(MRRL_1(G,x)=\prod\limits_{uv\in E(G)}x^{\left(  r_G^2(u)+r_G^2(v)+r_G(u)r_G(v)\right)} \)\\[3ex]
		
		10&	\multirow{2}{5cm}{\small The second multiplicative Revan RL exponential} & \(MRRL_2(G,x)=\prod\limits_{uv\in E(G)}x^{\left(  r_G^2(u)+r_G^2(v)-r_G(u)r_G(v)\right)} \)\\[3ex]
		
		11&	\multirow{2}{5cm}{\small The first multiplicative hyper Revan RL exponential} & \(MGRL_1(G,x)=\prod\limits_{uv\in E(G)}x^{\left(  r_G^2(u)+r_G^2(v)+r_G(u)r_G(v)\right)^2 }\)\\[5ex]
		
		12&	\multirow{2}{5cm}{\small The second multiplicative hyper Revan RL exponential} & \(MHRRL_2(G,x)=\prod\limits_{uv\in E(G)}x^{\left(  r_G^2(u)+r_G^2(v)-r_G(u)r_G(v)\right)^2} \)\\[3ex]
		\hline 
		13&	\multirow{2}{5cm}{The first multiplicative inverse Revan RL exponential} & \(\displaystyle MIRRL_1(G,x )=\prod\limits_{uv\in E(G)}x^{\dfrac{1}{\left(  r_G^2(u)+r_G^2(v)+r_G(u)r_G(v)\right)}} \)\\[3ex]
		
		14&	\multirow{2}{5cm}{The second multiplicative inverse Revan RL exponential} & \(\displaystyle MIRRL_2(G,x)=\prod\limits_{uv\in E(G)}x^{\dfrac{1}{\left(  r_G^2(u)+r_G^2(v)-r_G(u)r_G(v)\right)} }\)\\[5ex]
		
		15&	\multirow{2}{5cm}{The first multiplicative general Revan RL exponential} & \(\displaystyle MGRRL_1(G,x)=\prod\limits_{uv\in E(G)}x^{\left(  r_G^2(u)+r_G^2(v)+r_G(u)r_G(v)\right)^a} \)\\[5ex]
		
		16&	\multirow{2}{5cm}{The second multiplicative general Revan RL exponential} & \(\displaystyle MGRRL_2(G,x)=\prod\limits_{uv\in E(G)}x^{\left(  r_G^2(u)+r_G^2(v)-r_G(u)r_G(v)\right)^a} \)\\[6ex]
		
		\hline
		
	\end{longtable}
	
	\begin{longtable}[c]{| c | l | l |}
		\caption{Revan Rehan- Lanel  indices/exponentials \label{tablec3}}\\
		
		\hline
		\multicolumn{1}{|c}{\textbf{S.N.}} & 
		\multicolumn{1}{|c}{\textbf{Name of the index}} & \multicolumn{1}{|c|}{\textbf{Formula}} \\\hline
		\endfirsthead
		\hline
		\endlastfoot

		1&	The third hyper Revan RL index & \(HRRL_1(G)=\sum\limits_{uv\in E(G)}\left(  r_G(u)-r_G(v)+r_G(u)r_G(v)\right)^2\)\\[3ex]
		
		2&	The second  hyper Revan RL index & \(HRRL_2(G)=\sum\limits_{uv\in E(G)}\left(  \left|r_G(u)-r_G(v)\right|r_G(u)r_G(v)\right)^2 \)\\[3ex]
		
		3&	\multirow{1}{5cm}{The third inverse Revan RL index} & \(IRRL_1(G)=\sum\limits_{uv\in E(G)}\dfrac{1}{\left( r_G(u)-r_G(v)+r_G(u)r_G(v)\right)} \)\\[3ex]
		
		4&	\multirow{1}{6cm}{The second inverse Revan RL index} & \(IRRL_2(G)=\sum\limits_{uv\in E(G)}\dfrac{1}{\left( \left|r_G(u)-r_G(v)\right|r_G(u)r_G(v)\right)} \)\\[3ex]
		
		5&	\multirow{1}{6cm}{The third general Revan RL index} & \(GRRL_1(G)=\sum\limits_{uv\in E(G)}\left(  r_G(u)-r_G(v)+r_G(u)r_G(v)\right)^a \)\\[3ex]
		6&	\multirow{1}{6cm}{The second  general Revan RL index} & \(GRRL_2(G)=\sum\limits_{uv\in E(G)}\left(  \left|r_G(u)-r_G(v)\right|r_G(u)r_G(v)\right)^a \)\\[3ex]
		
		7&	\multirow{1}{6cm}{The third Revan RL exponential} & \(RRL_1(G,x)=\sum\limits_{uv\in E(G)}x^{\left(  r_G(u)-r_G(v)+r_G(u)r_G(v)\right)} \)\\[3ex]
		8&	\multirow{1}{6cm}{The second  Revan RL exponential} & \(RRL_2(G,x)=\sum\limits_{uv\in E(G)}x^{\left( \left|r_G(u)-r_G(v)\right|r_G(u)r_G(v)\right)} \)\\[3ex]
		
		9&	\multirow{1}{6cm}{The third general Revan RL exponential} & \(HRRL_1(G,x)=\sum\limits_{uv\in E(G)}x^{\left( r_G(u)-r_G(v)+r_G(u)r_G(v)\right)^a }\)\\[3ex]
		
		10&	\multirow{1}{6cm}{The second  general Revan RL exponential} & \(\displaystyle HRRL_2(G,x)=\sum\limits_{uv\in E(G)}x^{\left( \left|r_G(u)-r_G(v)\right|r_G(u)r_G(v)\right)^a} \)\\[3ex]
		
		11&	\multirow{1}{6cm}{The third inverse Revan RL exponential} & \(\displaystyle IRRL_1(G,x)=\sum\limits_{uv\in E(G)}x^{\dfrac{1}{\left(  r_G(u)-r_G(v)+r_G(u)r_G(v)\right)}} \)\\[3ex]
		
		12&	\multirow{1}{6cm}{The second inverse Revan RL exponential} & \(\displaystyle IRRL_2(G,x)=\sum\limits_{uv\in E(G)}x^{\dfrac{1}{\left(  \left|r_G(u)-r_G(v)\right|r_G(u)r_G(v)\right)} }\)\\[3ex]
		
		13&	\multirow{1}{6cm}{The third general Revan RL exponential} & \(GRRL_1(G,x)=\sum\limits_{uv\in E(G)}x^{\left( r_G(u)-r_G(v)+r_G(u)r_G(v)\right)^a} \)\\[3ex]
		
		14&	\multirow{1}{6cm}{The second  general Revan RL exponential} & \(GRRL_2(G,x)=\sum\limits_{uv\in E(G)}x^{\left(  \left|r_G(u)-r_G(v)\right|r_G(u)r_G(v)\right)^a} \)\\[3ex]

	\end{longtable}
	
	\begin{longtable}[c]{| c | l | l |}
		\caption{Multiplicative Revan Rehan- Lanel  indices/exponentials \label{tablec4}}\\
		
		\hline
		\multicolumn{1}{|c}{\textbf{S.N.}} & 
		\multicolumn{1}{|c}{\textbf{Name of the index}} & \multicolumn{1}{|c|}{\textbf{Formula}} \\\hline
		\endfirsthead
		
		\hline
		\endlastfoot

		1&	The third multiplicative Revan RL index & \(MRRL_1(G)=\prod\limits_{uv\in E(G)}\left(  r_G(u)-r_G(v)+r_G(u)r_G(v)\right) \)\\[3ex]
		
		2&	The second multiplicative Revan RL index & \(MRRL_2(G)=\prod\limits_{uv\in E(G)}\left(  \left|r_G(u)-r_G(v)\right|r_G(u)r_G(v)\right) \)\\[3ex]
		
		3&		\multirow{2}{4cm}{The third  multiplicative hyper Revan RL index} & \(MHRRL_1(G)=\prod\limits_{uv\in E(G)}\left( r_G(u)-r_G(v)+r_G(u)r_G(v)\right)^2 \)\\[5ex]
		4&		\multirow{2}{4cm}{The second multiplicative hyper Revan RL index} & \(MHRRL_2(G)=\prod\limits_{uv\in E(G)}\left( \left|r_G(u)-r_G(v)\right|r_G(u)r_G(v)\right)^2 \)\\[5ex]
		
		5&	\multirow{2}{4cm}{The third multiplicative inverse Revan RL index} & \(MIRRL_1(G)=\prod\limits_{uv\in E(G)}\dfrac{1}{\left(  r_G(u)-r_G(v)+r_G(u)r_G(v)\right)} \)\\[6ex]
		
		6&	\multirow{2}{5cm}{The second multiplicative inverse Revan RL index} & \(MIRRL_2(G)=\prod\limits_{uv\in E(G)}\dfrac{1}{\left(  \left|r_G(u)-r_G(v)\right|r_G(u)r_G(v)\right)} \)\\[6ex]
		
		7&	\multirow{2}{5cm}{The third general multiplicative Revan RL index} & \(MGRRL_1(G)=\prod\limits_{uv\in E(G)}\left(  r_G(u)-r_G(v)+r_G(u)r_G(v)\right)^a \)\\[6ex]
		8&	\multirow{2}{5cm}{\small The second  general multiplicative Revan RL index} & \(MGRRL_2(G)=\prod\limits_{uv\in E(G)}\left(  \left|r_G(u)-r_G(v)\right|r_G(u)r_G(v)\right)^a \)\\[6ex]
		
		9&	\multirow{2}{5cm}{\small The third multiplicative Revan RL exponential} & \(MRRL_1(G,x)=\prod\limits_{uv\in E(G)}x^{\left(  r_G(u)-r_G(v)+r_G(u)r_G(v)\right)} \)\\[3ex]
		10&	\multirow{2}{5cm}{\small The second multiplicative Revan RL exponential} & \(MRRL_2(G,x)=\prod\limits_{uv\in E(G)}x^{\left(  \left|r_G(u)-r_G(v)\right|r_G(u)r_G(v)\right)} \)\\[3ex]
		
		11&	\multirow{2}{5cm}{\small The third multiplicative hyper Revan RL exponential} & \(MGRL_1(G,x)=\prod\limits_{uv\in E(G)}x^{\left(  r_G(u)-r_G(v)+r_G(u)r_G(v)\right)^2 }\)\\[5ex]
		
		12&	\multirow{2}{5cm}{\small The second multiplicative hyper Revan RL exponential} & \(MHRRL_2(G,x)=\prod\limits_{uv\in E(G)}x^{\left(  \left|r_G(u)-r_G(v)\right|r_G(u)r_G(v)\right)^2} \)\\[3ex]
		
		13&	\multirow{2}{5cm}{The third multiplicative inverse Revan RL exponential} & \(\displaystyle MIRRL_1(G,x )=\prod\limits_{uv\in E(G)}x^{\dfrac{1}{\left(  r_G(u)-r_G(v)+r_G(u)r_G(v)\right)}} \)\\[3ex]
		
		14&	\multirow{2}{5cm}{The second multiplicative inverse Revan RL exponential} & \(\displaystyle MIRRL_2(G,x)=\prod\limits_{uv\in E(G)}x^{\dfrac{1}{\left(  \left|r_G(u)-r_G(v)\right|r_G(u)r_G(v)\right)} }\)\\[5ex]
		
		15&	\multirow{2}{5cm}{The third multiplicative general Revan RL exponential} & \(\displaystyle MGRRL_1(G,x)=\prod\limits_{uv\in E(G)}x^{\left(  r_G(u)-r_G(v)+r_G(u)r_G(v)\right)^a} \)\\[5ex]
		
		16&	\multirow{2}{5cm}{The fourth multiplicative general Revan RL exponential} & \(\displaystyle MGRRL_2(G,x)=\prod\limits_{uv\in E(G)}x^{\left(  \left|r_G(u)-r_G(v)\right|r_G(u)r_G(v)\right)^a} \)\\[6ex]
		
		\hline
		
	\end{longtable}

	Here, the first, second, third and fourth Revan Rehan-Lanel indices \((RRL_1, RRL_2, RRL_3, RRL_4)\) are calculated for \(r-\) regular graph, cycle, a complete graph, a path and a complete bipartite graph. Also, we compute the \(RRL_1, RRL_2, RRL_3\) and \(RRL_4\) for the wheel graph and sunflower graph.
	\begin{proposition}
		If \(G\) is an $r$-regular graph with \(n\) vertices and
		\(r \geq 2\), then
		
		\begin{enumerate}[label=\roman*.]
			\item \(RRL_1(G)=\dfrac{3nr^3}{2}\)
			\item \(RRL_2(G)=\dfrac{nr^3}{2}\)
			\item \(RRL_3(G)=\dfrac{nr^3}{2}\)
			\item\(RRL_4(G)=0.\)
		\end{enumerate} 
	\end{proposition}
	
	\begin{proof}
		Let \(G\) is an $r-$regular graph with \(n\) vertices and
		\(r \geq 2\).\\
		Now	\(d_G(u)=d_G(v)=r\), then \(\Delta(G)=\delta(G)=r\). So. \(r_G(u)=r_G(v)=r+r-r=r.\)
		
		Then \(G\) has \(\dfrac{nr}{2}\)	edges. 
		From definition we have
		$$ \begin{array}{ll}
			RRL_1(G)&=\sum\limits_{uv\in E(G)} (r^2+r^2+r\cdot r)\\
			&=\dfrac{nr}{2} (3r^2)\\
			&=\dfrac{3nr^3}{2}.
		\end{array}$$
		
		$$ \begin{array}{ll}
			RRL_2(G)&=\sum\limits_{uv\in E(G)} (r^2+r^2-r\cdot r)\\
			&=\dfrac{nr}{2} (r^2)\\
			&=\dfrac{nr^3}{2}.
		\end{array}$$ 
		$$ \begin{array}{ll}
			RRL_3(G)&=\sum\limits_{uv\in E(G)} (r-r+r\cdot r)\\
			&=\dfrac{nr}{2} (r^2)\\
			&=\dfrac{nr^3}{2}.
		\end{array}$$
		
		$$ \begin{array}{ll}
			RRL_4(G)&=\sum\limits_{uv\in E(G)} (|r-r|r\cdot r)\\
			&=0.
		\end{array}$$
	\end{proof}

	\begin{corollary}
		Let	\(C_n\)
		be a cycle with \(n \geq 3\) vertices. Then,
		\begin{enumerate}[label=\roman*.]
			\item \(RRL_1(C_n)=12n\)
			\item \(RRL_2(C_n)=4n\)
			\item \(RRL_3(C_n)=4n\)
			\item \(RRL_2(C_n)=0\).
		\end{enumerate}
	
	\end{corollary}

	\begin{corollary}
		Let \(K_n\) be a complete graph with \(n\geq 3\)
		vertices. Then,
		
		\begin{enumerate}[label=\roman*.]
			\item \(RRL_1(K_n)=\dfrac{3n(n-1)^3}{2}\)
			\item \(RRL_2(K_n)=\dfrac{n(n-1)^3}{2}\)
			\item \(RRL_3(K_n)=\dfrac{n(n-1)^3}{2}\)
			\item \(RRL_4(K_n)=0.\)
		\end{enumerate} 
	\end{corollary}
	\begin{proposition}
		Let \(P_n\) be a path with \(n\geq 3\) vertices. Then
		\begin{enumerate}[label=\roman*.]
			\item \(RRL_1(P_n)= 3n+5\)
			\item \(RRL_2(P_n)=n+3\)
			\item \(RRL_3(P_n)=n+3\)
			\item \(RRL_4(P_n)=4\).
		\end{enumerate}
	
	\end{proposition}
	
	\begin{proposition}
		Let $K_{m,n}$ be a complete bipartite graph with
		\(1\leq m\leq n\) and \(n \geq 2\). Then, 
		
		\begin{enumerate}[label=\roman*.]
			\item \(RRL_1(K_{m,n})= mn(m^2+n^2+mn)\)
			\item \(RRL_2(K_{m,n})= mn(m^2+n^2-mn)\)
			\item \(RRL_3(K_{m,n})= mn(n-m+mn)\)
			\item \(RRL_4(K_{m,n})= m^2n^2|n-m|.\)
		\end{enumerate}	
	\end{proposition}

	\begin{theorem}
		Let \(W_n\) be a wheel graph with \(n+1\) vertices and
		\(2n\) edges. Then, 
		\begin{enumerate}[label=\roman*.]
			\item \(RRL_1(W_n)=n(4n^2+3n+9).\)
			\item \(RRL_2(W_n)=n(2n^2-3n+9).\)
			 \item \(RRL_3(W_n)=n(n^2+4n-3)\)
			\item \(RRL_4(W_n)=3n^2|n-3|.\)
		\end{enumerate}

	\end{theorem}
	
	\begin{proof}

		In \(W_n\), two types of Revan edges can be defined based on the Revan degrees of the end vertices of each edge. These types are as follows: 
		$$\begin{array}{lll}
			RE_1&=\left\lbrace uv \in E(W_n)| r_G(v_i)=n\right\rbrace,& |E_1|=n.  \\
			RE_2&=\left\lbrace uv \in E(W_n)| r_G(v)=3, r_G(v_i)=n \right\rbrace &,|E_2|=n. 
		\end{array}$$

		We compute the First and second Revan Rehan-Lanel indices of the
		wheel graph \(W_n\).
		$$\begin{array}{ll}
			RRL_1(W_n)&=\sum\limits_{uv\in E(G)}r_G^2(u)+r_G^2(v)+r_G(u)r_G(v)\\
			&= n(n^2+n^2+n.n)+n(n^2+3^2+n\cdot 3)\\
			&=n(4n^2+3n+9).
		\end{array}$$
		
		$$\begin{array}{ll}
			RRL_2(W_n)&=\sum\limits_{uv\in E(G)}r_G^2(u)+r_G^2(v)-r_G(u)r_G(v)\\
			&= n(n^2+n^2-n.n)+n(n^2+3^2-n\cdot 3)\\
			&=n(2n^2-3n+9).
		\end{array}$$
		
			$$\begin{array}{ll}
			RRL_3(W_n)&=\sum\limits_{uv\in E(G)}r_G(u)-r_G(v)+r_G(u)r_G(v)\\
			&= n(n-n+n.n)+n(n-3+n\cdot 3)\\
			&=n(n^2+4n-3).
		\end{array}$$
		
		$$\begin{array}{ll}
			RRL_4(W_n)&=\sum\limits_{uv\in E(G)}\left|r_G(u)-r_G(v)\right|r_G(u)r_G(v)\\
			&= n(|n-n|n.n)+n(|n-3|n\cdot 3)\\
			&=3n^2|n-3|.
		\end{array}$$
		
	\end{proof}

	\begin{theorem}
		Let \(S\!f_n\) be a sunflower graph with \(3n+1\) vertices and \(5n\) edges.
		
		\begin{enumerate}[label=\roman*.]
			\item \(RRL_1(S\!f_n)=n(54n^2-42n+31).\)
			\item \(RRL_2(S\!f_n)=n(45n^2-36n+27). \)
			\item \(RRL_3(S\!f_n)=n(18n^2+3n+9).\)
			\item \(RRL_4(S\!f_n)=2n(3n-2)|4-3n|+6n^2(3n-2)+6n^2|2-3n|+2n(3n+1)|3n-1|. \)
		\end{enumerate}
	
	\end{theorem}
	
	\begin{proof}
		In \({Sf}_n\), there are five types of Revan edges as follows,
		$$\begin{array}{lll}
			RE_1&=\left\lbrace uv \in E(Sf_n)| r_G(v)=r_G(v)=3n-2 \right\rbrace, &|E_1|=n.\\
			RE_2&=\left\lbrace uv \in E(Sf_n)| r_G(v_0)=2 ,r_G(v)=3n-2 \right\rbrace , &|E_2|=n.\\
			RE_3&=\left\lbrace uv \in E(Sf_n)| r_G(u)=3n, r_G(v)=3n-2 \right\rbrace, &|E_3|=n.\\
			RE_4&=\left\lbrace uv \in E(Sf_n)| r_G(v_0)=2 ,r_G(v)=3n \right\rbrace , &|E_4|=n.\\
			RE_5&=\left\lbrace uv \in E(Sf_n)| r_G(w)=3n+1 ,r_G(v_0)=2 \right\rbrace , &|E_5|=n.
			
		\end{array}$$
		
		$$\begin{array}{ll}
			RRL_1(S\!f_n)&=\sum\limits_{uv\in E(G)}r_G^2(u)+r_G^2(v)+r_G(u)r_G(v)\\
			&=n\left( (3n-2)^2+(3n-2)^2+(3n-2)(3n-2)\right)+n\left( 2^2+(3n-2)^2+2\cdot (3n-2)\right)  \\ &+n\left( (3n)^2+(3n-2)^2+3n(3n-2)\right)+n\left( 2^2+(3n)^2+2\cdot (3n)\right)+n\left((3n+1)^2+2^2+(3n+1)\cdot 2 \right) \\
			&=n(54n^2-42n+31).  
		\end{array}$$
		
		$$\begin{array}{ll}
			RRL_2(S\!f_n)&=\sum\limits_{uv\in E(G)}r_G^2(u)+r_G^2(v)-r_G(u)r_G(v)\\
			&=n\left( (3n-2)^2+(3n-2)^2-(3n-2)(3n-2)\right)+n\left( 2^2+(3n-2)^2-2\cdot (3n-2)\right)  \\ &+n\left( (3n)^2+(3n-2)^2-3n(3n-2)\right)+n\left( 2^2+(3n)^2-2\cdot (3n)\right)+n\left((3n+1)^2+2^2-(3n+1)\cdot 2 \right) \\
			&=n(45n^2-36n+27). 
		\end{array}$$
		
			$$\begin{array}{ll}
			RRL_3(S\!f_n)&=\sum\limits_{uv\in E(G)}r_G(u)-r_G(v)+r_G(u)r_G(v)\\
			&=n\left( (3n-2)-(3n-2)+(3n-2)(3n-2)\right)+n\left( 2-(3n-2)+2\cdot (3n-2)\right)  \\ &+n\left( (3n)-(3n-2)+3n(3n-2)\right)+n\left( 2-(3n)+2\cdot (3n)\right)+n\left((3n+1)-2+(3n+1)\cdot 2 \right) \\
			&=n(18n^2+3n+9).  
		\end{array}$$
		
		$$\begin{array}{ll}
			RRL_4(S\!f_n)&=\sum\limits_{uv\in E(G)}\left|r_G(u)-r_G(v)\right|r_G(u)r_G(v)\\
			&=n\left| (3n-2)-(3n-2)\right|(3n-2)(3n-2)+n\left| 2-(3n-2)\right|2\cdot (3n-2)  \\ &+n\left| (3n)-(3n-2)\right|\cdot 3n(3n-2)+n\left| 2-(3n)\right|\cdot 2\cdot (3n)+n\left|(3n+1)-2\right|(3n+1)\cdot 2  \\
			&=2n(3n-2)|4-3n|+6n^2(3n-2)+6n^2|2-3n|+2n(3n+1)|3n-1|. 
		\end{array}$$
	\end{proof}
	
	Following a similar procedure, we can present the results for the rest of indices/exponentials which are listed in the Table \ref{tablec1}, \ref{tablec2}, \ref{tablec3} and \ref{tablec4}.
	\subsection{Domination RL indices}
	
	Let \(G\) be a graph and \(v \in V(G)\). Then the domination degree of vertex \(v\) is defined as the minimum number of vertices in the Minimal dominating sets containing \(v.\)
	
	
	The minimum domination and maximum domination degree of $G$ are denoted by \(\delta_d(G)=\delta_d\) and \(\Delta_d(G)=\Delta_d\) respectively, where
	\(\delta_d(G) =  \min \left\lbrace {d_d(v) :v \in V (G)}\right\rbrace  \ \text{and} \ \Delta_d = \max \left\lbrace {d_d(v): v \in V (G)}\right\rbrace \).	\\

	In \cite{article42}  VR Kulli introduced the Nirmala index of a molecular graph \(G\) which is defined as, \[N(G)=\sqrt{d_G(u)+d_G(v)}\]

	In \cite{article22} VR Kulli defined the Domination Nirmala index as follows, 
	
	\[DN(G)=\sqrt{d_d(u)+d_d(v)}\]
	
	Inspired by work on Nirmala indices, Domination Nirmala indices, we introduce the first and second Rehan-Lanel domination indices.
	
	\begin{definition}
		The first domination RL index and second  domination RL index of a graph \(G\) are defined as,
		
		\[DRL_1(G)=\sum\limits_{uv\in E(G)}\left(  d_d^2(u)+d_d^2(v)+d_d(u)d_d(v)\right) \ \text{and} \ DRL_2(G)=\sum\limits_{uv\in E(G)}\left(  d_d^2(u)+d_d^2(v)-d_d(u)d_d(v)\right).\]
		
	\end{definition}
	
	Considering the first and second Rehan-Lanel domination indices, we define the few more indices of a graph \(G\) as in table \ref{tabled1}, \ref{tabled2}, \ref{tabled3} and \ref{tabled4}.
	
	\begin{longtable}[c]{| c | l | l |}
		\caption{Rehan- Lanel domination indices/exponentials \label{tabled1}}\\

		\hline
		\multicolumn{1}{|c}{\textbf{S.N.}} & 
		\multicolumn{1}{|c}{\textbf{Name of the index}} & \multicolumn{1}{|c|}{\textbf{Formula}} \\\hline
		\endfirsthead
		
		\hline
		\endlastfoot

		1&	\multirow{2}{6cm}{The first  hyper domination RL index} & \(DRL_1(G)=\sum\limits_{uv\in E(G)}\left(  d_d^2(u)+d_d^2(v)+d_d(u)d_d(v)\right)^2\)\\[3ex]
		2&	\multirow{2}{6cm}{The second  hyper domination RL index} & \(DRL_2(G)=\sum\limits_{uv\in E(G)}\left(  d_d^2(u)+d_d^2(v)-d_d(u)d_d(v)\right)^2\)\\[3ex]
		
		3&	\multirow{2}{6cm}{The first domination RL exponential} & \(DRL_1(G,x)=\sum\limits_{uv\in E(G)}x^{\left(  d_d^2(u)+d_d^2(v)+d_d(u)d_d(v)\right)}\)\\[3ex]
		
		4&	\multirow{2}{6cm}{The second domination RL exponential} & \(DRL_2(G,x)=\sum\limits_{uv\in E(G)}x^{\left(  d_d^2(u)+d_d^2(v)-d_d(u)d_d(v)\right)}\)\\[3ex]

		5&	\multirow{2}{6cm}{The first  general domination RL index} & \(GDRL_1(G)=\sum\limits_{uv\in E(G)}\left(  d_d^2(u)+d_d^2(v)+d_d(u)d_d(v)\right)^a\)\\[3ex]
		
		6&	\multirow{2}{6cm}{The second  general domination RL index} & \(GDRL_2(G)=\sum\limits_{uv\in E(G)}\left(  d_d^2(u)+d_d^2(v)-d_d(u)d_d(v)\right)^a\)\\[3ex]
		7&	\multirow{2}{6cm}{The first hyper domination RL exponential} & \(HDRL_1(G,x)=\sum\limits_{uv\in E(G)}x^{\left(  d_d^2(u)+d_d^2(v)+d_d(u)d_d(v)\right)^2}\)\\[3ex]
		
		8&	\multirow{2}{6cm}{The second general domination RL exponential} & \(HDRL_2(G,x)=\sum\limits_{uv\in E(G)}x^{\left(  d_d^2(u)+d_d^2(v)-d_d(u)d_d(v)\right)^a}\)\\[3ex]
		9&	\multirow{2}{6cm}{The first general domination RL exponential} & \(HDRL_1(G,x)=\sum\limits_{uv\in E(G)}x^{\left(  d_d^2(u)+d_d^2(v)+d_d(u)d_d(v)\right)^a}\)\\[3ex]
		
		10&	\multirow{2}{6cm}{The second hyper domination RL exponential} & \(HDRL_2(G,x)=\sum\limits_{uv\in E(G)}x^{\left(  d_d^2(u)+d_d^2(v)-d_d(u)d_d(v)\right)^2}\)\\[3ex]	
		\hline 
		11&	\multirow{2}{6cm}{The first  inverse domination RL index} & \(IDRL_1(G)=\sum\limits_{uv\in E(G)}\dfrac{1}{\left(  d_d^2(u)+d_d^2(v)+d_d(u)d_d(v)\right)} \)\\[3ex]
		
		12&	\multirow{2}{7cm}{The second  inverse domination RL index} & \(IDRL_2(G)=\sum\limits_{uv\in E(G)}\dfrac{1}{\left(  d_d^2(u)+d_d^2(v)-d_d(u)d_d(v)\right)} \)\\[3ex]
		13&	\multirow{2}{7cm}{The first  inverse domination RL exponential} & \(IDRL_1(G,x)=\sum\limits_{uv\in E(G)}x^{\dfrac{1}{\left(  d_d^2(u)+d_d^2(v)+d_d(u)d_d(v)\right)}} \)\\[3ex]
		
		14&	\multirow{2}{7cm}{The second  inverse domination RL exponential} & \(IDRL_2(G,x)=\sum\limits_{uv\in E(G)}x^{\dfrac{1}{\left(  d_d^2(u)+d_d^2(v)-d_d(u)d_d(v)\right)}} \)\\[3ex]

	\end{longtable}
	
	\begin{longtable}[c]{| c | l | l |}
		\caption{Multiplicative Rehan- Lanel domination indices. \label{tabled2}}\\
		
		\hline
		\multicolumn{1}{|c}{\textbf{S.N.}} & 
		\multicolumn{1}{|c}{\textbf{Name of the index}} & \multicolumn{1}{|c|}{\textbf{Formula}} \\\hline
		\endfirsthead
		
		\hline
		\endlastfoot

		1&	\multirow{2}{7cm}{The first multiplicative domination RL index} & \(DRL_1(G)=\prod\limits_{uv\in E(G)}\left(  d_d^2(u)+d_d^2(v)+d_d(u)d_d(v)\right)\)\\[2ex]
		2&	\multirow{2}{7cm}{The second multiplicative domination RL index} & \(DRL_1(G)=\prod\limits_{uv\in E(G)}\left(  d_d^2(u)+d_d^2(v)-d_d(u)d_d(v)\right)\)\\[3ex]
		
		3&	\multirow{2}{7cm}{The first multiplicative hyper domination RL index} & \(DRL_1(G)=\prod\limits_{uv\in E(G)}\left(  d_d^2(u)+d_d^2(v)+d_d(u)d_d(v)\right)^2\)\\[3ex]
		
		4&	\multirow{2}{7cm}{The second multiplicative hyper domination RL index} & \(DRL_2(G)=\prod\limits_{uv\in E(G)}\left(  d_d^2(u)+d_d^2(v)-d_d(u)d_d(v)\right)^2\)\\[3ex]
		
		5&	\multirow{2}{7cm}{The first multiplicative domination RL exponential} & \(DRL_1(G,x)=\prod\limits_{uv\in E(G)}x^{\left(  d_d^2(u)+d_d^2(v)+d_d(u)d_d(v)\right)}\)\\[3ex]
		
		6&	\multirow{2}{7cm}{The second multiplicative domination RL exponential} & \(DRL_2(G,x)=\prod\limits_{uv\in E(G)}x^{\left(  d_d^2(u)+d_d^2(v)-d_d(u)d_d(v)\right)}\)\\[4ex]
		7&	\multirow{2}{7cm}{The first multiplicative general domination RL index} & \(GDRL_1(G)=\prod\limits_{uv\in E(G)}\left(  d_d^2(u)+d_d^2(v)+d_d(u)d_d(v)\right)^a\)\\[6ex]
		
		8&	\multirow{2}{7cm}{\small The second multiplicative general domination RL index} & \( GDRL_2(G)=\prod\limits_{uv\in E(G)}\left(  d_d^2(u)+d_d^2(v)-d_d(u)d_d(v)\right)^a\)\\[3ex]
		
		9&	\multirow{2}{7cm}{The first multiplicative hyper  domination RL exponential} & \(HDRL_1(G,x)=\prod\limits_{uv\in E(G)}x^{\left(  d_d^2(u)+d_d^2(v)+d_d(u)d_d(v)\right)^2}\)\\[3ex]
		
		10&	\multirow{2}{7cm}{\small The second multiplicative hyper domination RL exponential} & \( HDRL_2(G,x)=\prod\limits_{uv\in E(G)}x^{\left(  d_d^2(u)+d_d^2(v)-d_d(u)d_d(v)\right)^2}\)\\[5ex]
		
		11&	\multirow{2}{7cm}{\small The first multiplicative inverse domination RL index} & \(IDRL_1(G)=\prod\limits_{uv\in E(G)}\dfrac{1}{\left(  d_d^2(u)+d_d^2(v)+d_d(u)d_d(v)\right)} \)\\[3ex]
		
		12&	\multirow{2}{7cm}{The second multiplicative  inverse domination RL index} & \(\small IDRL_2(G)=\prod\limits_{uv\in E(G)}\dfrac{1}{\left(  d_d^2(u)+d_d^2(v)-d_d(u)d_d(v)\right)} \)\\[3ex]
		13&	\multirow{2}{7cm}{The first multiplicative general  domination RL exponential} & \(GDRL_1(G,x)=\prod\limits_{uv\in E(G)}x^{\left(  d_d^2(u)+d_d^2(v)+d_d(u)d_d(v)\right)^a}\)\\[3ex]
		
		14&	\multirow{2}{7cm}{\small The second multiplicative general domination RL exponential} & \( GDRL_2(G,x)=\prod\limits_{uv\in E(G)}x^{\left(  d_d^2(u)+d_d^2(v)-d_d(u)d_d(v)\right)^2}\)\\[5ex]
		\hline 
		15&	\multirow{2}{7cm}{\small The first multiplicative inverse domination RL exponential} & \(IDRL_1(G,x)=\prod\limits_{uv\in E(G)}x^{\dfrac{1}{\left(  d_d^2(u)+d_d^2(v)+d_d(u)d_d(v)\right)}} \)\\[3ex]
		
		16&	\multirow{2}{7cm}{The second multiplicative  inverse domination RL exponential} & \(\small IDRL_2(G,x)=\prod\limits_{uv\in E(G)}x^{\dfrac{1}{\left(  d_d^2(u)+d_d^2(v)-d_d(u)d_d(v)\right)} }\)\\[3ex]
		
	\end{longtable}
	
		\begin{longtable}[c]{| c | l | l |}
		\caption{Rehan- Lanel domination indices/exponentials \label{tabled3}}\\

		\hline
		\multicolumn{1}{|c}{\textbf{S.N.}} & 
		\multicolumn{1}{|c}{\textbf{Name of the index}} & \multicolumn{1}{|c|}{\textbf{Formula}} \\\hline
		\endfirsthead
		
		\hline
		\endlastfoot

		1&	\multirow{2}{6cm}{The third  hyper domination RL index} & \(DRL_3(G)=\sum\limits_{uv\in E(G)}\left(  d_d(u)-d_d(v)+d_d(u)d_d(v)\right)^2\)\\[3ex]
		
		2&	\multirow{2}{6cm}{The fourth  hyper domination RL index} & \(DRL_4(G)=\sum\limits_{uv\in E(G)}\left(\left| d_d(u)-d_d(v)\right|d_d(u)d_d(v)\right)^2\)\\[3ex]
		
		3&	\multirow{2}{6cm}{The third domination RL exponential} & \(DRL_3(G,x)=\sum\limits_{uv\in E(G)}x^{\left(  d_d(u)-d_d(v)+d_d(u)d_d(v)\right)}\)\\[3ex]
		
		4&	\multirow{2}{6cm}{The fourth domination RL exponential} & \(DRL_2(G,x)=\sum\limits_{uv\in E(G)}x^{\left( \left| d_d(u)-d_d(v)\right|d_d(u)d_d(v)\right)}\)\\[3ex]

		5&	\multirow{2}{6cm}{The third  general domination RL index} & \(GDRL_1(G)=\sum\limits_{uv\in E(G)}\left(  d_d(u)-d_d(v)+d_d(u)d_d(v)\right)^a\)\\[3ex]
		
		6&	\multirow{2}{6cm}{The fourth  general domination RL index} & \(GDRL_2(G)=\sum\limits_{uv\in E(G)}\left( \left| d_d(u)-d_d(v)\right|d_d(u)d_d(v)\right)^a\)\\[3ex]
		7&	\multirow{2}{6cm}{The third hyper domination RL exponential} & \(HDRL_1(G,x)=\sum\limits_{uv\in E(G)}x^{\left(  d_d(u)-d_d(v)+d_d(u)d_d(v)\right)^2}\)\\[3ex]
		
		8&	\multirow{2}{6cm}{The fourth general domination RL exponential} & \(HDRL_2(G,x)=\sum\limits_{uv\in E(G)}x^{\left(  \left| d_d(u)-d_d(v)\right|d_d(u)d_d(v)\right)^a}\)\\[3ex]
		9&	\multirow{2}{6cm}{The third general domination RL exponential} & \(HDRL_1(G,x)=\sum\limits_{uv\in E(G)}x^{\left(  d_d(u)-d_d(v)+d_d(u)d_d(v)\right)^a}\)\\[3ex]
		
		10&	\multirow{2}{6cm}{The fourth hyper domination RL exponential} & \(HDRL_2(G,x)=\sum\limits_{uv\in E(G)}x^{\left(  \left| d_d(u)-d_d(v)\right|d_d(u)d_d(v)\right)^2}\)\\[3ex]	
		11&	\multirow{2}{6cm}{The third  inverse domination RL index} & \(IDRL_1(G)=\sum\limits_{uv\in E(G)}\dfrac{1}{\left(  d_d(u)-d_d(v)+d_d(u)d_d(v)\right)} \)\\[3ex]
		
		12&	\multirow{2}{7cm}{The fourth  inverse domination RL index} & \(IDRL_2(G)=\sum\limits_{uv\in E(G)}\dfrac{1}{\left(  \left| d_d(u)-d_d(v)\right|d_d(u)d_d(v)\right)} \)\\[3ex]
		13&	\multirow{2}{7cm}{The third  inverse domination RL exponential} & \(IDRL_1(G,x)=\sum\limits_{uv\in E(G)}x^{\dfrac{1}{\left(  d_d(u)-d_d(v)+d_d(u)d_d(v)\right)}} \)\\[3ex]
		
		14&	\multirow{2}{7cm}{The fourth  inverse domination RL exponential} & \(IDRL_2(G,x)=\sum\limits_{uv\in E(G)}x^{\dfrac{1}{\left(  \left| d_d(u)-d_d(v)\right|d_d(u)d_d(v)\right)}} \)\\[3ex]

	\end{longtable}
	
	\begin{longtable}[c]{| c | l | l |}
		\caption{Multiplicative Rehan- Lanel domination indices. \label{tabled4}}\\
		
		\hline
		\multicolumn{1}{|c}{\textbf{S.N.}} & 
		\multicolumn{1}{|c}{\textbf{Name of the index}} & \multicolumn{1}{|c|}{\textbf{Formula}} \\\hline
		\endfirsthead
		
		\hline
		\endlastfoot

		1&	\multirow{2}{7cm}{The third multiplicative domination RL index} & \(DRL_1(G)=\prod\limits_{uv\in E(G)}\left(  d_d(u)-d_d(v)+d_d(u)d_d(v)\right)\)\\[2ex]
		\hline 
		2&	\multirow{2}{7cm}{The fourth multiplicative domination RL index} & \(DRL_1(G)=\prod\limits_{uv\in E(G)}\left(  \left| d_d(u)-d_d(v)\right|d_d(u)d_d(v)\right)\)\\[3ex]
		
		3&	\multirow{2}{7cm}{The third multiplicative hyper domination RL index} & \(DRL_1(G)=\prod\limits_{uv\in E(G)}\left(  d_d(u)-d_d(v)+d_d(u)d_d(v)\right)^2\)\\[3ex]
		
		4&	\multirow{2}{7cm}{The fourth multiplicative hyper domination RL index} & \(DRL_2(G)=\prod\limits_{uv\in E(G)}\left(  \left| d_d(u)-d_d(v)\right|d_d(u)d_d(v)\right)^2\)\\[3ex]
		
		5&	\multirow{3}{7cm}{The third multiplicative domination RL exponential} & \(DRL_1(G,x)=\prod\limits_{uv\in E(G)}x^{\left(  d_d(u)-d_d(v)+d_d(u)d_d(v)\right)}\)\\[3ex]
		
		6&	\multirow{3}{7cm}{The fourth multiplicative domination RL exponential} & \(DRL_2(G,x)=\prod\limits_{uv\in E(G)}x^{\left(  \left| d_d(u)-d_d(v)\right|d_d(u)d_d(v)\right)}\)\\[3ex]
		
		7&	\multirow{2}{7cm}{The third multiplicative general domination RL index} & \(GDRL_1(G)=\prod\limits_{uv\in E(G)}\left(  d_d(u)-d_d(v)+d_d(u)d_d(v)\right)^a\)\\[6ex]
		
		8&	\multirow{2}{7cm}{\small The fourth multiplicative general domination RL index} & \( GDRL_2(G)=\prod\limits_{uv\in E(G)}\left(  \left| d_d(u)-d_d(v)\right|d_d(u)d_d(v)\right)^a\)\\[3ex]
		9&	\multirow{2}{7cm}{The third multiplicative hyper  domination RL exponential} & \(HDRL_1(G,x)=\prod\limits_{uv\in E(G)}x^{\left(  d_d(u)-d_d(v)+d_d(u)d_d(v)\right)^2}\)\\[3ex]
		
		10&	\multirow{2}{7cm}{\small The fourth multiplicative hyper domination RL exponential} & \( HDRL_2(G,x)=\prod\limits_{uv\in E(G)}x^{\left(  \left| d_d(u)-d_d(v)\right|d_d(u)d_d(v)\right)^2}\)\\[5ex]
		
		11&	\multirow{2}{7cm}{\small The third multiplicative inverse domination RL index} & \(IDRL_1(G)=\prod\limits_{uv\in E(G)}\dfrac{1}{\left(  d_d(u)-d_d(v)+d_d(u)d_d(v)\right)} \)\\[3ex]
		
		12&	\multirow{2}{7cm}{The fourth multiplicative  inverse domination RL index} & \(\small IDRL_2(G)=\prod\limits_{uv\in E(G)}\dfrac{1}{\left(  \left| d_d(u)-d_d(v)\right|d_d(u)d_d(v)\right)} \)\\[3ex]
		13&	\multirow{2}{7cm}{The third multiplicative general  domination RL exponential} & \(GDRL_1(G,x)=\prod\limits_{uv\in E(G)}x^{\left(  d_d(u)-d_d(v)+d_d(u)d_d(v)\right)^a}\)\\[3ex]
		
		14&	\multirow{2}{7cm}{\small The fourth multiplicative general domination RL exponential} & \( GDRL_2(G,x)=\prod\limits_{uv\in E(G)}x^{\left(  \left| d_d(u)-d_d(v)\right|d_d(u)d_d(v)\right)^2}\)\\[5ex]
		15&	\multirow{2}{7cm}{\small The third multiplicative inverse domination RL exponential} & \(IDRL_1(G,x)=\prod\limits_{uv\in E(G)}x^{\dfrac{1}{\left(  d_d(u)-d_d(v)+d_d(u)d_d(v)\right)}} \)\\[3ex]
		
		16&	\multirow{2}{7cm}{The fourth multiplicative  inverse domination RL exponential} & \(\small IDRL_2(G,x)=\prod\limits_{uv\in E(G)}x^{\dfrac{1}{\left(  \left| d_d(u)-d_d(v)\right|d_d(u)d_d(v)\right)} }\)\\[3ex]
		
	\end{longtable}

	In this study, we calculate the first and second RL domination indices, together with their corresponding polynomials, for specific standard graphs such as a complete graph, star graph, double star graph, a complete bipartite graph and French windmill graph.
	
	\begin{proposition}
		Let \(K_n\) is a complete graph with \(n\) vertices. Then,
		
		\begin{enumerate}[label=\roman*.]
			\item \(DRL_1(K_n)=\dfrac{3n(n-1)}{2}.\)
			\item \(DRL_2(K_n)=\dfrac{n(n-1)}{2}.\)
			\item \(DRL_3(K_n)=\dfrac{n(n-1)}{2}\)
			\item \(DRL_4(K_n)=0\).
		\end{enumerate}
	
	\end{proposition}
	\begin{proof}
		Let \(G\) is an $r-$regular graph with \(n\) vertices and
		\(r \geq 2\).

		Then \(G\) has \(\dfrac{n(n-1)}{2}\)	edges and \(d_d(u)=1\).\\
		From definition, we have\\ 
		
		$$\begin{array}{ll}
			DRL_1(K_n) &=\sum\limits_{uv\in E(G)}\left(  d_d^2(u)+d_d^2(v)+d_d(u)d_d(v)\right)  \\
			& =\dfrac{n(n-1)}{2}(1^2+1^2+1\cdot 1)\\
			&=\dfrac{3n(n-1)}{2}.
		\end{array}$$
		
		$$\begin{array}{ll}
			DRL_2(K_n) &=\sum\limits_{uv\in E(G)}\left(  d_d^2(u)+d_d^2(v)-d_d(u)d_d(v)\right)  \\
			& =\dfrac{n(n-1)}{2}(1^2+1^2-1\cdot 1)\\
			&=\dfrac{n(n-1)}{2}.
		\end{array}$$
		
		$$\begin{array}{ll}
			DRL_3(K_n) &=\sum\limits_{uv\in E(G)}\left(  d_d(u)-d_d(v)+d_d(u)d_d(v)\right)  \\
			& =\dfrac{n(n-1)}{2}(1-1+1\cdot 1)\\
			&=\dfrac{n(n-1)}{2}.
		\end{array}$$
		
		$$\begin{array}{ll}
			DRL_4(K_n) &=\sum\limits_{uv\in E(G)}\left( \left| d_d(u)-d_d(v)\right|d_d(u)d_d(v)\right)  \\
			& =\dfrac{n(n-1)}{2}(|1-1|1\cdot 1)\\
			&=0.
		\end{array}$$
	\end{proof}
	
	\begin{proposition}
		Let \(S_{n+1}\) is a star graph with \(d_d(u)=1\). Then,
		\begin{enumerate}[label=\roman*.]
			\item \(DRL_1(S_{n+1})=3n\).
			\item \(DRL_2(S_{n+1})=n\).
				\item \(DRL_3(S_{n+1})=n\).
			\item \(DRL_4(S_{n+1})=0\).
		\end{enumerate}
	\end{proposition}
	\begin{proof}
		The proof can be done similarly as above.
	\end{proof}
	
	\begin{proposition}
		Let \(S_{p+1, q+1}\) is a double star graph with \(d_d(u)=2\). Then,
		\begin{enumerate}[label=\roman*.]
			\item \(DRL_1((S_{p+1, q+1})=12(p+q+1)\).
			\item \(DRL_2((S_{p+1, q+1})=4(p+q+1)\).
			\item \(DRL_3((S_{p+1, q+1})=4(p+q+1)\).
			\item \(DRL_4((S_{p+1, q+1})=0\).
		\end{enumerate}
	\end{proposition}
	
	\begin{proof}
		The proof can be done similarly as above.
	\end{proof}
	\begin{proposition}
		Let \(K_{m,n}\) be a complete bipartite graph with \(2 \leq m \leq n\). Then,
		\begin{enumerate}[label=\roman*.]
			\item \(DRL_1(K_{m,n})=mn(m^2+n^2+mn+3m+3n+3)\).
			\item \(DRL_2(K_{m,n})=mn(m^2+n^2-mn+m+n+1)\).
			\item \(DRL_3(K_{m,n})=mn(mn+2n+1)\).
				\item \(DRL_4(K_{m,n})=mn|n-m|(m+1)(n+1)\).
		\end{enumerate}
	\end{proposition}
	\begin{proof}
		The proof can be done similarly as above.
	\end{proof}

	In the following proposition, by using definition, we present the first and second Rehan-Lanel domination polynomials of \(K_n, S_{n+1},S_{p+1, q+1}, K_{m,n}. \)

	\begin{proposition}
		The first, second, third and fourth Rehan-Lanel domination polynomials of \(K_n, S_{n+1},S_{p+1, q+1}, K_{m,n}. \) are given by,
		
		\begin{enumerate}[label=\roman*.]
			\item \(DRL_1(K_n,x)=\dfrac{n(n-1)}{2}x^3\).
			\item \(DRL_2(K_n,x)=\dfrac{n(n-1)}{2}x\).
			\item \(DRL_3(K_n,x)=\dfrac{n(n-1)}{2}x\).
			\item \(DRL_4(K_n,x)=\dfrac{n(n-1)}{2}\).
			\item \(DRL_1(S_{n+1},x)=nx^3\).
			\item \(DRL_2(S_{n+1},x)=nx\).
			\item \(DRL_3(S_{n+1},x)=nx\).
			\item \(DRL_4(S_{n+1},x)=n\)
			\item \(DRL_1((S_{p+1, q+1},x)=(p+q+1)x^{12}\).
			\item \(DRL_2((S_{p+1, q+1},x)=(p+q+1)x^4\).
			\item \(DRL_1(K_{m,n},x)=mnx^{(m^2+n^2+mn+3m+3n+3)}.\)
			\item \(DRL_2(K_{m,n},x)=mnx^{(m^2+n^2-mn+m+n+1)}.\)
		\end{enumerate}
	\end{proposition}
	
	\begin{proof}
		The proof can be done similarly as above.
	\end{proof}

	\subsubsection{French Windmill graph}
	The French windmill graph \(F_n^m\) is determined by taking \(m\geq 3\) copies of \(K_n\), \(n\geq 3\) with a vertex in common as in the figure. The French windmill graph \(F_3^m\) is a friendship graph.

	\begin{figure}[h]
		\centering 
		\includegraphics[width=8cm, height=8cm]{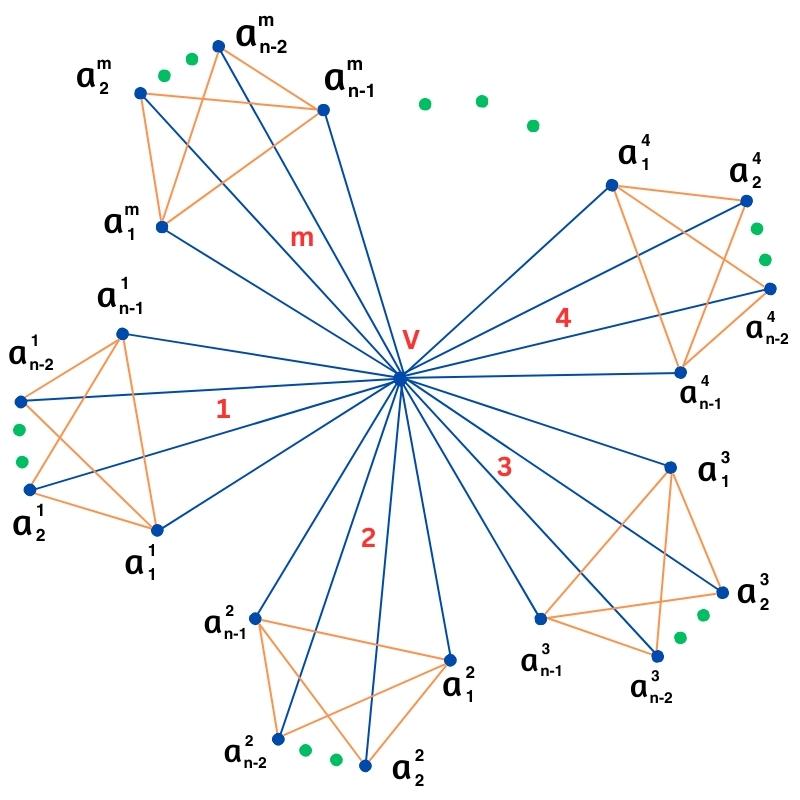}
		\caption{French Windmill graph}
		\label{fig:fw}
	\end{figure}

	Let \(F\) be a French windmill graph \(F_n^m\) . Then
	$$d_d(u)=\begin{cases}
		1,  \ \text{if} \ u \ \text{is a center},\\
		(n-1)^{m-1}, \ \text{otherwise.}
	\end{cases}$$
	
	\begin{theorem}
		Let \(F\) be a French windmill graph \(F_n^m\). Then,\begin{enumerate}[label=\roman*.]
			\item\(DRL_1(F,x)=m(n-1)\left[ (n-1)^{2(m-1)}+(n-1)^{m-1}+1\right] +3\left[ \dfrac{mn(n-1)(n-2)}{2}\right] \left[ (n-1)^{2(m-1)}\right] \).
			\item \(DRL_2(F,x)=m(n-1)\left[ (n-1)^{2(m-1)}-(n-1)^{m-1}+1\right] +\left[ \dfrac{mn(n-1)(n-2)}{2}\right] \left[ (n-1)^{2(m-1)}\right]\).
		\end{enumerate}
	\end{theorem}
	
	\begin{proof}
		In \(F\), there are two sets of edges. Let \(E_1\) be the set of all edges which are incident with the
		center vertex and \(E_2\) be the set of all edges of the complete graph. Then,
		
		$$	\begin{array}{rl}
			DRL_1(F)=&\sum\limits_{uv\in E(F)}\left(  d_d^2(u)+d_d^2(v)+d_d(u)d_d(v)\right)	 \\
			=& \sum\limits_{uv\in E_1(F)}\left(  d_d^2(u)+d_d^2(v)+d_d(u)d_d(v)\right)+\sum\limits_{uv\in E_2(F)}\left(  d_d^2(u)+d_d^2(v)+d_d(u)d_d(v)\right)\\
			=& m(n-1)\left[1^2+ ((n-1)^{(m-1)})^2+1.(n-1)^{m-1}\right]+\\
			&  \left[ \dfrac{mn(n-1)}{2}-m(n-1)\right] \left[ (n-1)^{(m-1)}+(n-1)^{(m-1)}+(n-1)^{(m-1)}(n-1)^{(m-1)}\right]\\
			=& m(n-1)\left[ (n-1)^{2(m-1)}+(n-1)^{m-1}+1\right] +3\left[ \dfrac{mn(n-1)(n-2)}{2}\right] \left[ (n-1)^{2(m-1)}\right]. 
		\end{array}$$
		
		$$	\begin{array}{rl}
			DRL_2(F)=&\sum\limits_{uv\in E(F)}\left(  d_d^2(u)+d_d^2(v)-d_d(u)d_d(v)\right)	 \\
			=& \sum\limits_{uv\in E_1(F)}\left(  d_d^2(u)+d_d^2(v)-d_d(u)d_d(v)\right)+\sum\limits_{uv\in E_2(F)}\left(  d_d^2(u)+d_d^2(v)-d_d(u)d_d(v)\right)\\
			=& m(n-1)\left[1^2+ ((n-1)^{(m-1)})^2-1.(n-1)^{m-1}\right]+\\
			&  \left[ \dfrac{mn(n-1)}{2}-m(n-1)\right] \left[ (n-1)^{2(m-1)}+(n-1)^{2(m-1)}-(n-1)^{(m-1)}(n-1)^{(m-1)}\right]\\
			=& m(n-1)\left[ (n-1)^{2(m-1)}-(n-1)^{m-1}+1\right] +\left[ \dfrac{mn(n-1)(n-2)}{2}\right] \left[ (n-1)^{2(m-1)}\right].
		\end{array}$$
		
	\end{proof}
	We can compute the results for all the other indices/exponentials listed in the Table  \ref{table1} and \ref{table2} similarly.
	\subsection{Temperature RL indices}
	The temperature of a vertex u of a graph G is defined by Fajtlowicz in \cite{article46} as \[T(u)=\dfrac{d_G(u)}{n-d_G(u)},\]
	
	where $n$ is the number of vertices of $G$.\\
	
	In \cite{article47}, Kulli introduced few temperature indices such as the first temperature index and the modified first temperature index  of a graph $G$ which are defined as
	\[T_1(G)=\sum_{u\in V(G)}T(u)^2\]
	\[^mT_1(G)=\sum_{u\in V(G)}\dfrac{1}{T(u)^2}.\]
	Recently, some temperature indices were introduced and studied, for example, in \cite{article48}\\
	
	Inspired by these approaches, here  we introduce the first and second Rehan-Lanel temperature indices as follows.
	
	\begin{definition}
		The first, second, third and fourth Rehan-Lanel temperature  indices of a graph \(G\) are defined by,
		\begin{enumerate}[label=\roman*.]
			\item \(TRL_1(G)=\sum\limits_{uv\in E(G)}\left(  T^2(u)+T^2(v)+T(u)T(v)\right) \)
			\item \(TRL_2(G)=\sum\limits_{uv\in E(G)}\left(  T^2(u)+T^2(v)-T(u)T(v)\right) .\)
			\item \(TRL_3(G)=\sum\limits_{uv\in E(G)}\left(  T(u)-T(v)+T(u)T(v)\right) \)
			\item \(TRL_4(G)=\sum\limits_{uv\in E(G)}\left|  T(u)-T(v)\right|T(u)T(v) .\)
		\end{enumerate}
	\end{definition}
	
	We introduce an additional indices/exponentials incorporating with these two indices as in the following Table \ref{tablee1}, \ref{tablee2}, \ref{tablee3} and \ref{tablee4}. 
	\begin{longtable}[c]{| c | l | l |}
		\caption{Temperature Rehan- Lanel  indices/exponentials\label{tablee1}}\\

		\hline
		\multicolumn{1}{|c}{\textbf{S.N.}} & 
		\multicolumn{1}{|c}{\textbf{Name of the index}} & \multicolumn{1}{|c|}{\textbf{Formula}} \\\hline
		\endfirsthead
		
		\hline
		\endlastfoot

		1&	The first hyper Temperature RL index & \(HTRL_1(G)=\sum\limits_{uv\in E(G)}\left(  T^2(u)+T^2(v)+T(u)T(v)\right)^2 \)\\[3ex]
		
		2&	The second  hyper Temperature RL index & \(HTRL_2(G)=\sum\limits_{uv\in E(G)}\left(  T^2(u)+T^2(v)-T(u)T(v)\right)^2 \)\\[3ex]
		
		3&	The first general Temperature RL index & \(GTRL_1(G)=\sum\limits_{uv\in E(G)}\left(  T^2(u)+T^2(v)+T(u)T(v)\right)^a \)\\[3ex]
		4&	The second  general Temperature RL index & \(GTRL_2(G)=\sum\limits_{uv\in E(G)}\left(  T^2(u)-T^2(v)+T(u)T(v)\right)^a \)\\[3ex]
		5&	\multirow{2}{7cm}{The first inverse Temperature RL index} & \(ITRL_1(G)=\sum\limits_{uv\in E(G)}\dfrac{1}{\left(  T^2(u)+T^2(v)+T(u)T(v)\right)} \)\\[3ex]
		
		6&	\multirow{2}{7cm}{The second inverse Temperature RL index} & \(ITRL_2(G)=\sum\limits_{uv\in E(G)}\dfrac{1}{\left(  T^2(u)+T^2(v)-T(u)T(v)\right)} \)\\[3ex]
		
		7&	\multirow{2}{7cm}{The first Temperature RL exponential} & \(TRL_1(G,x)=\sum\limits_{uv\in E(G)}x^{\left(  T^2(u)+T^2(v)+T(u)T(v)\right)} \)\\[3ex]
		8&	\multirow{2}{7cm}{The second  Temperature RL exponential} & \(TRL_2(G,x)=\sum\limits_{uv\in E(G)}x^{\left(  T^2(u)+T^2(v)-T(u)T(v)\right)} \)\\[3ex]
		9&	\multirow{2}{7cm}{The first hyper Temperature RL exponential} & \(HTRL_1(G,x)=\sum\limits_{uv\in E(G)}x^{\left(  T^2(u)+T^2(v)+T(u)T(v)\right)^2 }\)\\[3ex]
		
		10&	\multirow{2}{7cm}{The second  hyper Temperature RL exponential} & \(HTRL_2(G,x)=\sum\limits_{uv\in E(G)}x^{\left(  T^2(u)+T^2(v)-T(u)T(v)\right)^2} \)\\[5ex]
		\hline 
		11&	\multirow{2}{7cm}{The first inverse Temperature RL exponential} & \(\displaystyle ITRL_1(G,x)=\sum\limits_{uv\in E(G)}x^{\dfrac{1}{\left(  T^2(u)+T^2(v)+T(u)T(v)\right)}} \)\\[3ex]
		12&	\multirow{2}{7cm}{The second inverse Temperature RL exponential} & \(\displaystyle ITRL_2(G,x)=\sum\limits_{uv\in E(G)}x^{\dfrac{1}{\left(  T^2(u)+T^2(v)-T(u)T(v)\right)} }\)\\[3ex]
		
		13&	\multirow{2}{7cm}{The first general Temperature RL exponential} & \(\displaystyle GTRL_1(G,x)=\sum\limits_{uv\in E(G)}x^{\left(  T^2(u)+T^2(v)+T(u)T(v)\right)^a} \)\\[3ex]
		
		14&	\multirow{2}{7cm}{The second  general Temperature RL exponential} & \(\displaystyle GTRL_2(G,x)=\sum\limits_{uv\in E(G)}x^{\left(  T^2(u)+T^2(v)-T(u)T(v)\right)^a} \)\\[3ex]

	\end{longtable}
	
	\begin{longtable}[c]{| c | l | l |}
		\caption{Multiplicative Temperature Rehan- Lanel  indices/exponentials \label{tablee2}}\\

		\hline
		\multicolumn{1}{|c}{\textbf{S.N.}} & 
		\multicolumn{1}{|c}{\textbf{Name of the index}} & \multicolumn{1}{|c|}{\textbf{Formula}} \\\hline
		\endfirsthead
		\hline
		\endlastfoot

		1&	\multirow{2}{7cm}{The first multiplicative Temperature RL index }& \(MTRL_1(G)=\prod\limits_{uv\in E(G)}\left( T^2(u)+T^2(v)+T(u)T(v)\right) \)\\[3ex]
		
		2&	\multirow{2}{7cm}{The second multiplicative Temperature RL index }& \(MTRL_2(G)=\prod\limits_{uv\in E(G)}\left( T^2(u)+T^2(v)-T(u)T(v)\right) \)\\[3ex]
		
		3&		\multirow{2}{7cm}{The first  multiplicative hyper Temperature RL index} & \(MHTRL_1(G)=\prod\limits_{uv\in E(G)}\left(  T^2(u)+T^2(v)+T(u)T(v)\right)^2 \)\\[5ex]
		
		4&		\multirow{2}{7cm}{The second multiplicative hyper Temperature RL index} & \(MHTRL_2(G)=\prod\limits_{uv\in E(G)}\left(  T^2(u)+T^2(v)-T(u)T(v)\right)^2 \)\\[5ex]
		
		5&	\multirow{2}{7cm}{The first multiplicative inverse Temperature RL index} & \(MITRL_1(G)=\prod\limits_{uv\in E(G)}\dfrac{1}{\left(  T^2(u)+T^2(v)+T(u)T(v)\right)} \)\\[6ex]
		
		6&	\multirow{2}{7cm}{The second multiplicative inverse Temperature RL index} & \(MITRL_2(G)=\prod\limits_{uv\in E(G)}\dfrac{1}{\left(  T^2(u)+T^2(v)-T(u)T(v)\right)} \)\\[6ex]
		7&	\multirow{2}{7cm}{The first general multiplicative Temperature RL index} & \(MGTRL_1(G)=\prod\limits_{uv\in E(G)}\left( T^2(u)+T^2(v)+T(u)T(v)\right)^a \)\\[6ex]
		
		8&	\multirow{2}{7cm}{The second  general multiplicative Temperature RL index} & \(MGTRL_2(G)=\prod\limits_{uv\in E(G)}\left( T^2(u)+T^2(v)-T(u)T(v)\right) \)\\[6ex]
		
		9&	\multirow{2}{7cm}{The first multiplicative Temperature RL exponential} & \(MTRL_1(G,x)=\prod\limits_{uv\in E(G)}x^{\left(T^2(u)+T^2(v)+T(u)T(v)\right)} \)\\[3ex]
		
		10&	\multirow{2}{7cm}{The second multiplicative Temperature RL exponential} & \(MTRL_2(G,x)=\prod\limits_{uv\in E(G)}x^{\left(  T^2(u)+T^2(v)-T(u)T(v)\right)} \)\\[3ex]
		
		11&	\multirow{2}{7cm}{The first multiplicative general Temperature RL exponential} & \(MGTRL_1(G,x)=\prod\limits_{uv\in E(G)}x^{\left(  T^2(u)+T^2(v)+T(u)T(v)\right)^a }\)\\[5ex]
		12&	\multirow{2}{7cm}{The second multiplicative general Temperature RL exponential} & \(MGTRL_2(G,x)=\prod\limits_{uv\in E(G)}x^{\left( T^2(u)+T^2(v)-T(u)T(v)\right)} \)\\[3ex]
		\hline
		13&	\multirow{2}{7cm}{The first multiplicative inverse Temperature RL exponential} & \(MITRL_1(G,x )=\prod\limits_{uv\in E(G)}x^{\dfrac{1}{\left(  T^2(u)+T^2(v)+T(u)T(v)\right)}} \)\\[3ex]
		14&	\multirow{2}{7cm}{The second multiplicative inverse Temperature RL exponential} & \(MITRL_2(G,x)=\prod\limits_{uv\in E(G)}x^{\dfrac{1}{\left(  T^2(u)+T^2(v)-T(u)T(v)\right)} }\)\\[5ex]
		
		15&	\multirow{2}{7cm}{The first multiplicative general Temperature RL exponential} & \(MGTRL_1(G,x)=\prod\limits_{uv\in E(G)}x^{\left(  T^2(u)+T^2(v)+T(u)T(v)\right)^a} \)\\[5ex]
		
		16&	\multirow{2}{7cm}{The second multiplicative general Temperature RL exponential} & \(MGTRL_2(G,x)=\prod\limits_{uv\in E(G)}x^{\left(  T^2(u)+T^2(v)-T(u)T(v)\right)} \)\\[6ex]

	\end{longtable}
	
	\begin{longtable}[c]{| c | l | l |}
		\caption{Temperature Rehan- Lanel  indices/exponentials\label{tablee3}}\\

		\hline
		\multicolumn{1}{|c}{\textbf{S.N.}} & 
		\multicolumn{1}{|c}{\textbf{Name of the index}} & \multicolumn{1}{|c|}{\textbf{Formula}} \\\hline
		\endfirsthead
		
		\hline
		\endlastfoot

		1&	The third hyper Temperature RL index & \(HTRL_1(G)=\sum\limits_{uv\in E(G)}\left(   T(u)-T(v)+T(u)T(v)\right)^2 \)\\[3ex]
		
		2&	The fourth  hyper Temperature RL index & \(HTRL_2(G)=\sum\limits_{uv\in E(G)}\left( \left| T(u)-T(v)\right|T(u)T(v)\right)^2 \)\\[3ex]
		
		3&	The third general Temperature RL index & \(GTRL_1(G)=\sum\limits_{uv\in E(G)}\left(   T(u)-T(v)+T(u)T(v)\right)^a \)\\[3ex]
		4&	The fourth  general Temperature RL index & \(GTRL_2(G)=\sum\limits_{uv\in E(G)}\left( \left| T(u)-T(v)\right|T(u)T(v)\right)^a \)\\[3ex]
		
		5&	\multirow{2}{7cm}{The third inverse Temperature RL index} & \(ITRL_1(G)=\sum\limits_{uv\in E(G)}\dfrac{1}{\left(   T(u)-T(v)+T(u)T(v)\right)} \)\\[3ex]
		
		6&	\multirow{2}{7cm}{The fourth inverse Temperature RL index} & \(ITRL_2(G)=\sum\limits_{uv\in E(G)}\dfrac{1}{\left(  \left| T(u)-T(v)\right|T(u)T(v)\right)} \)\\[3ex]
		
		7&	\multirow{2}{7cm}{The third Temperature RL exponential} & \(TRL_1(G,x)=\sum\limits_{uv\in E(G)}x^{\left(   T(u)-T(v)+T(u)T(v)\right)} \)\\[3ex]
		8&	\multirow{2}{7cm}{The fourth  Temperature RL exponential} & \(TRL_2(G,x)=\sum\limits_{uv\in E(G)}x^{\left(  \left| T(u)-T(v)\right|T(u)T(v)\right)} \)\\[3ex]
		
		9&	\multirow{2}{7cm}{The third hyper Temperature RL exponential} & \(HTRL_1(G,x)=\sum\limits_{uv\in E(G)}x^{\left(   T(u)-T(v)+T(u)T(v)\right)^2 }\)\\[3ex]
		
		10&	\multirow{2}{7cm}{The fourth  hyper Temperature RL exponential} & \(HTRL_2(G,x)=\sum\limits_{uv\in E(G)}x^{\left( \left| T(u)-T(v)\right|T(u)T(v)\right)^2} \)\\[3ex]
		11&	\multirow{2}{7cm}{The third inverse Temperature RL exponential} & \(\displaystyle ITRL_1(G,x)=\sum\limits_{uv\in E(G)}x^{\dfrac{1}{\left(  T(u)-T(v)+T(u)T(v)\right)}} \)\\[3ex]
		
		12&	\multirow{2}{7cm}{The fourth inverse Temperature RL exponential} & \(\displaystyle ITRL_2(G,x)=\sum\limits_{uv\in E(G)}x^{\dfrac{1}{\left(  \left| T(u)-T(v)\right|T(u)T(v)\right)} }\)\\[3ex]
		
		13&	\multirow{2}{7cm}{The third general Temperature RL exponential} & \(\displaystyle GTRL_1(G,x)=\sum\limits_{uv\in E(G)}x^{\left(   T(u)-T(v)+T(u)T(v)\right)^a} \)\\[3ex]
		
		14&	\multirow{2}{7cm}{The fourth  general Temperature RL exponential} & \(\displaystyle GTRL_2(G,x)=\sum\limits_{uv\in E(G)}x^{\left(  \left| T(u)-T(v)\right|T(u)T(v)\right)^a} \)\\[3ex]

	\end{longtable}
	
	\begin{longtable}[c]{| c | l | l |}
		\caption{Multiplicative Temperature Rehan- Lanel  indices/exponentials \label{tablee4}}\\

		\hline
		\multicolumn{1}{|c}{\textbf{S.N.}} & 
		\multicolumn{1}{|c}{\textbf{Name of the index}} & \multicolumn{1}{|c|}{\textbf{Formula}} \\\hline
		\endfirsthead
		\hline
		\endlastfoot

		1&	\multirow{2}{7cm}{The third multiplicative Temperature RL index }& \(MTRL_1(G)=\prod\limits_{uv\in E(G)}\left(  T(u)-T(v)+T(u)T(v)\right) \)\\[3ex]
		
		2&	\multirow{2}{7cm}{The fourth multiplicative Temperature RL index }& \(MTRL_2(G)=\prod\limits_{uv\in E(G)}\left( \left| T(u)-T(v)\right|T(u)T(v)\right) \)\\[3ex]
		
		3&		\multirow{2}{7cm}{The third  multiplicative hyper Temperature RL index} & \(MHTRL_1(G)=\prod\limits_{uv\in E(G)}\left(   T(u)-T(v)+T(u)T(v)\right)^2 \)\\[5ex]
		
		4&		\multirow{2}{7cm}{The fourth multiplicative hyper Temperature RL index} & \(MHTRL_2(G)=\prod\limits_{uv\in E(G)}\left(  \left| T(u)-T(v)\right|T(u)T(v)\right)^2 \)\\[5ex]
		\hline
		5&	\multirow{2}{7cm}{The third multiplicative inverse Temperature RL index} & \(MITRL_1(G)=\prod\limits_{uv\in E(G)}\dfrac{1}{\left(   T(u)-T(v)+T(u)T(v)\right)} \)\\[6ex]
		
		6&	\multirow{2}{7cm}{The fourth multiplicative inverse Temperature RL index} & \(MITRL_2(G)=\prod\limits_{uv\in E(G)}\dfrac{1}{\left(  \left| T(u)-T(v)\right|T(u)T(v)\right)} \)\\[3ex]
		7&	\multirow{2}{7cm}{The third general multiplicative Temperature RL index} & \(MGTRL_1(G)=\prod\limits_{uv\in E(G)}\left(  T(u)-T(v)+T(u)T(v)\right)^a \)\\[6ex]
		
		8&	\multirow{2}{7cm}{The fourth  general multiplicative Temperature RL index} & \(MGTRL_2(G)=\prod\limits_{uv\in E(G)}\left( \left| T(u)-T(v)\right|T(u)T(v)\right) \)\\[6ex]
		
		9&	\multirow{2}{7cm}{The third multiplicative Temperature RL exponential} & \(MTRL_1(G,x)=\prod\limits_{uv\in E(G)}x^{\left( T(u)-T(v)+T(u)T(v)\right)} \)\\[3ex]
		
		10&	\multirow{2}{7cm}{The fourth multiplicative Temperature RL exponential} & \(MTRL_2(G,x)=\prod\limits_{uv\in E(G)}x^{\left(  \left| T(u)-T(v)\right|T(u)T(v)\right)} \)\\[3ex]
		
		11&	\multirow{2}{7cm}{The third multiplicative general Temperature RL exponential} & \(MGTRL_1(G,x)=\prod\limits_{uv\in E(G)}x^{\left(   T(u)-T(v)+T(u)T(v)\right)^a }\)\\[5ex]
		
		12&	\multirow{2}{7cm}{The fourth multiplicative general Temperature RL exponential} & \(MGTRL_2(G,x)=\prod\limits_{uv\in E(G)}x^{\left( \left| T(u)-T(v)\right|T(u)T(v)\right)} \)\\[3ex]
		13&	\multirow{2}{7cm}{The third multiplicative inverse Temperature RL exponential} & \(MITRL_1(G,x )=\prod\limits_{uv\in E(G)}x^{\dfrac{1}{\left(   T(u)-T(v)+T(u)T(v)\right)}} \)\\[3ex]
		14&	\multirow{2}{7cm}{The fourth multiplicative inverse Temperature RL exponential} & \(MITRL_2(G,x)=\prod\limits_{uv\in E(G)}x^{\dfrac{1}{\left(  \left| T(u)-T(v)\right|T(u)T(v)\right)} }\)\\[5ex]
		
		15&	\multirow{2}{7cm}{The third multiplicative general Temperature RL exponential} & \(MGTRL_1(G,x)=\prod\limits_{uv\in E(G)}x^{\left(   T(u)-T(v)+T(u)T(v)\right)^a} \)\\[5ex]
		
		16&	\multirow{2}{7cm}{The fourth multiplicative general Temperature RL exponential} & \(MGTRL_2(G,x)=\prod\limits_{uv\in E(G)}x^{\left(  \left| T(u)-T(v)\right|T(u)T(v)\right)} \)\\[6ex]

	\end{longtable}

	We compute the first, second, third and fourth temperature Rehan-Lanel indices\((TRL_1, TRL_2, TRL_4, TRL_4)\) for the \(r-\) regular graph, cycle, complete graph, path and complete bipartite graph. In addition, we present the results for wheel graph and sunflower graph.
	\begin{proposition}
		Let \(G\) is an $r$-regular graph with \(n\) vertices and
		\(r \geq 2\). Then, 
		\begin{enumerate}[label=\roman*.]
			\item \(TRL_1(G)=\dfrac{3nr^3}{2(n-r)^2}\)
			\item \(TRL_2(G)=\dfrac{nr^3}{2(n-r)^2}\)
			\item \(TRL_3(G)=\dfrac{3nr^3}{2(n-r)^2}\)
			\item \(TRL_4(G)=0.\)
		\end{enumerate}

	\end{proposition}  
	
	\begin{proof}
		Let \(G\) is an $r-$regular graph with \(n\) vertices and
		\(r \geq 2\).

		Then, \(G\) has \(\dfrac{nr}{2}\)	edges. Also, \(T(u)=\dfrac{r}{n-r}\).
		From definition we have
		$$ \begin{array}{ll}
			TRL_1(G)&=\sum\limits_{uv\in E(G)} \left[\left( \dfrac{r}{n-r}\right)^2+\left( \dfrac{r}{n-r}\right)^2+\left( \dfrac{r}{n-r}\right)\left( \dfrac{r}{n-r}\right)  \right] \\
			&=\dfrac{nr}{2} \left[3\left( \dfrac{r}{n-r}\right)^2 \right] \\
			&=\dfrac{3nr^3}{2(n-r)^2}.
		\end{array}$$
		
		$$ \begin{array}{ll}
			TRL_2(G)&=\sum\limits_{uv\in E(G)} \left[\left( \dfrac{r}{n-r}\right)^2+\left( \dfrac{r}{n-r}\right)^2-\left( \dfrac{r}{n-r}\right)\left( \dfrac{r}{n-r}\right)  \right] \\
			&=\dfrac{nr}{2} \left[\left( \dfrac{r}{n-r}\right)^2 \right] \\
			&=\dfrac{nr^3}{2(n-r)^2}.
		\end{array}$$
		
			$$ \begin{array}{ll}
			TRL_3(G)&=\sum\limits_{uv\in E(G)} \left[\left( \dfrac{r}{n-r}\right)-\left( \dfrac{r}{n-r}\right)+\left( \dfrac{r}{n-r}\right)\left( \dfrac{r}{n-r}\right)  \right] \\
			&=\dfrac{nr}{2} \left[\left( \dfrac{r}{n-r}\right)^2 \right] \\
			&=\dfrac{nr^3}{2(n-r)^2}.
		\end{array}$$
		
		$$ \begin{array}{ll}
			TRL_4(G)&=\sum\limits_{uv\in E(G)} \left|\left( \dfrac{r}{n-r}\right)-\left( \dfrac{r}{n-r}\right)\right|\left( \dfrac{r}{n-r}\right)\left( \dfrac{r}{n-r}\right)   \\
			
			&=0.
		\end{array}$$
	\end{proof}
	
	\begin{corollary}
		Let	\(C_n\)
		be a cycle with \(n \geq 3\) vertices. Then,
		\begin{enumerate}[label=\roman*.]
			\item \(TRL_1(C_n)=\dfrac{12n}{(n-2)^2}\)
			\item \(TRL_2(C_n)=\dfrac{4n}{(n-2)^2}\)
			\item \(TRL_3(C_n)=\dfrac{4n}{(n-2)^2}\)
			\item \( TRL_4(C_n)=0.\)
		\end{enumerate}
		
	\end{corollary}
	\begin{corollary}
		Let \(K_n\) be a complete graph with \(n\geq 3\)
		vertices. Then,
		
		\begin{enumerate}[label=\roman*.]
			\item \(TRL_1(K_n)=\dfrac{3n(n-1)^3}{2}\)
			\item \(TRL_2(K_n)=\dfrac{n(n-1)^3}{2}\)
			\item \(TRL_3(K_n)=\dfrac{n(n-1)^3}{2}\)
				\item \(TRL_4(K_n)=0\)
			\end{enumerate} 
	\end{corollary}

		\begin{proposition}
		Let \(P_n\) be a path with \(n\geq 3\) vertices. Then,
		
		\begin{enumerate}[label=\roman*.]
			\item \(TRL_1(P_n)= \dfrac{2\left[ 4(n-1)^2+(n-2)^2+2(n-1)(n-2)+6(n-3)(n-1)^2\right] }{(n-1)^2(n-2)^2}\).
			\item 
			\(TRL_2(P_n)=\dfrac{2\left[2(n-1)^2+(n-2)^2-2(n-1)(n-2)\right]}{(n-1)^2(n-2)^2}\)
			 \item \(TRL_3(P_n)= \dfrac{2\left[ 4(n-1)^2+(n-2)^2+2(n-1)(n-2)+6(n-3)(n-1)^2\right] }{(n-1)^2(n-2)^2}\).
			\item \(TRL_4(P_n)=\dfrac{2\left[2(n-1)^2+(n-2)^2-2(n-1)(n-2)\right]}{(n-1)^2(n-2)^2}\).
		\end{enumerate}
		
	\end{proposition}
	\begin{proof}
		The proof can be done similarly as the above.
	\end{proof}
	\begin{proposition}
		Let $K_{m,n}$ be a complete bipartite graph with
		\(1\leq m\leq n\) and \(n \geq 2\). Then 	\[TRL_1(K_{m,n})= mn(m^2+n^2+mn) \ \text{and} \ TRL_2(K_{m,n})= mn(m^2+n^2-mn). \]
	\end{proposition}
	\begin{proof}
		The proof can be done similarly as the above.\\
	\end{proof}
	
	\begin{theorem}
		Let \(W_n\) be a wheel graph. Then,
		
		\begin{enumerate}[label=\roman*.]
			\item \(TRL_1(W_n)=n\left[\dfrac{36}{(n-2)^2}+\dfrac{(n+1)n^2}{(n-2)}\right]\)
			\item \(TRL_2(W_n)=n\left[\dfrac{18}{(n-2)^2}-\dfrac{(n-5)n^2}{(n-2)}\right].\)
		\end{enumerate}

	\end{theorem}
	
	\begin{proof}
		In \(W_n\), \(E_1=\left\lbrace uv \in E(W_n)|T(u)=\dfrac{3}{n-2}=T(v)\right\rbrace \), \ \(|E_1|=n\)

		\(E_2=\left\lbrace uv \in E(W_n)|T(v_0)=n, T(v)=\dfrac{3}{n-2} \right\rbrace \), \ \(|E_2|=n\)\\
		
		$$\begin{array}{ll}
			TRL_1(W_n)&= \sum\limits_{uv\in E(G)}\left( T(u)^2+T(v)^2+T(u)T(v)\right)\\
			&=  n\left[\left( \dfrac{3}{n-2}\right)^2+ \left( \dfrac{3}{n-2}\right)^2+\left( \dfrac{3}{n-2}\right)\left( \dfrac{3}{n-2}\right) \right] +
			n\left[\left( n\right)^2+ \left( \dfrac{3}{n-2}\right)^2+n^2\left( \dfrac{3}{n-2}\right) \right]\\
			&=n\left[\dfrac{36}{(n-2)^2}+n^2 +\dfrac{3n^2}{(n-2)}\right] \\
			&=n\left[\dfrac{36}{(n-2)^2}+\dfrac{(n+1)n^2}{(n-2)}\right].
		\end{array}$$
		
		$$\begin{array}{ll}
			TRL_2(W_n)&= \sum\limits_{uv\in E(G)}\left( T(u)^2+T(v)^2-T(u)T(v)\right)\\
			&=  n\left[\left( \dfrac{3}{n-2}\right)^2+\left( \dfrac{3}{n-2}\right)^2-\left( \dfrac{3}{n-2}\right)\left( \dfrac{3}{n-2}\right) \right] +
			n\left[\left( n\right)^2+ \left( \dfrac{3}{n-2}\right)^2-n^2\left( \dfrac{3}{n-2}\right) \right]\\
			&=n\left[\dfrac{18}{(n-2)^2}+n^2 -\dfrac{3n^2}{(n-2)}\right] \\
			&=n\left[\dfrac{18}{(n-2)^2}-\dfrac{(n-5)n^2}{(n-2)}\right]. 
		\end{array}$$
		
		 $$\begin{array}{ll}
			TRL_3(W_n)&= \sum\limits_{uv\in E(G)}\left( T(u)-T(v)+T(u)T(v)\right)\\
			&=  n\left[\left( \dfrac{3}{n-2}\right)- \left( \dfrac{3}{n-2}\right)+\left( \dfrac{3}{n-2}\right)\left( \dfrac{3}{n-2}\right) \right] +
			n\left[\left( n\right)-\left( \dfrac{3}{n-2}\right)+n\left( \dfrac{3}{n-2}\right) \right]\\
			&=n\left[\dfrac{36}{(n-2)^2}+n^2 +\dfrac{3n^2}{(n-2)}\right] \\
			&=n\left[\dfrac{36}{(n-2)^2}+\dfrac{(n+1)n^2}{(n-2)}\right].
		\end{array}$$
		
		$$\begin{array}{ll}
			TRL_4(W_n)&= \sum\limits_{uv\in E(G)}\left(\left| T(u)-T(v)\right|T(u)T(v)\right)\\
			&=  n\left[\left|\left( \dfrac{3}{n-2}\right)-\left( \dfrac{3}{n-2}\right)\right|\left( \dfrac{3}{n-2}\right)\left( \dfrac{3}{n-2}\right) \right] +
			n\left[\left|\left( n\right)- \left( \dfrac{3}{n-2}\right)\right|n\left( \dfrac{3}{n-2}\right) \right]\\
			&=n\left[\dfrac{18}{(n-2)^2}+n^2 -\dfrac{3n^2}{(n-2)}\right] \\
			&=n\left[\dfrac{18}{(n-2)^2}-\dfrac{(n-5)n^2}{(n-2)}\right]. 
		\end{array}$$
	\end{proof}
	
	\begin{theorem}
		Let \(S\!f_n\) be a sunflower graph with \(3n+1\) vertices and \(5n\) edges. Then,
		\begin{enumerate}[label=\roman*.]
			\item \(TRL_1(S\!f_n)=n \left[ \dfrac{5}{(n-1)^2}+27n^2+\dfrac{8}{(3n-1)^2}+\dfrac{1}{9n^2}+\dfrac{6n}{3n-1}+\dfrac{3n}{n-1}+\dfrac{2}{(3n-1)(n-1)}+1\right]\).
			\item \(TRL_2(S\!f_n)=n\left[ \dfrac{3}{(n-1)^2}+27n^2+\dfrac{8}{(3n-1)^2}+\dfrac{1}{9n^2}-\dfrac{6n}{3n-1}-\dfrac{3n}{n-1}-\dfrac{2}{(3n-1)(n-1)}-1\right] \).
			
				\item \(TRL_3(S\!f_n)=n \left[ \dfrac{5}{(n-1)^2}+27n^2+\dfrac{8}{(3n-1)^2}+\dfrac{1}{9n^2}+\dfrac{6n}{3n-1}+\dfrac{3n}{n-1}+\dfrac{2}{(3n-1)(n-1)}+1\right]\).
			\item\(TRL_4(S\!f_n)=n\left[ \dfrac{3}{(n-1)^2}+27n^2+\dfrac{8}{(3n-1)^2}+\dfrac{1}{9n^2}-\dfrac{6n}{3n-1}-\dfrac{3n}{n-1}-\dfrac{2}{(3n-1)(n-1)}-1\right] \).
		\end{enumerate}
	\end{theorem}
	
	\begin{proof}
		In \(S\!f_n\), there are 5  types of temperature edges as follows.
		$$\begin{array}{lll}
			E_1=&\left\lbrace uv \in E(S\!f_n)|\ T(u)=\dfrac{1}{n-1}=T(v)\right\rbrace  , &|E_1|=n.\\
			E_2=& \left\lbrace uv \in E(S\!f_n)|\ T(v_0)=3n, T(v)=\dfrac{1}{n-1} \right\rbrace, & |E_2|=n.  \\
			E_3=&\left\lbrace uv \in E(S\!f_n)|\ T(u)=\dfrac{2}{3n-1}, T(v)=\dfrac{1}{n-1} \right\rbrace& |E_3|=n. \\
			E_4=&\left\lbrace uv \in E(S\!f_n)|\ T(u)=\dfrac{2}{3n-1}, T(v_0)=3n\right\rbrace &|E_4|=n.\\
			E_5=&\left\lbrace uv \in E(S\!f_n)|\ T(w)=\dfrac{1}{3n}, T(v_0)=3n\right\rbrace, &|E_5|=n.
		\end{array}$$

		
		
		

		$$\begin{array}{ll}
			TRL_1(Sf_n)&=\sum\limits_{uv\in E(G)}\left( T(u)^2+T(v)^2+T(u)T(v)\right)\\
			&=n\left[ \left( \dfrac{1}{n-1}\right)^2+\left( \dfrac{1}{n-1}\right)^2+\left( \dfrac{1}{n-1}\right)\left( \dfrac{1}{n-1}\right) \right] +n\left[ (3n)^2+\left( \dfrac{1}{n-1}\right)^2+3n\left( \dfrac{1}{n-1}\right) \right]\\&+n\left[\left( \dfrac{2}{3n-1}\right)^2+\left( \dfrac{1}{n-1}\right)^2+\left( \dfrac{2}{3n-1}\right)\left( \dfrac{1}{n-1}\right) \right]+n\left[\left( \dfrac{2}{3n-1}\right)^2+(3n)^2+\left( \dfrac{2}{3n-1}\right)(3n) \right] \\&+\left[ \left( \dfrac{1}{3n}\right)^2+(3n)^2+\left( \dfrac{1}{3n}\right)\cdot (3n)\right] \\
			&= n \left[ \dfrac{5}{(n-1)^2}+27n^2+\dfrac{8}{(3n-1)^2}+\dfrac{1}{9n^2}+\dfrac{6n}{3n-1}+\dfrac{3n}{n-1}+\dfrac{2}{(3n-1)(n-1)}+1\right] .
			
		\end{array}$$
		
		$$\begin{array}{ll}
			TRL_1(S\!f_n)&=\sum\limits_{uv\in E(G)}\left( T(u)^2+T(v)^2-T(u)T(v)\right)\\
			&=n\left[ \left( \dfrac{1}{n-1}\right)^2+\left( \dfrac{1}{n-1}\right)^2-\left( \dfrac{1}{n-1}\right)\left( \dfrac{1}{n-1}\right) \right] +n\left[ (3n)^2+\left( \dfrac{1}{n-1}\right)^2-3n\left( \dfrac{1}{n-1}\right) \right]\\&+n\left[\left( \dfrac{2}{3n-1}\right)^2+\left( \dfrac{1}{n-1}\right)^2-\left( \dfrac{2}{3n-1}\right)\left( \dfrac{1}{n-1}\right) \right]+n\left[\left( \dfrac{2}{3n-1}\right)^2+(3n)^2-\left( \dfrac{2}{3n-1}\right)(3n) \right] \\&+\left[ \left( \dfrac{1}{3n}\right)^2+(3n)^2-\left( \dfrac{1}{3n}\right)\cdot (3n)\right] \\
			&=  n\left[ \dfrac{3}{(n-1)^2}+27n^2+\dfrac{8}{(3n-1)^2}+\dfrac{1}{9n^2}-\dfrac{6n}{3n-1}-\dfrac{3n}{n-1}-\dfrac{2}{(3n-1)(n-1)}-1\right].
			
		\end{array}$$
		
		 $$\begin{array}{ll}
			TRL_3(Sf_n)&=\sum\limits_{uv\in E(G)}\left(  T(u)-T(v)+T(u)T(v)\right)\\
			&=n\left[ \left( \dfrac{1}{n-1}\right)-\left( \dfrac{1}{n-1}\right)+\left( \dfrac{1}{n-1}\right)\left( \dfrac{1}{n-1}\right) \right] +n\left[ (3n)^2+\left( \dfrac{1}{n-1}\right)-3n\left( \dfrac{1}{n-1}\right) \right]\\&+n\left[\left( \dfrac{2}{3n-1}\right)-\left( \dfrac{1}{n-1}\right)+\left( \dfrac{2}{3n-1}\right)\left( \dfrac{1}{n-1}\right) \right]+n\left[\left( \dfrac{2}{3n-1}\right)-(3n)+\left( \dfrac{2}{3n-1}\right)(3n) \right] \\&+\left[ \left( \dfrac{1}{3n}\right)-(3n)+\left( \dfrac{1}{3n}\right)\cdot (3n)\right] \\
			&= n \left[ \dfrac{5}{(n-1)^2}+27n^2+\dfrac{8}{(3n-1)^2}+\dfrac{1}{9n^2}+\dfrac{6n}{3n-1}+\dfrac{3n}{n-1}+\dfrac{2}{(3n-1)(n-1)}+1\right] .
			
		\end{array}$$
		
		$$\begin{array}{ll}
			TRL_4(S\!f_n)&=\sum\limits_{uv\in E(G)}\left| T(u)-T(v)\right|T(u)T(v)\\
			&=n\left[ \left|\left( \dfrac{1}{n-1}\right)-\left( \dfrac{1}{n-1}\right)\right|\left( \dfrac{1}{n-1}\right)\left( \dfrac{1}{n-1}\right) \right] +n\left[\left| (3n)-\left( \dfrac{1}{n-1}\right)\right|3n\left( \dfrac{1}{n-1}\right) \right]\\&+n\left[\left|\left( \dfrac{2}{3n-1}\right)-\left( \dfrac{1}{n-1}\right)\right|\left( \dfrac{2}{3n-1}\right)\left( \dfrac{1}{n-1}\right) \right]+n\left[\left|\left( \dfrac{2}{3n-1}\right)-(3n)\right|\left( \dfrac{2}{3n-1}\right)(3n) \right] \\&+\left[ \left|\left( \dfrac{1}{3n}\right)-(3n)\right|\left( \dfrac{1}{3n}\right)\cdot (3n)\right] \\
			&=  n\left[ \dfrac{3}{(n-1)^2}+27n^2+\dfrac{8}{(3n-1)^2}+\dfrac{1}{9n^2}-\dfrac{6n}{3n-1}-\dfrac{3n}{n-1}-\dfrac{2}{(3n-1)(n-1)}-1\right].
			
		\end{array}$$
	\end{proof}
	The results for the other \(RL\) indices/exponentials presented in the Table \ref{tablee1}, \ref{tablee2}, \ref{tablee3}, \ref{tablee4}  can be computed similarly as above.
	\subsection{KV indices}
	Consider $G$ to be a finite, simple linked graph with a vertex $V(G)$ and an edge set $E(G)$. The number of vertices adjacent to a vertex $v$ in $G$ determines its degree $d(v)$. Let the product of the degrees of all vertices adjacent to a vertex $v$ be denoted by $M_G(v)$.

	Kulli introduced the first, second, third and fourth KV indices in \cite{article51} , defined as, 
	
	\[KV_1(G)=\sum_{uv\in E(G)}^{}M_G(u)+M_G(v) \ \text{and} \ KV_2(G)=\sum_{uv\in E(G)}^{}M_G(u)M_G(v) . \]

	Inspired by these indices, here we introduce the first, second, third and fourth Rehan-Lanel KV indices as follows.
	
	\begin{definition}
		
		The first and second Rehan- Lanel KV indices of a graph \(G\) are defined as, 
		\begin{enumerate}[label=\roman*.]
			\item 	\(RLKV_1(G)=\sum\limits_{uv\in E(G)}\left(  M_G^2(u)+M_G^2(v)+M_G(u)M_G(v)\right)\).
			\item 	\(RLKV_2(G)=\sum\limits_{uv\in E(G)}\left(  M_G^2(u)+M_G^2(v)-M_G(u)M_G(v)\right)  \).
			\item 	\(RLKV_3(G)=\sum\limits_{uv\in E(G)}\left(  M_G(u)-M_G(v)+M_G(u)M_G(v)\right)\).
			\item 	\(RLKV_4(G)=\sum\limits_{uv\in E(G)}\left(  |M_G(u)-M_G(v)|M_G(u)M_G(v)\right)  \).
		\end{enumerate}

	\end{definition}
	
	In addition, we define few more KV indices/exponentials as in the table \ref{tablef1} and \ref{tablef2}.
	
	\begin{longtable}[c]{| c | l | l |}
		\caption{Rehan- Lanel KV indices/exponentials \label{tablef1}}\\

		\hline
		\multicolumn{1}{|c}{\textbf{S.N.}} & 
		\multicolumn{1}{|c}{\textbf{Name of the index}} & \multicolumn{1}{|c|}{\textbf{Formula}} \\\hline
		
		\endfirsthead
		
		\hline
		\endlastfoot

		1&		\multirow{2}{4cm}{The first hyper  RL KV index} & \(HRLKV_1(G)=\sum\limits_{uv\in E(G)}\left(  M_G^2(u)+M_G^2(v)+M_G(u)M_G(v)\right)^2 \)\\[5ex]
		
		2&		\multirow{2}{4cm}{The second hyper  RL KV index} & \(HRLKV_2(G)=\sum\limits_{uv\in E(G)}\left(  M_G^2(u)+M_G^2(v)-M_G(u)M_G(v)\right)^2 \)\\[5ex]
		3&	\multirow{2}{4cm}{The first  inverse  RL KV index} & \(IRLKV_1(G)=\sum\limits_{uv\in E(G)}\dfrac{1}{\left(  M_G^2(u)+M_G^2(v)+M_G(u)M_G(v)\right)} \)\\[6ex]
		
		4&	\multirow{2}{5cm}{The second  inverse RL KV index} & \(IRLKV_2(G)=\sum\limits_{uv\in E(G)}\dfrac{1}{\left(  M_G^2(u)+M_G^2(v)-M_G(u)M_G(v)\right)} \)\\[6ex]
		
		5&	\multirow{2}{5cm}{The first general  RL index} & \(GRLKV_1(G)=\sum\limits_{uv\in E(G)}\left(M_G^2(u)+M_G^2(v)+M_G(u)M_G(v)\right)^a \)\\[6ex]
		6&	\multirow{2}{5cm}{The second  general  RL index} & \(GRLKV_2(G)=\sum\limits_{uv\in E(G)}\left( M_G^2(u)+M_G^2(v)-M_G(u)M_G(v)\right)^a \)\\[6ex]
		7&	\multirow{2}{5cm}{The first  RL KV exponential} & \(RLKV_1(G,x)=\sum\limits_{uv\in E(G)}x^{\left(M_G^2(u)+M_G^2(v)+M_G(u)M_G(v)\right)} \)\\[3ex]
		
		8&	\multirow{2}{5cm}{The second  RL KV exponential} & \(RLKV_2(G,x)=\sum\limits_{uv\in E(G)}x^{\left(  M_G^2(u)+M_G^2(v)-M_G(u)M_G(v)\right)} \)\\[3ex]
		9&	\multirow{2}{5cm}{The first  hyper  RL KV exponential} & \(HRLKV_1(G,x)=\sum\limits_{uv\in E(G)}x^{\left(  M_G^2(u)+M_G^2(v)+M_G(u)M_G(v)\right)^2 }\)\\[5ex]
		
		10&	\multirow{2}{5cm}{The second  hyper RL KV exponential} & \(HRLKV_2(G,x)=\sum\limits_{uv\in E(G)}x^{\left( M_G^2(u)+M_G^2(v)-M_G(u)M_G(v)\right)^2} \)\\[3ex]
		
		11&	\multirow{2}{5cm}{The first  inverse  RLKV exponential} & \(IRLKV_1(G,x )=\sum\limits_{uv\in E(G)}x^{\dfrac{1}{\left( M_G^2(u)+M_G^2(v)+M_G(u)M_G(v)\right)}} \)\\[3ex]
		
		12&	\multirow{2}{5cm}{The second inverse  RLKV exponential} & \(IRLKV_2(G,x)=\sum\limits_{uv\in E(G)}x^{\dfrac{1}{\left(  M_G^2(u)+M_G^2(v)-M_G(u)M_G(v)\right)} }\)\\[5ex]
		
		13&	\multirow{2}{5cm}{The first general  RL KV exponential} & \(GRLKV_1(G,x)=\sum\limits_{uv\in E(G)}x^{\left(  M_G^2(u)+M_G^2(v)+M_G(u)M_G(v)\right)^a} \)\\[5ex]
		
		14&	\multirow{2}{5cm}{The second  general RL KV exponential} & \(GRLKV_2(G,x)=\sum\limits_{uv\in E(G)}x^{\left(  M_G^2(u)+M_G^2(v)-M_G(u)M_G(v)\right)^a} \)\\[6ex]

	\end{longtable}
	\begin{longtable}[c]{| c | l | l |}
		\caption{Multiplicative Rehan- Lanel KV indices/RLKV indices \label{tablef2}}\\

		\hline
		\multicolumn{1}{|c}{\textbf{S.N.}} & 
		\multicolumn{1}{|c}{\textbf{Name of the index}} & \multicolumn{1}{|c|}{\textbf{Formula}} \\\hline
		
		\endfirsthead
		
		\hline
		\endlastfoot
		
		1&	\multirow{2}{5cm}{The first multiplicative  RL KV index }& \(MRLKV_1(G)=\prod\limits_{uv\in E(G)}\left(M_G^2(u)+M_G^2(v)+M_G(u)M_G(v)\right) \)\\[3ex]
		2&	\multirow{2}{5cm}{The second multiplicative  RL KV index }& \(MRLKV_2(G)=\prod\limits_{uv\in E(G)}\left( M_G^2(u)+M_G^2(v)-M_G(u)M_G(v)\right) \)\\[3ex]
		
		3&		\multirow{2}{4cm}{The first  multiplicative hyper RL KV index} & \(MHRLKV_1(G)=\prod\limits_{uv\in E(G)}\left(  M_G^2(u)+M_G^2(v)+M_G(u)M_G(v)\right)^2 \)\\[5ex]
		
		4&		\multirow{2}{4cm}{The second multiplicative hyper  RL KV index} & \(MHRLKV_2(G)=\prod\limits_{uv\in E(G)}\left(  M_G^2(u)+M_G^2(v)-M_G(u)M_G(v)\right)^2 \)\\[5ex]
		
		5&	\multirow{2}{4cm}{The first multiplicative inverse  RL KV index} & \(MIRLKV_1(G)=\prod\limits_{uv\in E(G)}\dfrac{1}{\left(  M_G^2(u)+M_G^2(v)+M_G(u)M_G(v)\right)} \)\\[6ex]
		
		6&	\multirow{2}{5cm}{The second multiplicative inverse  RL KV index} & \(MIRLKV_2(G)=\prod\limits_{uv\in E(G)}\dfrac{1}{\left(  M_G^2(u)+M_G^2(v)-M_G(u)M_G(v)\right)} \)\\[6ex]
		
		7&	\multirow{2}{5cm}{The first general multiplicative  RL KV index} & \(GRLKV_1(G)=\prod\limits_{uv\in E(G)}\left( M_G^2(u)+M_G^2(v)+M_G(u)M_G(v)\right)^a \)\\[6ex]
		
		8&	\multirow{2}{5cm}{The second  general multiplicative  RL KV index} & \(MGRLKV_2(G)=\prod\limits_{uv\in E(G)}\left( M_G^2(u)+M_G^2(v)-M_G(u)M_G(v)\right) \)\\[6ex]
		9&	\multirow{2}{5cm}{The first multiplicative  RL KV exponential} & \(BRL_1(G,x)=\prod\limits_{uv\in E(G)}x^{\left(M_G^2(u)+M_G^2(v)+M_G(u)M_G(v)\right)} \)\\[3ex]
		10&	\multirow{2}{5cm}{The second multiplicative  RL KV exponential} & \(MRLKV_2(G,x)=\prod\limits_{uv\in E(G)}x^{\left(  M_G^2(u)+M_G^2(v)-M_G(u)M_G(v)\right)} \)\\[3ex]
		
		11&	\multirow{2}{5cm}{The first multiplicative hyper RL KV exponential} & \(MHRLKV_1(G,x)=\prod\limits_{uv\in E(G)}x^{\left(  M_G^2(u)+M_G^2(v)+M_G(u)M_G(v)\right)^2 }\)\\[5ex]
		
		12&	\multirow{2}{5cm}{The second multiplicative hyper  RL KV exponential} & \(MHRLKV_2(G,x)=\prod\limits_{uv\in E(G)}x^{\left( M_G^2(u)+M_G^2(v)-M_G(u)M_G(v)\right)^2} \)\\[3ex]
		13&	\multirow{2}{5cm}{The first multiplicative inverse  RL KV exponential} & \(IBRL_1(G,x )=\prod\limits_{uv\in E(G)}x^{\dfrac{1}{\left(  M_G^2(u)+M_G^2(v)+M_G(u)M_G(v)\right)}} \)\\[3ex]
		
		14&	\multirow{2}{5cm}{The second multiplicative inverse  RL KV exponential} & \(MIRLKV_2(G,x)=\prod\limits_{uv\in E(G)}x^{\dfrac{1}{\left(  M_G^2(u)+M_G^2(v)-M_G(u)M_G(v)\right)} }\)\\[5ex]
		
		15&	\multirow{2}{5cm}{The first multiplicative general  RL KV exponential} & \(MGRLKV_1(G,x)=\prod\limits_{uv\in E(G)}x^{\left( M_G^2(u)+M_G^2(v)+M_G(u)M_G(v)\right)^a} \)\\[5ex]

		16&	\multirow{2}{5cm}{The second multiplicative general  RL KV exponential} & \(MGRLKV_2(G,x)=\prod\limits_{uv\in E(G)}x^{\left(  M_G^2(u)+M_G^2(v)-M_G(u)M_G(v)\right)^a} \)\\[6ex]

	\end{longtable}
	
		\begin{longtable}[c]{| c | l | l |}
		\caption{Rehan- Lanel KV indices/exponentials \label{tablef3}}\\

		\hline
		\multicolumn{1}{|c}{\textbf{S.N.}} & 
		\multicolumn{1}{|c}{\textbf{Name of the index}} & \multicolumn{1}{|c|}{\textbf{Formula}} \\\hline
		
		\endfirsthead
		
		\hline
		\endlastfoot

		1&		\multirow{2}{4cm}{The third hyper  RL KV index} & \(HRLKV_1(G)=\sum\limits_{uv\in E(G)}\left(  M_G(u)-M_G(v)+M_G(u)M_G(v)\right)^2 \)\\[5ex]
		
		2&		\multirow{2}{4cm}{The fourth hyper  RL KV index} & \(HRLKV_2(G)=\sum\limits_{uv\in E(G)}\left( |M_G(u)-M_G(v)|M_G(u)M_G(v)\right)^2 \)\\[3ex]
		
		3&	\multirow{2}{4cm}{The third  inverse  RL KV index} & \(IRLKV_1(G)=\sum\limits_{uv\in E(G)}\dfrac{1}{\left(  M_G(u)-M_G(v)+M_G(u)M_G(v)\right)} \)\\[3ex]
		4&	\multirow{1}{5cm}{The fourth  inverse RL KV index} & \(IRLKV_2(G)=\sum\limits_{uv\in E(G)}\dfrac{1}{\left(  |M_G(u)-M_G(v)|M_G(u)M_G(v)\right)} \)\\[6ex]
		5&	\multirow{1}{5cm}{The third general  RL index} & \(GRLKV_1(G)=\sum\limits_{uv\in E(G)}\left(M_G(u)-M_G(v)+M_G(u)M_G(v)\right)^a \)\\[6ex]
		6&	\multirow{1}{5cm}{The fourth  general  RL index} & \(GRLKV_2(G)=\prod\limits_{uv\in E(G)}\left( |M_G(u)-M_G(v)|M_G(u)M_G(v)\right)^a \)\\[6ex]
		7&	\multirow{1}{5cm}{The third  RL KV exponential} & \(RLKV_1(G,x)=\sum\limits_{uv\in E(G)}x^{\left(M_G(u)-M_G(v)+M_G(u)M_G(v)\right)} \)\\[3ex]
		
		8&	\multirow{1}{5cm}{The fourth  RL KV exponential} & \(RLKV_2(G,x)=\sum\limits_{uv\in E(G)}x^{\left(  |M_G(u)-M_G(v)|M_G(u)M_G(v)\right)} \)\\[3ex]
		9&	\multirow{2}{5cm}{The third  hyper  RL KV exponential} & \(HRLKV_1(G,x)=\sum\limits_{uv\in E(G)}x^{\left(  M_G(u)-M_G(v)+M_G(u)M_G(v)\right)^2 }\)\\[5ex]
		
		10&	\multirow{2}{5cm}{The fourth  hyper RL KV exponential} & \(HRLKV_2(G,x)=\sum\limits_{uv\in E(G)}x^{\left( |M_G(u)-M_G(v)|M_G(u)M_G(v)\right)^2} \)\\[3ex]
		
		11&	\multirow{2}{5cm}{The third  inverse  RLKV exponential} & \(IRLKV_1(G,x )=\sum\limits_{uv\in E(G)}x^{\dfrac{1}{\left( M_G(u)-M_G(v)+M_G(u)M_G(v)\right)}} \)\\[3ex]
		
		12&	\multirow{2}{5cm}{The fourth inverse  RLKV exponential} & \(IRLKV_2(G,x)=\sum\limits_{uv\in E(G)}x^{\dfrac{1}{\left(  |M_G(u)-M_G(v)|M_G(u)M_G(v)\right)} }\)\\[5ex]
		
		13&	\multirow{2}{5cm}{The third general  RL KV exponential} & \(GRLKV_1(G,x)=\sum\limits_{uv\in E(G)}x^{\left(  M_G(u)-M_G(v)+M_G(u)M_G(v)\right)^a} \)\\[5ex]
		
		14&	\multirow{2}{5cm}{The fourth  general RL KV exponential} & \(GRLKV_2(G,x)=\sum\limits_{uv\in E(G)}x^{\left( |M_G(u)-M_G(v)|M_G(u)M_G(v)\right)^a} \)\\[6ex]

	\end{longtable}
	
	\begin{longtable}[c]{| c | l | l |}
		\caption{Multiplicative Rehan- Lanel KV indices/RLKV indices \label{tablef4}}\\

		\hline
		\multicolumn{1}{|c}{\textbf{S.N.}} & 
		\multicolumn{1}{|c}{\textbf{Name of the index}} & \multicolumn{1}{|c|}{\textbf{Formula}} \\\hline
		
		\endfirsthead
		
		\hline
		\endlastfoot
		
		1&	\multirow{2}{5cm}{The third multiplicative  RL KV index }& \(MRLKV_1(G)=\prod\limits_{uv\in E(G)}\left(M_G(u)-M_G(v)+M_G(u)M_G(v)\right) \)\\[3ex]
		
		2&	\multirow{2}{5cm}{The fourth multiplicative  RL KV index }& \(MRLKV_2(G)=\prod\limits_{uv\in E(G)}\left( |M_G(u)-M_G(v)|M_G(u)M_G(v)\right) \)\\[3ex]
		\hline 
		3&		\multirow{2}{4cm}{The third  multiplicative hyper RL KV index} & \(MHRLKV_1(G)=\prod\limits_{uv\in E(G)}\left(  M_G(u)-M_G(v)+M_G(u)M_G(v)\right)^2 \)\\[5ex]
		
		4&		\multirow{2}{4cm}{The fourth multiplicative hyper  RL KV index} & \(MHRLKV_2(G)=\prod\limits_{uv\in E(G)}\left(  |M_G(u)-M_G(v)|M_G(u)M_G(v)\right)^2 \)\\[5ex]
		
		5&	\multirow{2}{4cm}{The third multiplicative inverse  RL KV index} & \(MIRLKV_1(G)=\prod\limits_{uv\in E(G)}\dfrac{1}{\left(  M_G(u)-M_G(v)+M_G(u)M_G(v)\right)} \)\\[6ex]
		6&	\multirow{2}{5cm}{The fourth multiplicative inverse  RL KV index} & \(MIRLKV_2(G)=\prod\limits_{uv\in E(G)}\dfrac{1}{\left(  |M_G(u)-M_G(v)|M_G(u)M_G(v)\right)} \)\\[3ex]
		7&	\multirow{2}{5cm}{The fourth general multiplicative  RL KV index} & \(GRLKV_1(G)=\prod\limits_{uv\in E(G)}\left( M_G(u)-M_G(v)+M_G(u)M_G(v)\right)^a \)\\[6ex]
		
		8&	\multirow{2}{5cm}{The fourth  general multiplicative  RL KV index} & \(MGRLKV_2(G)=\prod\limits_{uv\in E(G)}\left( |M_G(u)-M_G(v)|M_G(u)M_G(v)\right) \)\\[6ex]
		9&	\multirow{2}{5cm}{The third multiplicative  RL KV exponential} & \(BRL_1(G,x)=\prod\limits_{uv\in E(G)}x^{\left(M_G(u)-M_G(v)+M_G(u)M_G(v)\right)} \)\\[3ex]
		10&	\multirow{2}{5cm}{The fourth multiplicative  RL KV exponential} & \(MRLKV_2(G,x)=\prod\limits_{uv\in E(G)}x^{\left(  |M_G(u)-M_G(v)|M_G(u)M_G(v)\right)} \)\\[3ex]
		
		11&	\multirow{2}{5cm}{The third multiplicative hyper RL KV exponential} & \(MHRLKV_1(G,x)=\prod\limits_{uv\in E(G)}x^{\left(  M_G(u)-M_G(v)+M_G(u)M_G(v)\right)^2 }\)\\[5ex]
		
		12&	\multirow{2}{5cm}{The fourth multiplicative hyper  RL KV exponential} & \(MHRLKV_2(G,x)=\prod\limits_{uv\in E(G)}x^{\left( |M_G(u)-M_G(v)|M_G(u)M_G(v)\right)^2} \)\\[3ex]
		13&	\multirow{2}{5cm}{The third multiplicative inverse  RL KV exponential} & \(IBRL_1(G,x )=\prod\limits_{uv\in E(G)}x^{\dfrac{1}{\left(  M_G(u)-M_G(v)+M_G(u)M_G(v)\right)}} \)\\[3ex]
		
		14&	\multirow{2}{5cm}{The fourth multiplicative inverse  RL KV exponential} & \(MIRLKV_2(G,x)=\prod\limits_{uv\in E(G)}x^{\dfrac{1}{\left(  |M_G(u)-M_G(v)|M_G(u)M_G(v)\right)} }\)\\[5ex]
		
		15&	\multirow{2}{5cm}{The third multiplicative general  RL KV exponential} & \(MGRLKV_1(G,x)=\prod\limits_{uv\in E(G)}x^{\left( M_G(u)-M_G(v)+M_G(u)M_G(v)\right)^a} \)\\[5ex]

		16&	\multirow{2}{5cm}{The fourth multiplicative general  RL KV exponential} & \(MGRLKV_2(G,x)=\prod\limits_{uv\in E(G)}x^{\left(|M_G(u)-M_G(v)|M_G(u)M_G(v)\right)^a} \)\\[6ex]

	\end{longtable}
	We calculate the first, second, third and fourth Rehan-Lanel KV indices (\(RLKV_1\), \(RLKV_2\), \(RLKV_3\), \(RLKV_4\)) for the standard graphs such as regular graph, a cycle, complete graph, a path and a complete bipartite graph. Furthermore, We determine (\(RLKV_1\), \(RLKV_2\), \(RLKV_3\), \(RLKV_4\)) for wheel graph and the sunflower graph.
	\begin{proposition}
		Let \(G\) is an $r$-regular graph with \(n\) vertices and
		\(r \geq 2\). Then,
		\begin{enumerate}[label=\roman*.]
			\item \(RLKV_1(G)=\dfrac{3nr^{3r}}{2}\)
			\item \(RLKV_2(G)=\dfrac{nr^{3r}}{2}\)
			\item \(RLKV_3(G)=\dfrac{3nr^{2r+1}}{2}\)
			\item \(RLKV_4(G)=0.\)
		\end{enumerate}

	\end{proposition}
	
	\begin{proof}
		Let \(G\) is an $r-$regular graph with \(n\) vertices and
		\(r \geq 2\).\\
		Now	\(M_G(u)=M_G(v)=r^r\).
		
		Then \(G\) has \(\dfrac{nr}{2}\)	edges. 
		From definition we have
		$$ \begin{array}{ll}
			RLKV_1(G)&=\sum\limits_{uv\in E(G)} (r^{2r}+r^{2r}+r^r\cdot r^r)\\
			&=\dfrac{nr}{2} (3r^{2r})\\
			&=\dfrac{3nr^{3r}}{2}.
		\end{array}$$
		
		$$ \begin{array}{ll}
			RL_2(G)&=\sum\limits_{uv\in E(G)} (r^{2r}+r^{2r}-r^r\cdot r^r)\\
			&=\dfrac{nr}{2} (r^{2r})\\
			&=\dfrac{nr^{3r}}{2}
		\end{array}$$
		$$ \begin{array}{ll}
			RLKV_3(G)&=\sum\limits_{uv\in E(G)} (r^r-r^r+r^r\cdot r^r)\\
			&=\dfrac{nr}{2} (r^{2r})\\
			&=\dfrac{nr^{2r+1}}{2}.
		\end{array}$$
		
			$$ \begin{array}{ll}
			RLKV_4(G)&=\sum\limits_{uv\in E(G)} (|r^r-r^r|r^r\cdot r^r)\\
			&=0.\\
			
		\end{array}$$
	\end{proof}
	\begin{corollary}
		Let	\(C_n\)
		be a cycle with \(n \geq 3\) vertices. Then, 
		\begin{enumerate}[label=\roman*.]
			\item \(RLKV_1(C_n)=96n\)
			\item \(RLKV_2(C_n)=16n\)
			\item \(RLKV_3(C_n)=16n\)
			\item \(RLKV_3(C_n)=0\).
		\end{enumerate}

	\end{corollary}
	\begin{corollary}
		Let \(K_n\) be a complete graph with \(n\geq 3\)
		vertices. Then,
		\begin{enumerate}[label=\roman*.]
			\item \(RLKV_1(K_n)=\dfrac{3n(n-1)^{n-1}}{2}\)
			\item \(RLKV_2(K_n)=\dfrac{n(n-1)^{n-1}}{2}\)
			\item  \(RLKV_3(G)=\dfrac{n(n-1)^{2(n-1)+1}}{2}\)
			\item \(RLKV_4(G)=0\)
		\end{enumerate}

	\end{corollary}
	
	\begin{proposition}
		Let \(P_n\) be a path with \(n\geq 3\) vertices. Then
		\begin{enumerate}[label=\roman*.]
			\item \(RLKV_1(P_n)= 24(2n-5)\)
			\item \(RLKV_2(P_n)=8(2n-5)\)
			\item \(RLKV_3(P_n)= 16n-40 \)
			\item \(RLKV_4(P_n)=0.\)
		\end{enumerate}
	
	\end{proposition}
	
	\begin{proof}
		The proof can be done similarly as above.
	\end{proof}
	\begin{proposition}
		Let $K_{m,n}$ be a complete bipartite graph with
		\(1\leq m\leq n\) and \(n \geq 2\). Then,
		
		\begin{enumerate}[label=\roman*.]
			\item \(RLKV_1(K_{m,n})= mn(m^{2n}+n^{2m}+m^nn^m)\)
			\item \(RLKV_2(K_{m,n})= mn(m^{2n}+n^{2m}-m^nn^m)\)
			\item \(RLKV_3(K_{m,n})= mn(m^{n}-n^{m}+m^nn^m)\)
			\item \(RLKV_4(K_{m,n})= |m^{n}-n^{m}|m^{n+1}n^{m+1}.\)
		\end{enumerate}	
		
	\end{proposition}
	\begin{proof}
		The proof can be done similarly as above.
	\end{proof}
	\begin{theorem}
		Let \(W_n\) be a wheel graph. Then,
		\begin{enumerate}[label=\roman*.]
			\item \(RLKV_1(W_n)=n(324n^2+3^n(3^n+9n))\)
			\item \(RLKV_2(W_n)=n(162n^2+3^n(3^n-9n))\)
			\item \(RLKV_3(W_n)=n(81n^2+(9n+1)3^n-9n)\)
			\item \(RLKV_4(W_n)=n(|3^n-9n|\cdot 3^n \cdot 9n)\).
		\end{enumerate}

	\end{theorem}
	\begin{proof}
		In \(W_n\), \(E_1=\left\lbrace uv \in E(G)| M_G(v)=9n=M_G(v)\right\rbrace \), \(|E_1|=n\).\\

		\(E_2=\left\lbrace uv \in E(G)|M_G(v_0)=3^n, M_G(v)=9n \right\rbrace \), \(|E_2|=n\).
		
		$$\begin{array}{ll}
			RLKV_1(W_n)&=\sum\limits_{uv\in E(W_n)}\left( M^2_G(u)+M^2_G(v)+M_G(u)M_G(v)\right)\\
			&=n\left( (9n)^2+(9n)^2+9n.9n \right) +n\left( (3^n)^2+(9n)^2+3^n.9n\right) \\
			&=n(324n^2+3^n(3^n+9n)).
		\end{array}$$

		$$\begin{array}{ll}
			RLKV_2(W_n)&=\sum\limits_{uv\in E(W_n)}\left( M^2_G(u)+M^2_G(v)-M_G(u)M_G(v)\right)\\
			&=n\left( (9n)^2+(9n)^2-9n.9n \right) +n\left( (3^n)^2+(9n)^2-3^n.9n\right) \\
			&=n(162n^2+3^n(3^n-9n)).
		\end{array}$$
		
			$$\begin{array}{ll}
			RLKV_3(W_n)&=\sum\limits_{uv\in E(W_n)}\left( M_G(u)-M_G(v)+M_G(u)M_G(v)\right)\\
			&=n\left( (9n)-(9n)+9n.9n \right) +n\left( (3^n)-(9n)+3^n.9n\right) \\
			&=n(81n^2+(9n+1)3^n-9n).
		\end{array}$$

		$$\begin{array}{ll}
			RLKV_4(W_n)&=\sum\limits_{uv\in E(W_n)}\left| M_G(u)+M_G(v)\right|M_G(u)M_G(v)\\
			&=n\left( |(9n)-(9n)|9n.9n \right) +n\left( |(3^n)-(9n)|3^n.9n\right) \\
			&=n(|3^n-9n|\cdot 3^n \cdot 9n).
		\end{array}$$
	\end{proof}

	\begin{theorem}
		Let \(S\!f_n\) be a sunflower graph. Then,
		\begin{enumerate}[label=\roman*.]
			\item  \(RLKV_1(S\!f_n)=n\left[ 47520n^2+3 \cdot 2^{6n}+111n\cdot 2^{3n}\right] \)
			\item \(RLKV_2(S\!f_n)=n\left[26793n^2+3 \cdot 2^{6n}-111n\cdot 2^{3n}\right]. \)
			
		\end{enumerate}

	\end{theorem}
	
	\begin{proof}
		In \(S\!f_n\), there are five types of KV edges,

		There are five types of KV edges of a sunflower graph as follows.
		$$\begin{array}{lll}
			E_1&=\left\lbrace uv \in E(S\!f_n)| M_G(v)=M_G(v)=2^5 \cdot 3n \right\rbrace, &|E_1|=n.\\
			E_2&=\left\lbrace uv \in E(S\!f_n)| M_G(v_0)=2^{3n} ,M_G(v)=2^5\cdot 3n \right\rbrace ,& |E_2|=n.\\
			E_3&=\left\lbrace uv \in E(S\!f_n)| M_G(u)=2^2\cdot 3n, M_G(v)=2^5 \cdot 3n \right\rbrace,& |E_3|=n.\\
			E_4&=\left\lbrace uv \in E(Sf_n)| M_G(v_0)=2^{3n} ,M_G(u)= 2^2 \cdot 3n \right\rbrace , &|E_4|=n.\\
			E_5&=\left\lbrace uv \in E(S\!f_n)| M_G(w)=3n ,M_G(v_0)=2^{3n} \right\rbrace ,& |E_5|=n.
			
		\end{array}$$

		$$\begin{array}{ll}
			RLKV_1(S\!f_n)&=\sum\limits_{uv\in E(Sf_n)}\left( M^2_G(u)+M^2_G(v)+M_G(u)M_G(v)\right)\\
			&=n\left( (2^5 \cdot 3n)^2+(2^5 \cdot 3n)^2+\left( 2^5 \cdot 3n\right) \left(2^5 \cdot 3n \right)  \right) +n\left( (2^{3n})^2+(2^5\cdot 3n)^2+2^{3n}\cdot 2^5 \cdot 3n\right) \\ & +n\left((2^2\cdot 3n)^2+(2^5\cdot 3n)^2 +2^2\cdot 3n \cdot 2^5 \cdot 3n\right)+n\left( (2^{3n})^2+(2^2\cdot 3n)^2+(2^{3n}\cdot 2^2 \cdot 3n )\right)  \\ & +n\left( (3n)^2+(2^{3n})^2+3n \cdot 2^{3n}\right) \\
			&=n(47520n^2+3 \cdot 2^{6n}+111n\cdot 2^{3n}).
		\end{array}$$

		$$\begin{array}{ll}
			RLKV_2(S\!f_n)&=\sum\limits_{uv\in E(Sf_n)}\left( M^2_G(u)+M^2_G(v)-M_G(u)M_G(v)\right)\\
			&=n\left( (2^5 \cdot 3n)^2+(2^5 \cdot 3n)^2-\left( 2^5 \cdot 3n\right) \left(2^5 \cdot 3n \right)  \right) +n\left( (2^{3n})^2+(2^5\cdot 3n)^2-2^{3n}\cdot 2^5 \cdot 3n\right) \\ & +n\left((2^2\cdot 3n)^2+(2^5\cdot 3n)^2 -2^2\cdot 3n \cdot 2^5 \cdot 3n\right)+n\left( (2^{3n})^2+(2^2\cdot 3n)^2-(2^{3n}\cdot 2^2 \cdot 3n )\right)  \\ & +n\left( (3n)^2+(2^{3n})^2-3n \cdot 2^{3n}\right) \\
			&=n(26793n^2+3 \cdot 2^{6n}-111n\cdot 2^{3n}).
		\end{array}$$
		
%
%
	\end{proof}
	
	Similarly, we can compute the results for the other indices/exponentials listed in the above pair of tables.
	\subsection{Nbd indices}

	Let $S_G(u)$ denote the sum of the degrees of all neighborhood vertices of a vertex $u$ in $G$.
	
	\[S_G(u)=\sum_{u\in N_G(v)} d_G(u)\]
	
	In \cite{article40}	Graovac et al. defined the following indices:
	
	\[N_1(G)=\sum\limits_{uv\in E(G)}(S_G(u)+S_G(v))\]
	\[N_2(G)=\sum\limits_{uv\in E(G)}(S_G(u)S_G(v)).\]
	Motivated by these definitions, here we introduce the neighborhood version of Rehan- Lanel indices indices.
	
	\begin{definition}
		The first and second Neighborhood Rehan-Lanel indices of a graph \(G\) are defined by,
		\begin{enumerate}[label=\roman*.]
			\item \(NRL_1(G)=\sum_{uv\in E(G)}(S_G^2(u)+S_G^2(v)+S_G(u)S_G(v)).\)
			\item \( NRL_2(G)=\sum_{uv\in E(G)}(S_G^2(u)+S_G^2(v)-S_G(u)S_G(v))\).
						\item \(NRL_3(G)=\sum_{uv\in E(G)}(S_G(u)-S_G(v)+S_G(u)S_G(v)).\)
			\item \( NRL_4(G)=\sum_{uv\in E(G)}(|S_G(u)-S_G(v)|S_G(u)S_G(v))\).
		\end{enumerate}
	
	\end{definition}
	
	In addition, another new indices will be introduced based on the above two indices as in the table \ref{tableg1},\ref{tableg2}, \ref{tableg3},\ref{tableg4}.

	\begin{longtable}[c]{| c | l | l |}
		\caption{Rehan- Lanel Nbd indices/exponentials \label{tableg1}}\\

		\hline
		\multicolumn{1}{|c}{\textbf{S.N.}} & 
		\multicolumn{1}{|c}{\textbf{Name of the index}} & \multicolumn{1}{|c|}{\textbf{Formula}} \\\hline
		\endfirsthead
		
		\hline
		\endlastfoot

		1&		\multirow{2}{4cm}{The first hyper Nbd RL  index} & \(HNRL_1(G)=\sum\limits_{uv\in E(G)}\left(  S_G^2(u)+S_G^2(v)+S_G(u)S_G(v)\right)^2 \)\\[5ex]
		
		2&		\multirow{2}{4cm}{The second hyper Nbd RL  index} & \(HNRL_2(G)=\sum\limits_{uv\in E(G)}\left(  S_G^2(u)+S_G^2(v)-S_G(u)S_G(v)\right)^2 \)\\[5ex]
		
		3&	\multirow{2}{4cm}{The first  inverse Nbd RL index} & \(INRL_1(G)=\sum\limits_{uv\in E(G)}\dfrac{1}{\left(  S_G^2(u)+S_G^2(v)+S_G(u)S_G(v)\right)} \)\\[6ex]
		
		4&	\multirow{2}{5cm}{The second  inverse Nbd RL index} & \(INRL_2(G)=\sum\limits_{uv\in E(G)}\dfrac{1}{\left( S_G^2(u)+S_G^2(v)-S_G(u)S_G(v)\right)} \)\\[6ex]
		
		5&	\multirow{2}{5cm}{The first general Nbd RL index} & \(GNRL_1(G)=\sum\limits_{uv\in E(G)}\left(S_G^2(u)+S_G^2(v)+S_G(u)S_G(v)\right)^a \)\\[6ex]
		
		6&	\multirow{2}{5cm}{The second  general Nbd RL index} & \(GNRL_2(G)=\prod\limits_{uv\in E(G)}\left( S_G^2(u)+S_G^2(v)-S_G(u)S_G(v)\right)^a \)\\[6ex]
		7&	\multirow{2}{5cm}{The first Nbd RL  exponential} & \(NRL_1(G,x)=\sum\limits_{uv\in E(G)}x^{\left(S_G^2(u)+S_G^2(v)+S_G(u)S_G(v)\right)} \)\\[3ex]
		
		8&	\multirow{2}{5cm}{The second Nbd RL exponential} & \(NRL_2(G,x)=\sum\limits_{uv\in E(G)}x^{\left( S_G^2(u)+S_G^2(v)-S_G(u)S_G(v)\right)} \)\\[3ex]
		
		9&	\multirow{2}{5cm}{The first  hyper Nbd RL  exponential} & \(HNRL_1(G,x)=\sum\limits_{uv\in E(G)}x^{\left(  S_G^2(u)+S_G^2(v)+S_G(u)S_G(v)\right)^2 }\)\\[5ex]
		
		10&	\multirow{2}{5cm}{The second  hyper Nbd RL  exponential} & \(HNRL_2(G,x)=\sum\limits_{uv\in E(G)}x^{\left( S_G^2(u)+S_G^2(v)-S_G(u)S_G(v)\right)^2} \)\\[3ex]
		11&	\multirow{2}{5cm}{The first  inverse Nbd RL exponential} & \(INRL_1(G,x )=\sum\limits_{uv\in E(G)}x^{\dfrac{1}{\left( S_G^2(u)+S_G^2(v)+S_G(u)S_G(v)\right)}} \)\\[3ex]
		12&	\multirow{2}{5cm}{The second inverse Nbd RL exponential} & \(INRL_2(G,x)=\sum\limits_{uv\in E(G)}x^{\dfrac{1}{\left(  S_G^2(u)+S_G^2(v)-S_G(u)S_G(v)\right)} }\)\\[5ex]
		
		13&	\multirow{2}{5cm}{The first general Nbd RL  exponential} & \(GNRL_1(G,x)=\sum\limits_{uv\in E(G)}x^{\left(  S_G^2(u)+S_G^2(v)+S_G(u)S_G(v)\right)^a} \)\\[5ex]
		14&	\multirow{2}{5cm}{The second  general Nbd RL  exponential} & \(GNRL_2(G,x)=\sum\limits_{uv\in E(G)}x^{\left(  S_G^2(u)+S_G^2(v)-S_G(u)S_G(v)\right)^a} \)\\[6ex]

	\end{longtable}
	
	\begin{longtable}[c]{| c | l | l |}
		\caption{Multiplicative Rehan- Lanel Nbd indices/exponentials \label{tableg2}}\\

		\hline
		\multicolumn{1}{|c}{\textbf{S.N.}} & 
		\multicolumn{1}{|c}{\textbf{Name of the index}} & \multicolumn{1}{|c|}{\textbf{Formula}} \\\hline
		\endfirsthead
		
		\hline
		\endlastfoot

		1&	\multirow{2}{5cm}{The first multiplicative Nbd RL index }& \(MNRL_1(G)=\prod\limits_{uv\in E(G)}\left(S_G^2(u)+S_G^2(v)+S_G(u)S_G(v)\right) \)\\[3ex]
		2&	\multirow{2}{5cm}{The second multiplicative Nbd RL index }& \(MNRL_2(G)=\prod\limits_{uv\in E(G)}\left( S_G^2(u)+S_G^2(v)-S_G(u)S_G(v)\right) \)\\[3ex]
		
		3&		\multirow{2}{4cm}{The first  multiplicative hyper Nbd RL index} & \(MHNRL_1(G)=\prod\limits_{uv\in E(G)}\left(  S_G^2(u)+S_G^2(v)+S_G(u)S_G(v)\right)^2 \)\\[5ex]
		
		4&		\multirow{2}{4cm}{The second multiplicative hyper  Nbd RL index} & \(MHNRL_2(G)=\prod\limits_{uv\in E(G)}\left(  S_G^2(u)+S_G^2(v)-S_G(u)S_G(v)\right)^2 \)\\[5ex]
		
		5&	\multirow{2}{4cm}{The first multiplicative inverse Nbd RL index} & \(MINRL_1(G)=\prod\limits_{uv\in E(G)}\dfrac{1}{\left(  S_G^2(u)+S_G^2(v)+S_G(u)S_G(v)\right)} \)\\[6ex]
		
		6&	\multirow{2}{5cm}{The second multiplicative inverse Nbd RL index} & \(MINRL_2(G)=\prod\limits_{uv\in E(G)}\dfrac{1}{\left(  S_G^2(u)+S_G^2(v)-S_G(u)S_G(v)\right)} \)\\[6ex]
		
		7&	\multirow{2}{5cm}{The first general multiplicative Nbd  RL index} & \(GMNRL_1(G)=\prod\limits_{uv\in E(G)}\left( S_G^2(u)+S_G^2(v)+S_G(u)S_G(v)\right)^a \)\\[6ex]
		
		8&	\multirow{2}{5cm}{The second  general multiplicative Nbd RL  index} & \(GMNRL_2(G)=\prod\limits_{uv\in E(G)}\left( S_G^2(u)+S_G^2(v)-S_G(u)S_G(v)\right)^a \)\\[6ex]
		
		9&	\multirow{2}{5cm}{The first multiplicative Nbd  RL  exponential} & \(MNRL_1(G,x)=\prod\limits_{uv\in E(G)}x^{\left(S_G^2(u)+S_G^2(v)+S_G(u)S_G(v)\right)} \)\\[3ex]
		
		10&	\multirow{2}{5cm}{The second multiplicative Nbd  RL exponential} & \(MNRL_2(G,x)=\prod\limits_{uv\in E(G)}x^{\left(  S_G^2(u)+S_G^2(v)-S_G(u)S_G(v)\right)} \)\\[3ex]
		
		11&	\multirow{2}{5cm}{The first multiplicative hyper Nbd RL exponential} & \(MHNRL_1(G,x)=\prod\limits_{uv\in E(G)}x^{\left(  S_G^2(u)+S_G^2(v)+S_G(u)S_G(v)\right)^2 }\)\\[5ex]
		
		12&	\multirow{2}{5cm}{The second multiplicative hyper  Nbd RL exponential} & \(MHNRL_2(G,x)=\prod\limits_{uv\in E(G)}x^{\left( S_G^2(u)+S_G^2(v)-S_G(u)S_G(v)\right)^2} \)\\[3ex]
		13&	\multirow{2}{5cm}{The first multiplicative inverse Nbd RL  exponential} & \(IBRL_1(G,x )=\prod\limits_{uv\in E(G)}x^{\dfrac{1}{\left(  S_G^2(u)+S_G^2(v)+S_G(u)S_G(v)\right)}} \)\\[3ex]
		14&	\multirow{2}{5cm}{The second multiplicative inverse  Nbd RL  exponential} & \(MINRL_2(G,x)=\prod\limits_{uv\in E(G)}x^{\dfrac{1}{\left( S_G^2(u)+S_G^2(v)-S_G(u)S_G(v)\right)} }\)\\[5ex]
		15&	\multirow{2}{5cm}{The first multiplicative general Nbd RL exponential} & \(MGNRL_1(G,x)=\prod\limits_{uv\in E(G)}x^{\left( S_G^2(u)+S_G^2(v)+S_G(u)S_G(v)\right)^a} \)\\[5ex]

		16&	\multirow{2}{5cm}{The second multiplicative general Nbd RL exponential} & \(MGNRL_2(G,x)=\prod\limits_{uv\in E(G)}x^{\left( S_G^2(u)+S_G^2(v)-S_G(u)S_G(v)\right)^a} \)\\[6ex]

	\end{longtable}
	\begin{longtable}[c]{| c | l | l |}
		\caption{Rehan- Lanel Nbd indices/exponentials \label{tableg3}}\\

		\hline
		\multicolumn{1}{|c}{\textbf{S.N.}} & 
		\multicolumn{1}{|c}{\textbf{Name of the index}} & \multicolumn{1}{|c|}{\textbf{Formula}} \\\hline
		\endfirsthead
		
		\hline
		\endlastfoot

		1&		\multirow{2}{4cm}{The third hyper Nbd RL  index} & \(HNRL_1(G)=\sum\limits_{uv\in E(G)}\left(  S_G(u)-S_G(v)+S_G(u)S_G(v)\right)^2 \)\\[5ex]
		
		2&		\multirow{2}{4cm}{The fourth hyper Nbd RL  index} & \(HNRL_2(G)=\sum\limits_{uv\in E(G)}\left(  \left|S_G(u)-S_G(v)\right|S_G(u)S_G(v)\right)^2 \)\\[5ex]
		
		3&	\multirow{2}{4cm}{The third  inverse Nbd RL index} & \(INRL_1(G)=\sum\limits_{uv\in E(G)}\dfrac{1}{\left(  S_G(u)-S_G(v)+S_G(u)S_G(v)\right)} \)\\[6ex]
		
		4&	\multirow{2}{5cm}{The fourth  inverse Nbd RL index} & \(INRL_2(G)=\sum\limits_{uv\in E(G)}\dfrac{1}{\left( \left|S_G(u)-S_G(v)\right|S_G(u)S_G(v)\right)} \)\\[6ex]
		
		5&	\multirow{2}{5cm}{The third general Nbd RL index} & \(GNRL_1(G)=\sum\limits_{uv\in E(G)}\left(S_G(u)-S_G(v)+S_G(u)S_G(v)\right)^a \)\\[6ex]
		
		6&	\multirow{2}{5cm}{The fourth  general Nbd RL index} & \(GNRL_2(G)=\prod\limits_{uv\in E(G)}\left( \left|S_G(u)-S_G(v)\right|S_G(u)S_G(v)\right)^a \)\\[6ex]
		7&	\multirow{2}{5cm}{The third Nbd RL  exponential} & \(NRL_1(G,x)=\sum\limits_{uv\in E(G)}x^{\left(S_G(u)-S_G(v)+S_G(u)S_G(v)\right)} \)\\[3ex]
		
		8&	\multirow{2}{5cm}{The fourth Nbd RL exponential} & \(NRL_2(G,x)=\sum\limits_{uv\in E(G)}x^{\left(\left|S_G(u)-S_G(v)\right|S_G(u)S_G(v)\right)} \)\\[3ex]
		
		9&	\multirow{2}{5cm}{The third  hyper Nbd RL  exponential} & \(HNRL_1(G,x)=\sum\limits_{uv\in E(G)}x^{\left(  S_G(u)-S_G(v)+S_G(u)S_G(v)\right)^2 }\)\\[5ex]
		
		10&	\multirow{2}{5cm}{The fourth  hyper Nbd RL  exponential} & \(HNRL_2(G,x)=\sum\limits_{uv\in E(G)}x^{\left( \left|S_G(u)-S_G(v)\right|S_G(u)S_G(v)\right)^2} \)\\[3ex]
		
		11&	\multirow{2}{5cm}{The third  inverse Nbd RL exponential} & \(INRL_1(G,x )=\sum\limits_{uv\in E(G)}x^{\dfrac{1}{\left( S_G(u)-S_G(v)+S_G(u)S_G(v)\right)}} \)\\[3ex]
		12&	\multirow{2}{5cm}{The fourth inverse Nbd RL exponential} & \(INRL_2(G,x)=\sum\limits_{uv\in E(G)}x^{\dfrac{1}{\left(  \left|S_G(u)-S_G(v)\right|S_G(u)S_G(v)\right)} }\)\\[5ex]
		
		13&	\multirow{2}{5cm}{The third general Nbd RL  exponential} & \(GNRL_1(G,x)=\sum\limits_{uv\in E(G)}x^{\left(  S_G(u)-S_G(v)+S_G(u)S_G(v)\right)^a} \)\\[5ex]
		14&	\multirow{2}{5cm}{The fourth  general Nbd RL  exponential} & \(GNRL_2(G,x)=\sum\limits_{uv\in E(G)}x^{\left(  \left|S_G(u)-S_G(v)\right|S_G(u)S_G(v)\right)^a} \)\\[6ex]

	\end{longtable}
	
	\begin{longtable}[c]{| c | l | l |}
		\caption{Multiplicative Rehan- Lanel Nbd indices/exponentials \label{tableg4}}\\

		\hline
		\multicolumn{1}{|c}{\textbf{S.N.}} & 
		\multicolumn{1}{|c}{\textbf{Name of the index}} & \multicolumn{1}{|c|}{\textbf{Formula}} \\\hline
		\endfirsthead
		
		\hline
		\endlastfoot

		1&	\multirow{2}{5cm}{The third multiplicative Nbd RL index }& \(MNRL_1(G)=\prod\limits_{uv\in E(G)}\left(S_G(u)-S_G(v)+S_G(u)S_G(v)\right) \)\\[3ex]
		\hline 	
		2&	\multirow{2}{5cm}{The fourth multiplicative Nbd RL index }& \(MNRL_2(G)=\prod\limits_{uv\in E(G)}\left( \left|S_G(u)-S_G(v)\right|S_G(u)S_G(v)\right) \)\\[3ex]
		
		3&		\multirow{2}{4cm}{The third  multiplicative hyper Nbd RL index} & \(MHNRL_1(G)=\prod\limits_{uv\in E(G)}\left(  S_G(u)-S_G(v)+S_G(u)S_G(v)\right)^2 \)\\[5ex]
		
		4&		\multirow{2}{4cm}{The fourth multiplicative hyper  Nbd RL index} & \(MHNRL_2(G)=\prod\limits_{uv\in E(G)}\left(  \left|S_G(u)-S_G(v)\right|S_G(u)S_G(v)\right)^2 \)\\[5ex]
		
		5&	\multirow{2}{4cm}{The third multiplicative inverse Nbd RL index} & \(MINRL_1(G)=\prod\limits_{uv\in E(G)}\dfrac{1}{\left(  S_G(u)-S_G(v)+S_G(u)S_G(v)\right)} \)\\[6ex]
		
		6&	\multirow{2}{5cm}{The fourth multiplicative inverse Nbd RL index} & \(MINRL_2(G)=\prod\limits_{uv\in E(G)}\dfrac{1}{\left(  \left|S_G(u)-S_G(v)\right|S_G(u)S_G(v)\right)} \)\\[6ex]
		
		7&	\multirow{2}{5cm}{The third general multiplicative Nbd  RL index} & \(GMNRL_1(G)=\prod\limits_{uv\in E(G)}\left( S_G(u)-S_G(v)+S_G(u)S_G(v)\right)^a \)\\[6ex]
		
		8&	\multirow{2}{5cm}{The fourth  general multiplicative Nbd RL  index} & \(GMNRL_2(G)=\prod\limits_{uv\in E(G)}\left(\left|S_G(u)-S_G(v)\right|S_G(u)S_G(v)\right)^a \)\\[6ex]
		
		9&	\multirow{2}{5cm}{The third multiplicative Nbd  RL  exponential} & \(MNRL_1(G,x)=\prod\limits_{uv\in E(G)}x^{\left(S_G(u)-S_G(v)+S_G(u)S_G(v)\right)} \)\\[3ex]
		
		10&	\multirow{2}{5cm}{The fourth multiplicative Nbd  RL exponential} & \(MNRL_2(G,x)=\prod\limits_{uv\in E(G)}x^{\left( \left|S_G(u)-S_G(v)\right|S_G(u)S_G(v)\right)} \)\\[3ex]
		
		11&	\multirow{2}{5cm}{The third multiplicative hyper Nbd RL exponential} & \(MHNRL_1(G,x)=\prod\limits_{uv\in E(G)}x^{\left(  S_G(u)-S_G(v)+S_G(u)S_G(v)\right)^2 }\)\\[5ex]
		
		12&	\multirow{2}{5cm}{The fourth multiplicative hyper  Nbd RL exponential} & \(MHNRL_2(G,x)=\prod\limits_{uv\in E(G)}x^{\left( \left|S_G(u)-S_G(v)\right|S_G(u)S_G(v)\right)^2} \)\\[3ex]
		
		13&	\multirow{2}{5cm}{The third multiplicative inverse Nbd RL  exponential} & \(IBRL_1(G,x )=\prod\limits_{uv\in E(G)}x^{\dfrac{1}{\left(  S_G(u)-S_G(v)+S_G(u)S_G(v)\right)}} \)\\[3ex]
		14&	\multirow{2}{5cm}{The fourth multiplicative inverse  Nbd RL  exponential} & \(MINRL_2(G,x)=\prod\limits_{uv\in E(G)}x^{\dfrac{1}{\left( \left|S_G(u)-S_G(v)\right|S_G(u)S_G(v)\right)} }\)\\[5ex]
		15&	\multirow{2}{5cm}{The third multiplicative general Nbd RL exponential} & \(MGNRL_1(G,x)=\prod\limits_{uv\in E(G)}x^{\left( S_G(u)-S_G(v)+S_G(u)S_G(v)\right)^a} \)\\[5ex]

		16&	\multirow{2}{5cm}{The fourth multiplicative general Nbd RL exponential} & \(MGNRL_2(G,x)=\prod\limits_{uv\in E(G)}x^{\left( \left|S_G(u)-S_G(v)\right|S_G(u)S_G(v)\right)^a} \)\\[6ex]

	\end{longtable}
	We calculate the first Neighborhood Rehan-Lanel index \(NRL_1\) and the second Rehan-Lanel index \(NRL_2\) for few standard graphs such as \(r-\) regular graph, cycle, complete graph, a path and complete biparttite graph. Additionally, we compute \(NRL_1\) and \(NRL_2\) for the wheel graph, sunflower graph. 
	
	\begin{proposition}
		If \(G\) is an $r$-regular graph with \(n\) vertices and
		\(r \geq 2\), then,
		\begin{enumerate}[label=\roman*.]
			\item \(NRL_1(G)=\dfrac{3nr^3(n-1)^2}{2}\)
			\item \( NRL_2(G)=\dfrac{nr^3(n-1)^2}{2}\)
			\item \(NRL_1(G)=\dfrac{nr^3(n-1)^2}{2}\)
			\item \( NRL_2(G)=0.\)
		\end{enumerate}
		 
	\end{proposition}
	
	\begin{proof}
		Let \(G\) is an $r-$regular graph with \(n\) vertices and \(r \geq 2\).\\

		Then, \(G\) has \(\dfrac{nr}{2}\)	edges. Also \(S_G(u)=r(n-1)\).\\
		
		From definition, we get
		$$ \begin{array}{ll}
			RL_1(G)&=\displaystyle\sum\limits_{uv\in E(G)} (\left( r(n-1)\right)^2+\left( r(n-1)\right)^2+r(n-1)\cdot r(n-1))\\
			&=\dfrac{nr}{2} (3r^2(n-1)^2)\\
			&=\dfrac{3nr^3(n-1)^2}{2}.
		\end{array}$$
		
		$$ \begin{array}{ll}
			RL_2(G)&=\displaystyle\sum\limits_{uv\in E(G)} (\left( r(n-1)\right)^2+\left( r(n-1)\right)^2-r(n-1)\cdot r(n-1))\\
			&=\dfrac{nr}{2} (r^2(n-1)^2)\\
			&=\dfrac{nr^3(n-1)^2}{2}.
		\end{array}$$
		
			$$ \begin{array}{ll}
			RL_3(G)&=\displaystyle\sum\limits_{uv\in E(G)} (\left( r(n-1)\right)-\left( r(n-1)\right)+r(n-1)\cdot r(n-1))\\
			&=\dfrac{nr}{2} (r^2(n-1)^2)\\
			&=\dfrac{nr^3(n-1)^2}{2}.
		\end{array}$$
		
		$$ \begin{array}{ll}
			RL_4(G)&=\displaystyle\sum\limits_{uv\in E(G)} (|\left( r(n-1)\right)-\left( r(n-1)\right)|\cdot r(n-1)\cdot r(n-1))\\
			&=0.
		\end{array}$$
	\end{proof}
	
	\begin{corollary}
		Let	\(C_n\)
		be a cycle with \(n \geq 3\) vertices. Then,
		\begin{enumerate}[label=\roman*.]
			\item \(NRL_1(C_n)=12n(n-1)^2\).
			\item \(NRL_2(C_n)=4n(n-1)^2.\)
			  \item \(NRL_3(C_n)=4n(n-1)^2\).
			\item \(NRL_4(C_n)=0.\)
		\end{enumerate}

	\end{corollary}
	\begin{corollary}
		Let \(K_n\) be a complete graph with \(n\geq 3\)
		vertices. Then,
		\begin{enumerate}[label=\roman*.]
			\item \(NRL_1(K_n)=\dfrac{3n(n-1)^5}{2}\).
			\item \(NRL_2(K_n)=\dfrac{n(n-1)^5}{2}.\)
			 \item \(NRL_3(K_n)=\dfrac{4n(n-1)^5}{2}\).
			\item \(NRL_4(K_n)=0.\)
		\end{enumerate}

	\end{corollary}
	
	\begin{proposition}
		Let \(P_n\) be a path with \(n\geq 3\) vertices. Then,
		\begin{enumerate}[label=\roman*.]
			\item \(NRL_1(P_n)= 48n-106.\) 
			\item \(NRL_2(P_n)=16n-34\).
			  \item \(NRL_3(P_n)= 48n-106.\) 
			\item \(NRL_4(P_n)=16n-34\).
		\end{enumerate}
		
	\end{proposition}
	
	\begin{proof}
		The proof can be done similarly as above.
	\end{proof}
	
	\begin{proposition}
		Let $K_{m,n}$ be a complete bipartite graph with
		\(1\leq m\leq n\) and \(n \geq 2\). Then, 
		\begin{enumerate}[label=\roman*.]
			\item \(NRL_1(K_{m,n})= 3(mn)^3.\) 
			\item \(NRL_2(K_{m,n})= (mn)^3.\)
			   \item \(NRL_3(K_{m,n})= m^2n^2.\) 
			\item \(NRL_4(K_{m,n})= 0.\)
		\end{enumerate}	
		
	\end{proposition}
	\begin{proof}
		The proof can be done similarly as above.
	\end{proof}
	
	\begin{theorem}
		Let \(W_n\) be a wheel graph. Then,
		
		\begin{enumerate}[label=\roman*.]
			\item \(NRL_1(W_n)=n(16n^2+ 66n+144)\).
			\item \(NRL_2(W_n)=n(8n^2+ 6n+72)\).
			 \item \(NRL_3(W_n)=n(4n^2+ 28n+42)\).
			\item \(NRL_4(W_n)=6n^2\left| 3-n\right|(n+6)\).
		\end{enumerate}
		
	\end{theorem}
	
	\begin{proof}
		In \(W_n\), There are two types of edges as follows \(d(v)=n\) and \(d(v_i)=3\),
		$$\begin{array}{lll}
			E_1=&\left\lbrace uv \in E(W_n)| \ d(u)=d(v)=3 \right\rbrace, &|E_1|=n   \\
			E_2=&\left\lbrace uv \in E(W_n)|\ d(u)=3 ,d(v)=n \right\rbrace , &|E_2|=n 
		\end{array}$$

		Therefore, the edge partition of the wheel graph for neighborhood RL indices.
		
		$$\begin{array}{lll}
			E_1&=\left\lbrace uv \in E(W_n)| S_G(u)=S_G(v)=n+6 \right\rbrace, &|E_1|=n.  \\
			E_2&=\left\lbrace uv \in E(W_n)| S_G(u)=n+6 ,S_G(v)=3n \right\rbrace , &|E_2|=n.
		\end{array}$$

		$$\begin{array}{ll}
			NRL_1(W_n)&=\sum\limits_{uv\in E(G)}(S_G^2(u)+S_G^2(v)+S_G(u)S_G(v))\\
			&=n\left[ (n+6)^2+(n+6)^2+(n+6)(n+6)\right] +n\left[(n+6)^2+(3n)^2+(n+6)(3n) \right] \\
			&=n\left[ 4(n+6)^2+9n^2+3n(n+6)\right]\\
			&=n(16n^2+ 66n+144).
		\end{array}$$
		
		$$\begin{array}{ll}
			NRL_2(W_n)&=\sum\limits_{uv\in E(G)}(S_G^2(u)+S_G^2(v)-S_G(u)S_G(v))\\
			&=n\left[ (n+6)^2+(n+6)^2-(n+6)(n+6)\right] +n\left[(n+6)^2+(3n)^2-(n+6)(3n) \right] \\
			&=n\left[ 2(n+6)^2+9n^2-3n(n+6)\right]\\
			&=n(8n^2+ 6n+72).
		\end{array}$$
		
		$$\begin{array}{ll}
			NRL_3(W_n)&=\sum\limits_{uv\in E(G)}(S_G(u)-S_G(v)+S_G(u)S_G(v))\\
			&=n\left[ (n+6)-(n+6)+(n+6)(n+6)\right] +n\left[(n+6)-(3n)+(n+6)(3n) \right] \\
			&=n(4n^2+ 28n+42).
		\end{array}$$
		
		$$\begin{array}{ll}
			NRL_4(W_n)&=\sum\limits_{uv\in E(G)}(\left|S_G(u)-S_G(v)\right|S_G(u)S_G(v))\\
			&=n\left|(n+6)-(n+6)\right|(n+6)(n+6) +n\left|(n+6)-(3n)\right|(n+6)(3n) \\
			&=6n^2\left| 3-n\right|(n+6).
		\end{array}$$
	\end{proof}
	
	\begin{theorem}
		Let \(S\!f_n\) be a sunflower graph. Then,
		
		\begin{enumerate}[label=\roman*.]
			\item \(NRL_1(S\!f_n)=n(328n^2+406n+504).\)
			\item \(NRL_2(S\!f_n)=2n(4n^2+3n+36).\)
		\end{enumerate}
		
	\end{theorem}
	
	\begin{proof}

		In \({S\!f}_n\), there are five types of nbd edges as follows,
		$$\begin{array}{lll}
			NE_1  & =\left\lbrace uv \in E(Sf_n)| S_G(v)=S_G(v)=3n+10 \right\rbrace, & |E_1|=n. \\
			NE_2   & =\left\lbrace uv \in E(Sf_n)| S_G(v_0)=7n ,S_G(v)=3n+10 \right\rbrace , & |E_2|=n.\\
			NE_3 & =\left\lbrace uv \in E(Sf_n)| S_G(u)=3n+4, d(v)=3n+10 \right\rbrace,& |E_3|=n.\\
			NE_4& =\left\lbrace uv \in E(Sf_n)| S_G(v_0)=7n ,S_G(u)=3n+4 \right\rbrace , & |E_4|=n.\\
			NE_5&=\left\lbrace uv \in E(Sf_n)| S_G(w)=3n ,S_G(v_0)=7n \right\rbrace , & |E_5|=n.
		\end{array}$$	
		
		$$\begin{array}{ll}
			NRL_1(S\!f_n)&= \sum\limits_{uv\in E(G)}(S_G^2(u)+S_G^2(v)+S_G(u)S_G(v))\\
			&=n\left((3n+10)^2+(3n+10)^2+(3n+10)(3n+10) \right)+n\left( (7n)^2+(3n+10)^2+7n\cdot( 3n+10)\right) \\&+n\left((3n+4)^2+(3n+10)^2+(3n+4)(3n+10) \right) +n\left((7n)^2+(3n+4)^2+7n \cdot (3n+4)\right)\\&+n\left( (3n)^2+(7n)^2+3n\cdot 7n\right)  \\
			&=n(328n^2+406n+504).
		\end{array}$$
		
		$$\begin{array}{ll}
			NRL_2(S\!f_n)&= \sum\limits_{uv\in E(G)}(S_G^2(u)+S_G^2(v)-S_G(u)S_G(v))\\
			&=n\left((3n+10)^2+(3n+10)^2-(3n+10)(3n+10) \right)+n\left( (7n)^2+(3n+10)^2-7n\cdot( 3n+10)\right) \\&+n\left((3n+4)^2+(3n+10)^2-(3n+4)(3n+10) \right) +n\left((7n)^2+(3n+4)^2-7n \cdot (3n+4)\right)\\&+n\left( (3n)^2+(7n)^2-3n\cdot 7n\right)  \\
			&=2n(4n^2+3n+36).
		\end{array}$$
		
%
	\end{proof}

	Similarly, we can derive the results for all the other indices and the exponentials listed in the Table \ref{tablef1}, \ref{tablef2}, \ref{tablef3} and \ref{tablef4}.
	
	Additionally, we aspire to introduce several more indices by incorporating \(ve-\) degree, \(HDR\) degree, K- Banhatti degree as a future work.
	
	\subsection{The fifth and sixth Rehan-Lanel indices}
	
		Using the concept of the exponential of a one number of the other, here we define the two more indices as follows.
		
		\begin{definition}
			The fifth and sixth Rehan-Lanel indices are defined by,
			
			\begin{enumerate}[label=\roman*.]
				\item \(RL_5(G)=\sum\limits_{uv\in E(G)}{d_G(u)}^{d_G(v)}\).
				\item \(RL_6(G)=\sum\limits_{uv\in E(G)}{d_G(u)}^{d_G(v)}+{d_G(v)}^{d_G(u)}\).
			\end{enumerate}
		\end{definition}
		\subsection{Closeness Centrality}
	
	Closeness centrality, in the context of a connected graph, refers to a metric that quantifies the centrality of a node inside a network. It is computed by taking the reciprocal of the sum of the shortest path lengths between the node and all other nodes in the graph.
	
	The normalized form of the closeness centrality is defined as,
	\[c(u)=\dfrac{n-1}{\sum\limits_{v\in V(G)}d(v,u)}\]
	where \(d(u,v)\) is the distance between the verices \(u\) and \(v\) and also \(n\) denotes the number of the vertices of the graph \(G.\)
	
	Inspired by indices like Zagreb indices, we define the first closeness centrality index(the seventh Rehan-Lanel index) \(C_1(G)/RL_7(G)\)  and the second closeness centrality index(the eighth Rehan-Lanel index)\(C_2(G)/RL_8(G)\) as follows;
	
	\begin{definition}
		The first closeness centrality index(the seventh Rehan-Lanel index) \(C_1(G)/RL_7(G)\)  and the second closeness centrality index(the eighth Rehan-Lanel index) of a graph \(G\) are defined by,
		\[RL_7(G)=C_1(G)=\sum_{uv\in E(G)}c(u)+c(v)\]
		and 
		\[RL_8(G)=C_2(G)=\sum_{uv\in E(G)}c(u)c(v).\]

	\end{definition}

	Also, using the first and second Zagreb indices, Furtula
	and Gutman \cite{article82} introduced forgotten topological index
	(also called F-index) which was defined as
	\[F(G)=\sum_{u\in V(G)}d^3_G(u)=\sum_{uv\in E(G)}d^2_G(u)+d^2_G(v).\]

	Based on the fogotten index and the fifth and sixth Rehan-Lanel indices, we define the fogotten closeness centrality index/the seventh Rehan-Lanel index(\(FC(G)/RL_7(G)\)) of a graph \(G\) as follows:
	
	\begin{definition}
		The  fogotten closeness centrality index/the ninth Rehan-Lanel index(\(FC(G)/RL_9(G)\)) of a graph \(G\) is defined by, 
		\[RL_9(G)=FC(G)=\sum_{uv\in E(G)}c^2(u)+c^2(v).\]
	\end{definition}

	In 2020, Ivan Gutmn \cite{article83} introduced the Sombor index (SO) of a graph \(G\) as follows, 
	
	\[SO(G)=\sum_{uv \in E(G)} \sqrt{d^2_G(u)+d^2_G(v)}.\]
	
	Motivated by Sombor index, we define the Centrality Sombor index/the tenth Rehan-Lanel index(\(CSO(G)/RL_{10}(G)\)) as follows:
	
	\begin{definition}
		The centrality sombor index/the tenth Rehan-Lanel index(\(CSO(G)/RL_8(G)\)) is defined by, 
		\[CSO(G)=RL_{10}(G)=\sum_{uv \in E(G)} \sqrt{c^2_G(u)+c^2_G(v)}.\]
	\end{definition}

	Inspired by work on Sombor indices, VR Kulli\cite{article84}, an indian scholar introduced the Nirmala index of a graph as follows:
	\[N(G)=\sum_{uv \in E(G)} \sqrt{d_G(u)+d_G(v)}.\]

	Inspired by works on nirmala index, we define the eleventh Rehan-Lanel index/ the centrality nirmala index of a graph \(G)\) as follows:
	
	\begin{definition}
		The eleventh Rehan-Lanel index/ the centrality nirmala index of a graph \(G)\) is defined by,
		\[CN(G)=RL_{11}(G)=\sum_{uv \in E(G)} \sqrt{c_G(u)+c_G(v)}.\]
		
	\end{definition}
	
	Based on the Albertson irregularity index and the closeness centrality of a vertex of a graph \(G\), we define the Albertson centrality index/the tenth Rehan-Lanel index as follows:

	\begin{definition}
		The Albertson centrality index/the twelfth Rehan-Lanel index of a connected graph \(G\) is defined by,
		
		\[AC(G)=RL_{12}(G)=\sum_{uv \in E(G)} |c_G(u)-c_G(v)|.\]
	\end{definition}

	In addition, we introduce a new degree definition of of a vertex of  a graph \(G\) based on the absolute deviations of  the degree of a vertex from the its adjacent vertices. We name this new degree as CL degree of a vertex of  a graph \(G.\)

	\begin{definition}
		The Chandana-Lanel degree (CL degree) of a vertex \(u\) of a graph \(G\) is defined by, \[CL(u)=\max\left\lbrace \left|d_G(u)-d_G(v) \right|: v \text{ is a vertex adjacent to   } u \right\rbrace \]
	\end{definition}

	Using this definition and the Zagreb  indices, we introduce the thirteenth and fourteenth Rehan-Lanel indices as follows,
	
	\begin{definition}
		The thirteenth and fourteenth Rehan-Lanel indices of a graph \(G\) are defined by,
		
		\[RL_{13}(G)=\sum_{uv \in E(G)} CL_1(u)+CL_2(v) \ \text{and} \ RL_{14}(G)=\sum_{uv \in E(G)} CL_1(u)CL_2(v). \]
	\end{definition}  
	
	\begin{definition}
		The Forgotten Rehan-Lanel index(fifteenth RL) of a graph \(G\) is defined by,
		
		\[FRL(G)=RL_{15}(G)=\sum_{uv \in E(G)} CL^2_1(u)+CL^2_2(v). \]
	\end{definition}

	\begin{definition}
	The Sombor CL of a graph \(G\) is defined by,
	
	\[SCL(G)=RL_{16}(G)=\sum_{uv \in E(G)} \sqrt{CL^2_1(u)+CL^2_2(v)}. \]
\end{definition}  

	\begin{definition}
	The Nirmala CL of a graph \(G\) is defined by,
	
	\[NCL(G)=RL_{17}(G)=\sum_{uv \in E(G)} \sqrt{CL_1(u)+CL_2(v)}. \]
\end{definition} 
The Heronian mean \(H\) of two non-negative real numbers \(a\) and \(b\) is given by, \(H=\dfrac{1}{3}(a+\sqrt{ab}+b).\)

Motivated by definition, here we define the Heronian Rehan-Lanel index of a graph \(G\) as follows,

\begin{definition}
 The Heronian Rehan-Lanel index of a graph \(G\) by,

\[HRL(G)=\sum_{uv \in E(G)}d_G(u)+\sqrt{d_G(u)d_G(v)}+d_G(v).\]
\end{definition}

	\section{Conclusion}

	This study introduces four new indices, namely the first, second, third and fourth Rehan-Lanel index. By utilizing these indices, we present an additional 448 indices/exponentials by including Banhatti degree, Revan degree, Temperature, Domination degree, KV indices, and Neighborhood degree. Afterwards, we conducted comprehensive calculations for several indices, on typical graph types such as \(r-\) regular graph, complete graph, cycle, path, and complete bipartite graph.Additionally, we extend our analysis to the wheel graph and sunflower graph, excluding \(DRL_1, DRL_2, DRL_3, DRL_4\), for which results are presented exclusively for the French windmill graph. Furthermore, we are currently exploring the chemical applications of the aforementioned indices/exponentials. In addition, based of the exponential of the degree of a vertex and closeness centrality concept, we defined another 8 indices. Furthermore, we introduced a new degree of a vertex and named as Chandana-Lanel degree of a vertex(CL degree). Using this degree we defined 5 more indices. Also, using heronian mean of two numbers, we defined the Heronian Rehan-Lanel index.   As part of our future work, we aim to introduce distance-based indices associated with the 462 indices discussed in this study. 
	
	\bibliography{myref}

\end{document}